\documentclass{amsbook}

\newtheorem{prop}{Proposition}[chapter]
\newtheorem{theo}[prop]{Theorem}
\newtheorem{lemm}[prop]{Lemma}
\newtheorem{coro}[prop]{Corollary}
\newtheorem{prob}[prop]{Problem}

\newtheorem{conj}[prop]{Conjecture}

\theoremstyle{definition}
\newtheorem{defi}[prop]{Definition}
\newtheorem{exam}[prop]{Example}

\theoremstyle{remark}
\newtheorem{rema}[prop]{Remark}

\numberwithin{section}{chapter}
\numberwithin{equation}{chapter}

\usepackage{amsmath, amscd, amssymb}
\usepackage[frame,cmtip,arrow,matrix,line,graph,curve]{xy}
\usepackage{graphpap, color}
\usepackage[mathscr]{eucal}
\usepackage{cancel}
\usepackage{verbatim}

\def\sO{{\mathscr O}}

\def\sE{\mathscr{E}}
\def\sF{\mathscr{F}}

\newcommand{\CC}{\mathbb{C}}

\newcommand{\PP}{\mathbb{P}}
\newcommand{\QQ}{\mathbb{Q}}
\newcommand{\RR}{\mathbb{R}}
\newcommand{\ZZ}{\mathbb{Z}}

\newcommand{\bL}{\mathbf{L}}

\newcommand{\bm}{\mathbf{m}}

\newcommand{\cal}{\mathcal}

\def\cA{{\cal A}}
\def\cB{{\cal B}}

\def\cD{{\cal D}}
\def\cE{{\cal E}}
\def\cF{{\cal F}}
\def\cG{{\cal G}}
\def\cH{{\cal H}}

\def\cL{{\cal L}}
\def\cM{{\cal M}}

\def\cO{{\cal O}}

\def\cQ{{\cal Q}}

\def\cT{{\cal T}}
\def\cU{{\cal U}}
\def\cV{{\cal V}}
\def\cW{{\cal W}}
\def\cW{{\cal W}}

\def\cI{{\cal I}}

\def\cX{\mathcal{X} }

\def\fC{\mathfrak{C}}

\def\fV{\mathfrak{V}}

\def\fm{\mathfrak{m}}

\newcommand{\tz}{\tilde{z}}

\def\and{\quad{\rm and}\quad}

\def\mapright#1{\,\smash{\mathop{\lra}\limits^{#1}}\,}
\def\mapleft#1{\,\smash{\mathop{\longleftarrow}\limits^{#1}}\,}
\DeclareMathOperator{\pr}{pr} 
 
\DeclareMathOperator{\End}{End}  \DeclareMathOperator{\Hom}{Hom}
\DeclareMathOperator{\diag}{diag} 
\DeclareMathOperator{\image}{Im} 
\DeclareMathOperator{\id}{id} 
 
\DeclareMathOperator{\tr}{tr} 

\DeclareMathOperator{\rank}{rank}
\DeclareMathOperator{\spec}{Spec}

\def\dbar{{\overline{\partial}}}

\def\rat{\mathrm{rat}}

\def\lalpbe{_{\alpha\beta}}
\def\bul{^\bullet}

\def\lbe{_\beta}
\def\OadE{\Omega^{0,1}(\End E)_k}

\def\Lap{\bigtriangleup}

\def\vsp{\vskip5pt}

\def\PP{\mathbb{P}}
\def\CC{\mathbb{C}}

\def\lra{\longrightarrow}
\def\mapright#1{\,\smash{\mathop{\lra}\limits^{#1}}\,}
\def\cO{\mathcal{O}}

\def\Spec{\mathrm{Spec}}

\def\coker{\mathrm{coker} }
\def\ker{\mathrm{ker} }

\def\lsta{_\ast}

\def\ti {\tilde }
\def\ud{^{\bullet} }

\def\bbH{\mathbb{H} }
\def\stax{^{\ast} }
\def\sta{^{\ast}}
\def\sub{\subset}

\def\fa{\mathfrak a}
\def\fb{\mathfrak b}
\def\defeq{\!:=}
\def\bl{\bigl(}
\def\br{\bigr)}
\def\dual{^\vee}

\def\upmo{^{-1}}

\let\eps=\epsilon

\def\bR{\mathbf R}

\def\veps{\varepsilon}
\def\fc{\mathfrak c}
\def\lalp{_\alpha}
\def\bM{\mathbf M}
\def\udag{^\dag}
\def\dpri{^{\prime\prime}}

\def\fa{\mathfrak{a} }
\def\fX{{X} }
\def\bpar{\bar{\partial} }

\def\black{\color{black}}

\def\beq{\begin{equation}}
\def\eeq{\end{equation}}
\def\Pii{{\pi}}

\def\ured{^{\mathrm{red}} }

\begin{document}
\frontmatter
\title{Categorification of Donaldson-Thomas invariants via perverse sheaves}

\author{Young-Hoon Kiem}
\address{Department of Mathematics and Research Institute of Mathematics, Seoul National University, Seoul 151-747, Korea}
\email{kiem@snu.ac.kr}
\thanks{The first named author was supported in part by Korea NRF Grant 2011-0027969.}

\author{Jun Li}
\address{%Shanghai Institute for Mathematical Sciences, Fudan University, China; 
Department of Mathematics,  Stanford University, CA 94305, USA}
\email{jli@math.stanford.edu}
\thanks{The second named author was partially supported by NSF grant DMS-1104553 and DMS-1159156.}

\date{November, 2015}

\begin{abstract}
We show that there is a perverse sheaf on a fine moduli space of stable sheaves on a smooth projective Calabi-Yau 3-fold, 
which is locally the perverse sheaf of vanishing cycles for a local Chern-Simons functional. This perverse sheaf lifts to a mixed Hodge module and gives us a cohomology theory which enables us to define the Gopakumar-Vafa invariants mathematically. This paper is a completely rewritten version of the authors' 2012 preprint, arXiv:1212.6444.
\end{abstract}

\maketitle

\setcounter{page}{4}
\tableofcontents

%-----------------------------------------------------------------------------
% Beginning of preface.tex
%-----------------------------------------------------------------------------
%
% AMS-LaTeX 1.2 sample file for a monograph, based on amsbook.cls.
% This is a data file input by chapter.tex.
%%%%%%%%%%%%%%%%%%%%%%%%%%%%%%%%%%%%%%%%%%%%%%%%%%%%%%%%%%%%%%%%%%%%%%%%

\chapter*{Preface}

The Donaldson-Thomas invariant (DT invariant, for short) is a virtual count of stable sheaves on a smooth projective Calabi-Yau 3-fold $Y$ over $\CC$ which was defined as the degree of the virtual fundamental class of the moduli space $\fX$ of stable sheaves  on $Y$ (\cite{Tho}). Using microlocal analysis, Behrend showed that the DT invariant is in fact the Euler number of the moduli space, weighted by a constructible function $\nu_{\fX}$, called the Behrend function (\cite{Beh}). Since the ordinary Euler number is the alternating sum of Betti numbers of cohomology groups, it is reasonable to ask if the DT invariant is in fact the Euler number of a cohomology theory on $\fX$. 

On the other hand, it is known that the moduli space is locally the critical locus of a holomorphic function, called a local Chern-Simons functional (\cite{JoSo}). Given a holomorphic function $f$ on a complex manifold $V$, one has the perverse sheaf 
$\phi_f(\QQ[\dim V-1])$ of vanishing cycles supported on the critical locus. So the moduli space is covered
by charts each of which comes with the perverse sheaf of vanishing cycles $\phi_f(\QQ[\dim V-1])$.

% and the Euler number of this perverse sheaf at a point $x$ equals $\nu_{\fX}(x)$. 

This motivated Joyce and Song to raise the following question (\cite[Question 5.7]{JoSo}).

\medskip 
\noindent {\it Let $\fX$ be a moduli space of simple coherent sheaves on $Y$. Does there exist a natural perverse sheaf $P$ on the underlying analytic space which is locally isomorphic to the sheaf $\phi_f(\QQ[\dim V-1])$ of vanishing cycles for $(V,f)$ above?}

\medskip
The purpose of this paper is to provide an affirmative answer.  
\begin{theo}\label{thmInt} %(Theorem \ref{vm16})\\
Let $\fX$ be a moduli space of simple sheaves on a smooth projective Calabi-Yau 3-fold $Y$. Suppose the reduced scheme $X\ured$ of $X$ is of finite type and admits a tautological family. Then the collection of perverse sheaves
$\phi_{f}(\QQ[\dim V-1])$ geometrically glue to a perverse sheaf $P$ on $X$. The same holds for polarizable mixed Hodge modules.
%an \'etale Galois cover  $$\rho:X^\dagger\to X=X^\dagger/G$$ and a perverse sheaf $P\ud$ on $X^\dagger$, which is locally isomorphic to the perverse sheaf  $\phi_f(\QQ[\dim V-1])$ of vanishing cycles for the local Chern-Simons functional $f$. In fact, for any \'etale Galois cover  $\rho:X^\dagger\to X$, there exists such a perverse sheaf $P\ud$ if and only if the line bundle $\rho^*\det( \Ext_{\pi}\bul(\cE,\cE))$ admits a square root on $X^\dagger$ where $\pi:X\times Y\to X$ is the projection and $ \Ext_{\pi}\bul(\cE,\cE)=R\pi_*R\cH om(\cE,\cE)$.
%(2) If $\fX$ parmaeterizes one dimensional sheaves, there are an \'etale Galois cover $\rho:X^\dagger\to X$ and a perverse sheaf $P\ud$ on $X^\dagger$ which is locally the perverse sheaf of vanishing cycles of a CS functional.
%(2) Let $M=[\cX/G]$ be the moduli stack of Gieseker semistable sheaves on $Y$, which we write as the quotient of the semistable part $\cX$ of the Quot scheme by a reductive group $G$. Then there is a $G$-equivariant perverse sheaf $P\udG\in Perv^G(X)$ whose restriction to the Joyce-Song slice is isomorphic to the equivariant perverse sheaf of vanishing cycles for the local Chern-Simons functional $f$ with degrees shifted by $\dim G/H$ where $H$ is the stabilizer at the center of the slice.
\end{theo}
%We will also prove the same for mixed Hodge modules (Theorem \ref{thmMHM}). Namely, there is a mixed Hodge module $M\ud$ on $X$ whose underlying perverse sheaf is $\rat(M\ud)=P\ud$. (See \S\ref{sec9}.)

%The perverse sheaf $P$ is unique up to tensoring by a $\ZZ_2$-local system on $X$. 

\medskip

%Suppose $\det\Ext_{\pi}\bul(\cE,\cE)$ has a square root after pulling back to an \'etale Galois $G$-cover $\rho:X^\dagger\to X$ as in the case of one dimensional sheaves. 
As an application of Theorem \ref{thmInt}, we use the hypercohomology $\bbH^i(X,P)$ of $P$ to deduce 
the \emph{DT (Laurent) polynomial}
$$DT_t^Y(X)=\sum_i t^i\dim \bbH^i_c(X,P)
$$
such that $DT_{-1}^Y(X)$ is the ordinary DT invariant by \cite{Beh}. 
%Here the shift $[-\dim M^s(P,L)]$ comes from the degree convention in the theory of perverse sheaves which puts the middle cohomology at degree $0$. 
%In general, if we write the moduli stack of semistable sheaves on $Y$ as the quotient $[\cX/G]$,  the equivariant perverse sheaf $A\udG$ gives us a power series
%\[ DT_t^Y(P,L)=\sum_i t^i\dim \bbH^i_G(X; A\udG[-\dim G])\]
%by taking the equivariant cohomology of $A\udG[-\dim G]$, i.e. the hypercohomology of  $\tilde{A}\ud[-\dim G]$ on the homotopy quotient $EG\times_G X$ where $\tilde{A}\ud$ is the image of $A\udG$ by the forgetful functor $D^b_G(X)\to D^b(EG\times_GX)$. Here $X$ denotes the analytic space underlying $\cX$. 

Another application is a mathematical theory of Gopakumar-Vafa invariants (GV invariants, for short) proposed in \cite{GoVa}. Let $\fX$ be a moduli space of stable sheaves supported on curves of homology class $\beta\in H_2(Y,\ZZ)$.
The GV invariants are integers $n_h(\beta)$ for $h\in \ZZ_{\ge 0}$ defined by an $sl_2\times sl_2$ action on \emph{some cohomology} of $\fX$ such that $n_0(\beta)$ is the DT invariant of $\fX$ and that they give all genus Gromov-Witten invariants $N_g(\beta)$ of $Y$ via a relation in \cite{GoVa} (cf. \eqref{7.1}).   
Applying Theorem \ref{thmInt} to $\fX$, %the moduli of stable sheaves of pure dimension one on $Y$ of $\chi=1$, 
%when $\det \Ext_{\pi}\bul(\cE,\cE)$ admits a square root, 
we obtain a perverse sheaf $P$ on $X$ which is locally the perverse sheaf of vanishing cycles. 
By the relative hard Lefschetz theorem for the morphism to the Chow scheme (\cite{Sai88}), we have an action of $sl_2\times sl_2$ on $\bbH^*(X,\hat{P})$ where $\hat{P}$ is the gradation of $P$ by the weight filtration of $P$. This gives us a geometric theory of GV invariants which we conjecture to give all the Gromov-Witten invariants $N_g(\beta)$.

\bigskip

Our proof of Theorem \ref{thmInt} uses gauge theory. By the Seidel-Thomas twist (\cite[Chapter 8]{JoSo}), it suffices to consider only vector bundles on $Y$. Let $\cB_{si}=\cA_{si}/\cG$ be the space of semiconnections on a hermitian vector bundle $E$ modulo the gauge group action %;  let $\cB_{si}$ be the open subset of simple sheaves,
and let $\fX\subset \cB_{si}$ be a locally closed complex analytic subspace parameterizing integrable semiconnections. 
Let $CS:\cA_{si}\to \CC$ be the holomorphic Chern-Simons functional. We call a finite dimensional complex submanifold $V$ of $\cA$  a \emph{CS chart} if the critical locus of 
$f=CS|_{V}$ is an open complex analytic subspace of $\fX$. 
%By \cite{JoSo}, at each $x\in X=\fX\lred$, we have a CS chart $V$ with $T_xV=T_x\fX$, which we call the Joyce-Song chart (JS chart, for short). Thus we have a perverse sheaf $P\ud|_V$ on $V\cap X$.  
%One of the difficulties in gluing the local perverse sheaves $P\ud|_V$ is that the dimensions of the JS charts $V$ vary from point to point. 

Our proof of Theorem \ref{thmInt} consists of three parts:
\bigskip

\noindent
I. We show that there are \begin{enumerate}
\item an open cover $\fX=\cup_\alpha \fX_\alpha$;
\item $\fX\lalp=(df\lalp=0)$ for an $r$-dimensional CS chart $(V\lalp,f\lalp)$, with $V\lalp\sub \cA$ and
$f\lalp=CS|_{V\lalp}$;
%a (continuous) family $\cV_\alpha\to U_\alpha$ of CS charts of constant dimension $r$, each of which contains the JS chart;
\item for any pair $(\alpha,\beta)$, there are open neighborhoods $V\lalpbe\sub V\lalp$ and $V_{\beta\alpha}\sub V\lbe$ of 
$\fX\lalp\cap\fX\lbe$, and
a biholomorphic map $\varphi\lalpbe:V\lalpbe\to V_{\beta\alpha}$ with $f_\beta\circ \varphi\lalpbe =f_\alpha$.
%$(\cV_{\alpha\beta},f\lalp|_{\cV\lalpbe})$ is equivalent to $(\cV_{\beta\alpha},f\lbe|_{\cV_{\beta\alpha}})$ as CS charts.
\end{enumerate}
%where  and 
%Here the equivalence of CS charts is defined in the obvious manner (cf. Definition \ref{loc-eq}). 
An analytic space $\fX$ equipped with the above (1)-(3) will be called a \emph{critical virtual manifold}. See Definition \ref{1.5} for the precise statement.
%In particular, they give $\fX$ a Chern-Simons-Virtual-Manifold structure (see Definition \ref{vm1}.)  

\bigskip

\noindent II. We prove that given a critical virtual manifold structure, we can extract perverse sheaves $P_\alpha$ on $X_\alpha$ for all $\alpha$ (i.e. the perverse sheaves of vanishing cycles) and gluing isomorphisms 
$\sigma_{\alpha\beta}:P_\alpha|_{X_{\alpha\beta}}\to P_\beta|_{X_{\alpha\beta}}$ of geometric origin
(cf. Definition \ref{2.15}) 
where $X_{\alpha\beta}=X_\alpha\cap X_\beta$. The same hold for polarizable mixed Hodge modules.

\bigskip

\noindent
III. We show that the 2-cocycle obstruction for gluing $\{P_\alpha\}$ to a global perverse sheaf is $\ZZ_2$-valued
and vanishes if there is a square root of the determinant line bundle $\det Rp_*R\cH om(\cE,\cE)|_{X\ured}$, where $\cE$ denotes (local) tautological family while $p:\fX\ured\times Y\to \fX\ured$ is the projection. Therefore the local perverse sheaves $\{P_\alpha\}$ glue to a global perverse sheaf if there is a square root of  $\det Rp_*R\cH om(\cE,\cE)|_{X\ured}$ in $\mathrm{Pic}(X\ured)$ whose existence is guaranteed by an argument of Okounkov when $X\ured$ admits a global tautological family.

%We deal only with the case of global quotient stacks because the theory of perverse sheaves on Artin stacks seems to be incomplete at this stage. As one can find in \cite{BL}, the correct notion of the derived category of a global quotient $[\cX/G]$ is not the derived category of $G$-equivariant sheaves on $\cX$. One has to go through the homotopy quotient $EG\times_GX$ but it is not clear to us how to formulate it for Artin stacks in general. 

%We remark that our proof of item (2) in Theorem \ref{thmInt}  doesn't use geometric invariant theory in any significant way. All we need is that there are locally closed analytic complex analytic subspaces $U_\alpha$ and reductive subgroups $H_\alpha$ of $G$ such that $G\times_{H_\alpha}U_\alpha\to GU_\alpha\subset \cX$ form an \'etale cover and that there are local Chern-Simons functionals $f_\alpha$ whose critical loci are $U_\alpha$ with compatible obstruction assignments. Therefore our proof applies to open substacks of derived category objects as soon as one can check the compatibility of obstruction theories of local Chern-Simons functionals in \cite{JoCS}.

\bigskip 

This paper consists of two parts and 
a detailed layout of each part is provided at the beginning. 

\bigskip

This paper is a completely revised version of the authors' preprint posted on the arXiv in December 2012. All the ideas and arguments are essentially the same. Some related results were independently obtained by Chris Brav, Vittoria Bussi, Delphine Dupont, Dominic Joyce and Balazs Szendroi in \cite{BBDJS}. We are grateful to Dominic Joyce for his comments and suggestions. We thank Martin Olsson and Yan Soibelman for their comments. We also thank Takuro Mochizuki and Masaki Kashiwara for answering questions on mixed Hodge modules.

%\bigskip

%\noindent\textbf{Acknowledgement}. The first named author gratefully acknowledges the support from an FRG grant and Stanford math department during his visit to Stanford. 

%\noindent \textbf{Notations}. A complex analytic space is a second countable paracompact local ringed space which is covered by open sets, each of which is isomorphic to a ringed space defined by an ideal of holomorphic functions on an analytic open subset of $\CC^n$ for some $n>0$, and whose transition maps preserve the sheaves of holomorphic functions.  
%A complex analytic variety is a reduced complex analytic space.
%Throughout this paper, a complex analytic space $\fX$ is a local ringed space, locally defined by an ideal of holomorphic functions on an open set of $\CC^n$. A (complex) analytic variety  $X$ is a reduced complex analytic space, i.e. a complex analytic space defined by a reduced ideal sheaf. 
%We will denote the variety underlying a complex analytic space $\fX$ by $X$. 
%We will use smooth functions to mean $C^\infty$ functions. 
%In case the space is singular, with a stratification by smooth strata, smooth functions are continuous functions that
%are smooth along each stratum. We use analytic functions to mean continuous functions that locally have
%power series expansions in the real and imaginary parts of coordinate variables. We will work with analytic topology unless otherwise mentioned.

%\aufm{Author Name}

%-----------------------------------------------------------------------------
% End of preface.tex
%-----------------------------------------------------------------------------

\mainmatter

\def\fA{\mathfrak{A} }
\def\cD{\mathcal{D} }
\def\bbL{\mathbb{L} } 
\def\gaal{{\gamma\alpha}}
\def\bega{{\beta\gamma}}
\def\albe{{\alpha\beta}}
\def\beal{{\beta\alpha}}
\def\gabe{{\gamma\beta}}
\def\alga{{\alpha\gamma}}
\def\albega{{\alpha\beta\gamma}}
\def\zero{\mathrm{zero} }
\def\Crit{\mathrm{Crit}}
\def\virt{^{\mathrm{vir}}}

\part{Critical virtual manifolds and perverse sheaves}

A critical virtual manifold is a (complex) analytic space which is locally the critical locus of a holomorphic function on a complex manifold of fixed dimension (Definition \ref{1.5}). One may say that a critical virtual manifold is the gluing of Landau-Ginzburg models. Usual complex manifolds are special cases of critical virtual manifolds. As we will show in Part 2, moduli spaces of sheaves on Calabi-Yau 3-folds are critical virtual manifolds. In Part 1, we investigate several interesting structures on critical virtual manifolds, such as orientability, semi-perfect obstruction theory, virtual fundamental class, Donaldson-Thomas type invariant, weighted Euler characteristic, local perverse sheaves of vanishing cycles, their gluing isomorphisms and mixed Hodge modules.

Since a critical virtual manifold $X$ is locally the critical locus $X\lalp$ of a holomorphic function $f\lalp$ on a complex manifold $V\lalp$, it comes with a natural symmetric obstruction theory $E_\alpha=[T_{V\lalp}|_{X\lalp}\to \Omega_{V\lalp}|_{X\lalp}]\to  \bbL_{X\lalp}$ on each $X\lalp$. 
These local symmetric obstruction theories form a semi-perfect obstruction theory (Proposition \ref{1.44}). By \cite{ChLiS}, we then have a virtual fundamental class $[X]\virt\in A_0(X)$ whose degree gives us the Donaldson-Thomas type invariant $DT(X)$, when $X$ is compact. 
By adapting the arguments in \cite{Beh}, we find that $DT(X)$ is indeed the Euler characteristic $\chi_\nu(X)$ of $X$ weighted by the Behrend function $\nu$ (Theorem \ref{1.53}). Then it is natural to consider the categorification problem (Problem \ref{1.55}) which asks if there is a cohomology theory on $X$ whose Euler characteristic equals $DT(X)$.

The local holomorphic functions $f\lalp:V\lalp\to\CC$ with critical loci $X\lalp$ give us local perverse sheaves $P\lalp$ of vanishing cycles on $X\lalp$ (Definition \ref{2.1}) whose Euler characteristic at each point equals the value of the Behrend function (cf. \eqref{2.5a}).
Therefore, if there is a perverse sheaf $P$ on $X$ whose restriction to $X\lalp$ is isomorphic to $P\lalp$, then its hypercohomology has Euler characteristic $\chi_\nu(X)=DT(X)$ (cf. Proposition \ref{2.6}).
We show that there is an isomorphism  $\sigma_\albe:P_\alpha|_{X_\albe}\to P_\beta|_{X_\albe}$ whenever $X_\albe=X\lalp\cap X_\beta\ne \emptyset$ (cf. Corollary \ref{2.11}). 

It is well known that perverse sheaves glue (cf. Proposition \ref{2.14}). 
For a critical virtual manifold $X=\cup X\lalp$ with charts $f\lalp:V\lalp\to\CC$, 
the line bundles $\det T_{V\lalp}|_{X\lalp^{\mathrm{red}}}$ on $X\lalp\ured$ glue to a globally defined line bundle on the reduced space $X\ured$ 
exactly when a 2-cocycle obstruction class $\xi\in H^2(X,\ZZ_2)$ vanishes. We say the critical virtual manifold is orientable if the obstruction class $\xi$ is trivial (Definition \ref{1.17}). 
We prove that the 2-cocycle obstruction for gluing the local perverse sheaves $P\lalp$ to a globally defined perverse sheaf $P$ coincides with the obstruction for the orientability of $X$ (cf. Corollary \ref{2.21}). Therefore an orientable critical virtual manifold $X=\cup X\lalp$ has a perverse sheaf $P$ which is locally the perverse sheaf $P\lalp$ of vanishing cycles on $X\lalp$ (cf. Theorem \ref{2.23}). 

The hypercohomology of a perverse sheaf on an analytic space satisfies the usual nice properties of the cohomology of a compact K\"ahler manifold, such as the Hodge decomposition and the hard Lefschetz property, if the perverse sheaf underlies a polarizable Hodge module (cf. \cite{Sai88, Sai90}). We further prove that an orientable critical virtual manifold has a polarizable mixed Hodge module $\cM$ whose underlying perverse sheaf $\rat(\cM)=P$ is the perverse sheaf above (Theorem \ref{2.32}). The gradation $\hat{P}=\mathrm{gr}^W P$ of $P$ with respect to the weight filtration is a perverse sheaf which underlies a direct sum of polarizable Hodge modules. It is easy to see that the Euler characteristics of the hypercohomology of $P$ and $\hat{P}$ are the same. So we have two solutions to the categorification problem given by $P$ and $\hat{P}$. The latter $\hat{P}$ will be useful in Part 2 when we develop a mathematical theory of the Gopakumar-Vafa invariant.

There is a related notion of a $d$-critical locus defined by Joyce in \cite{Joyce}. We prove that critical virtual manifolds are $d$-critical loci and vice versa (cf. \S\ref{S1.5}).

\bigskip

Here is the layout of Part 1.
In Chapter \ref{ch1}, we introduce critical virtual manifolds, discuss orientability, semi-perfect obstruction theory, Behrend's theorem and the categorification problem. In Chapter \ref{ch2}, we prove the gluing theorem for perverse sheaves and mixed Hodge modules. Chapter \ref{ch3} is devoted to a proof of Theorem \ref{2.20} which is essential for the gluing of perverse sheaves.

\chapter{Critical virtual manifolds}\label{ch1}

In this chapter, we introduce the notion of a critical virtual manifold (Definition \ref{1.5}) and the orientability (Definition \ref{1.17}). We prove that a critical virtual manifold $\fX$ has a natural symmetric semi-perfect obstruction theory (Proposition \ref{1.44}) in the sense of \cite{ChLiS}, which gives us the virtual fundamental class $[\fX]\virt\in A_0(\fX)$ (Definition \ref{1.49}) whose degree $$DT(\fX)=\deg [\fX]\virt$$ is shown to equal the Euler characteristic $$\chi_\nu(\fX)=\sum_{n\in \ZZ}n\cdot \chi_c(\nu^{-1}(n))$$ weighted by the Behrend function $\nu$ (Theorem \ref{1.53}). Then we pose a categorification problem (Problem \ref{1.55}) which will be solved in the subsequent chapter. We end this chapter by compairing the notion of a critical virtual manifold with a $d$-critical locus from \cite{Joyce}.

\bigskip

\section{Background}\label{s1.1}
Recall that a complex manifold is locally modeled on open subsets in $\CC^n$, $n>0$. 
\begin{defi}\label{1.1} \cite[page 14]{GrHa} A \emph{complex manifold} $X$ is a second countable  paracompact Hausdorff space together with an open cover $\{X_\alpha\}$ and homeomorphisms $\varphi_\alpha:X_\alpha\to \varphi_\alpha(X_\alpha)\subset \CC^n$ onto open sets such that $\varphi_\albe:=\varphi_\beta\circ \varphi_\alpha^{-1}$ is holomorphic on $\varphi_\alpha(X_\albe)\subset \CC^n$ for all $\alpha, \beta$. Here $X_\albe=X_\alpha\cap X_\beta$. 
\end{defi}

When $X$ is a smooth projective variety (or more generally a compact K\"ahler manifold), its cohomology $H^*(X,\CC)$ satisfies nice properties like the Hodge decomposition and the Lefschetz decomposition, collectively known as the K\"ahler package \cite{GrHa, KiWo}. These properties are well known to be fundamental tools in complex and algebraic geometry. However when $X$ is singular, these nice properties do not hold unless we use more sophisticated tools.  

\medskip

Our discussion on singular spaces will be based on the notion of analytic spaces. 
\begin{defi}\label{1.2}
A \emph{local ringed space} is a topological space $X$ together with a sheaf $\sO_X$ of rings whose stalks are local rings.
\end{defi}
Given holomorphic functions $f_1,\cdots, f_m$ on an open $V\subset \CC^n$, the common vanishing locus $Y=\zero(f_1,\cdots, f_m)\subset V$ equipped with the sheaf $\sO_Y=\sO_V/(f_1,\cdots,f_m)$ of holomorphic functions on $Y$ is a local ringed space. 
\begin{defi}\label{1.3}
An \emph{analytic space} is a second countable paracompact Hausdorff local ringed space $(X,\sO_\fX)$ together with an open covering $\{X_\alpha\}$ and an isomorphism $\varphi_\alpha:(X_\alpha, \sO_X|_{X_\alpha})\mapright{\cong} (Y_\alpha, \sO_{Y_\alpha})$ onto a local ringed space $(Y_\alpha,\sO_{Y_\alpha})$ in an open $V\lalp\subset\CC^n$ defined by finitely many holomorphic functions for each $\alpha$.
%, such that $\varphi_\albe=\varphi_\beta\circ\varphi_\alpha^{-1}: \varphi_\alpha(X_\albe)\to \varphi_\beta(X_\albe)$ is an isomorphism of local ringed spaces where $X_\albe=X_\alpha\cap X_\beta$. 
A morphism of analytic spaces is a morphism of local ringed spaces, i.e. a continuous map $\varphi:X\to Y$ together with  a sheaf homomorphism $\varphi^*:\varphi^{-1}\sO_Y\to \sO_X$ which induces a local homomorphism of stalks.
\end{defi}

Given an analytic space $(X,\sO_{\fX})$, its reduced space $\fX\ured$ is given by the same topological space and the reduced sheaf $\sO_{\fX\ured}$ which is the quotient of $\sO_{\fX}$ by its radical ideal. %By abuse of notation, we will sometimes denote the reduced space $\fX\lred$ of an analytic space $\fX=(X,\sO_\fX)$ by $X$. 

Although the K\"ahler package may fail for the ordinary cohomology $H^*(X,\CC)$, if we use the hypercohomology $\bbH^*(X,P)$ of a perverse sheaf $P$ on $X$ which underlies a polarizable Hodge module, then we still have the nice properties, thanks to the theory of perverse sheaves and mixed Hodge modules due to Kashiwara, Goresky, MacPherson, Beilinson, Berstein, Deligne and Saito, to name a few \cite{BBD, GM, Sai88, Sai90}.
So it is reasonable to raise the following question:

\medskip

\noindent \emph{How do we get a natural perverse sheaf $P$
on an analytic space $\fX$ underlying a polarizable Hodge module?}

\medskip

One obvious way is to use the intersection cohomology complex $IC_X$ which extends the constant variation of Hodge structure $\QQ_U$ on an open dense smooth subset $U\subset X$. For the Donaldson-Thomas theory which is the focus of this paper, the spaces we deal with are locally the critical locus $\mathrm{Crit}(f)$ of a holomorphic function $f:V\to \CC$ on a complex manifold $V$ (cf. \cite{JoSo}). In this case, the perverse sheaf $P_f=  \phi_f\QQ[\dim V-1]$ of vanishing cycles is much more relevant than the constant sheaf $\QQ$ or the intersection cohomology complex $IC_X$ (cf. \cite{Beh}). 

In Part 1 of this paper, we will introduce the notion of a critical virtual manifold and show that there are natural choices of perverse sheaves underlying polarizable mixed Hodge modules. In Part 2, we will prove that  moduli spaces for the Donaldson-Thomas theory are critical virtual manifolds. 
We thus will obtain a cohomology theory for the Donaldson-Thomas invariant which moreover gives us a mathematical theory of the Gopakumar-Vafa invariant. 

\bigskip

\section{Definition of Critical Virtual Manifolds}\label{S1.2}

A critical virtual manifold is an analytic space locally modeled on the critical locus $\Crit(f)$ of a holomorphic function $f:V\to \CC$ on a complex manifold $V$.
\begin{defi}\label{1.4}
(1) An \emph{LG pair} is a complex manifold $V$ together with a holomorphic function $f:V\to \CC$ which has only one critical value $0$.

(2) Two LG pairs $(V_1,f_1)$ and $(V_2,f_2)$ are called \emph{equivalent} if there exists a biholomorphic map $\varphi:V_1\to V_2$ satisfying $f_2\circ\varphi=f_1$.

(3) The \emph{critical locus} of an LG pair $(V,f)$ is the analytic space in $V$ defined by the ideal $(df)$ generated by the partial derivatives of $f$. We will sometimes denote the critical locus $\Crit(f)$ by $\fX_f$.

(4) The reduced analytic space $\fX_f\ured$ of the critical locus $\fX_f$ is the analytic space defined by the radical ideal of $(df)$. %We denote $(\fX_f)\lred$ sometimes by $X_f$. 
\end{defi}
Here LG stands for Landau-Ginzburg. % $(V,f)$ as a Landau-Ginzburg model with superpotential $f$. 
Perhaps a reader familiar with mirror symmetry may think of $(V,f)$ in Definition \ref{1.4} (1) as a \emph{Landau-Ginzburg model} with superpotential $f$. A critical virtual manifold is the gluing of Landau-Ginzburg models.

%We are now ready to define a critical virtual manifold.
\begin{defi}\label{1.5}
A \emph{critical virtual manifold} is an analytic space $(X,\sO_\fX)$ together with an open covering $\fX=\cup\lalp \fX_\alpha$ and a closed embedding $\varphi_\alpha:\fX_\alpha\hookrightarrow V_\alpha$ into a complex manifold $V_\alpha$ for each $\alpha$ 
satisfying the following:
\begin{enumerate}
\item $\varphi_\alpha$ is an isomorphism onto the critical locus $\fX_{f_\alpha}$ for an LG pair $(V\lalp,f\lalp)$;
\item for each pair $(\alpha,\beta)$ of indices and
$\fX_\albe=\fX_\alpha\cap\fX_\beta$, we have an open neighborhood $V_\albe$ (resp. $V_\beal$) of $\varphi_\alpha(\fX_\albe)$ in $V_\alpha$ (resp. $\varphi_\beta(\fX_{\beta\alpha})$ in $V_\beta$) and a biholomorphic map $\varphi_\albe:V_\albe\to V_\beal$ that fit into the commutative diagram
\beq\label{1.0}\xymatrix{
\fX_\albe\ar@{=}[dd]\ar[r]^{\varphi_\alpha}&V_\albe\ar[dd]_{\varphi_\albe}\ar@{^(->}[r]\ar[dr]& V_\alpha \ar[d]^{f_\alpha}\ar[d] \\
 && \CC\\
\fX_\beal\ar[r]_{\varphi_\beta}&V_\beal \ar[ur]\ar@{^(->}[r]& V_\beta \ar[u]_{f_\beta}
}\eeq
i.e. $\varphi_\albe\circ \varphi_\alpha|_{\fX_{\albe}}=\varphi_\beta|_{\fX_{\beal}}$ and $f_\beta\circ\varphi_{\albe}=f_\alpha|_{V_\albe}$;
\item $\varphi_\albe=\varphi_{\beta\alpha}^{-1}$ for each pair $(\alpha,\beta)$ of indices and $\varphi_{\alpha\alpha}=\id_{V\lalp}$.
\end{enumerate}
\end{defi}
In other words, we have LG pairs $(V_\alpha,f_\alpha)$ whose critical loci are $\fX_\alpha$ and the restrictions of $(V\lalp,f\lalp)$ and $(V_\beta,f_\beta)$ to some open neighborhoods of $\fX_\albe$ are equivalent. 
%We will sometimes denote by $X=\fX\lred$, $X_\alpha=(\fX_\alpha)\lred$ and $X_\albe=(\fX_\albe)\lred$ the reduced analytic spaces of $\fX$, $\fX_\alpha$ and $\fX_\albe$ respectively.

\begin{defi}\label{1.6}
A diagram $(X\lalp\mapright{\varphi\lalp} V\lalp\mapright{f\lalp}\CC)$ from Definition \ref{1.5} (1) is called a \emph{chart} of the critical virtual manifold $\fX$. When $\dim V\lalp=r$ for each $\alpha$, we say the critical virtual manifold $\fX$ has an \emph{atlas of dimension} $r$.
\end{defi}

\begin{rema}\label{1.7}
When $\fX_\albega=\fX_\alpha\cap \fX_\beta\cap \fX_\gamma\ne \emptyset$, we can find an open neighborhood
$V_\albega$ of $\varphi\lalp(\fX_\albega)$ in $V_\alpha$ such that the composition
\beq\label{mapcoc} \varphi_\albega=\varphi_\gaal\circ\varphi_\bega\circ\varphi_\albe:V_\albega\lra V_\alpha\eeq
is biholomorphic onto its image. 
For instance, we may let
$$V_\albega=\varphi_\beal(V_\beal\cap V_\bega\cap \varphi_\gabe(V_\gaal\cap V_\gabe)).$$
By Definition \ref{1.5} (2), $\varphi_\albega\circ \varphi\lalp|_{\fX_\albega}=\varphi\lalp|_{\fX_\albega}$ as morphisms of analytic spaces. Note that we do not require the cocycle condition that $\varphi_\albega$ should be the identity map of $V_\albega$. The cocycle condition $\varphi_\albega=\mathrm{id}$ is guaranteed to hold only at $\varphi_\alpha(\fX_\albega)$.
\end{rema}

\begin{exam}\label{1.8} Complex manifolds (Definition \ref{1.1}) are critical virtual manifolds by letting $\fX\lalp=V\lalp$, $f\lalp=0$ for any open cover $X=\cup_\alpha X\lalp$. 
\end{exam}

Nontrivial examples arise from the Landau-Ginzburg theory.
\begin{exam}\label{1.9}
Let $E$ be a holomorphic vector bundle over a complex manifold $W$ equipped with a holomorphic section $s:W\to E$. 
Let $\pi:V=E^*\to W$ be the dual vector bundle of $E$ over $W$, which is a complex manifold. Then $s$ defines a holomorphic function $f:V\to \CC$ defined by $f(v)=\langle s(\pi(v)),v\rangle$.  Obviously the critical locus $\fX_f=\Crit(f)$ is a critical virtual manifold. 

Suppose the analytic space $\zero(s)$ defined by the vanishing of $s$ is smooth, i.e. the Jacobian matrices $(\frac{\partial s_i}{\partial x_j})$ have full rank along $\zero(s)=\zero(s_1,\cdots,s_r)$, where $\{x_j\}$ are local coordinates of $W$ and $\{s_i\}_{1\le i\le r}$ are the coordinate functions of the section $s$ after choosing a local trivialization $E|_U\cong U\times \CC^r$ over an open $U\subset W$. Then the critical locus $\fX_f=\Crit(f)$ is precisely the complete intersection $\zero(s)\subset W$. Indeed, in local coordinates, 
$$f=\langle s,v\rangle =\sum_{i=1}^r s_iv_i$$
where $\{v_i\}$ are the vertical coordinates of $\pi$. Hence 
$$df=\sum_{i,j}v_i\frac{\partial s_i}{\partial x_j}dx_j + \sum_{i=1}^r s_i dv_i=0$$
amounts to $v_i=0$ and $s_i=0$ for all $i$. Thus $\Crit(f)=\zero(s)\subset W$

In this way, all smooth complete intersections in complex manifolds come with nontrivial critical virtual manifold structures.  
\end{exam}

\begin{exam}\label{1.9e}
Let $$\Crit(yz^2)=\{(y,z)\,|\, yz=z^2=0\}\subset \CC^2=\{(y,z)\,|\,y,z\in \CC\}$$ and $X$ be the closure of $\Crit(yz^2)$ in the projective plane $\PP^2$. In fact, $X$ is the projective line $\PP^1=X\ured$ with an embedded point at $(0,0)\in \CC^2\subset\PP^2$. We claim that $X$ is a critical virtual manifold.
Let $(x:y)\in \PP^1$  be the homogeneous coordinates and $((x:y),z)$ denote the coordinates of the line bundle $$\sO_{\PP^1}(1)\to \PP^1,\quad ((x:y),z)\mapsto (x:y).$$ 
Let $D_1\subset\CC$ denote the unit disk centered at $0$ and consider the charts 
$$U_0=\{((1:y),z)\,|\, yz=0=z^2, |y|<1\}=\Crit(yz^2)\subset D_1\times\CC=V_0,$$
$$U_+=\{ ((1:y),z)\,|\, y\notin \RR_{\ge 0}, z=0\}=\Crit(z^2)\subset (\CC-\RR_{\ge 0})\times \CC=V_+,$$
$$U_-=\{ ((1:y),z)\,|\, y\notin \RR_{\le 0}, z=0\}=\Crit(z^2)\subset (\CC-\RR_{\le 0})\times \CC=V_-,$$
$$U_\infty=\{((x:1),z)\,|\,z=0, |x|<1\}\subset \Crit(z^2)\subset D_1\times\CC=V_\infty.$$
Then the holomorphic maps $\varphi_{\pm,\infty}=\id$, $\varphi_{+,-}=\id$, and 
$$\varphi_{0,\pm}:V_{0,\pm}\lra V_{\pm,0}, \quad (y,z)\mapsto (y,\sqrt{y}z)$$
give us a critical manifold structure with $f_0=yz^2$, $f_{\pm}=z^2$, $f_\infty=z^2$.
\end{exam}

In Part 2, we will prove that if $Y$ is a Calabi-Yau 3-fold (a smooth projective variety of dimension 3 whose canonical line bundle $K_Y$ is trivial), then moduli spaces of simple sheaves on $Y$ are all critical virtual manifolds.

\begin{rema}\label{1.10} 
In \cite{Joyce}, Joyce introduced the notion of a $d$-critical locus and 
in \cite{BrBuJo}, Brav, Bussi and Joyce proved that a $(-1)$-shifted symplectic derived scheme
is an algebraic $d$-critical locus. In \S\ref{S1.5}, we will prove that critical virtual manifolds are $d$-critical loci and vice versa. In particular, the analytic space associated to each $(-1)$-shifted symplectic derived scheme admits a critical virtual manifold structure. 
\end{rema}
%\begin{rema}
%In \cite{ChKi}, we will introduce the notion of an algebraic critical virtual manifold. Then we will prove that for orientable algebraic critical virtual manifolds, we have the motivic invariant weighted by the motivic Milnor fiber.
%\end{rema}

\begin{rema}
Although it was not explicitly stated as a definition, the notion of a critical virtual manifold was the key in our solution to the categorification problem of the Donaldson-Thomas invariants in our previous draft in 2012. 
\end{rema}

\medskip

How can we compare two critical virtual manifold structures of different dimensions? 
\begin{defi}\label{1.63}
Let $X$ be an analytic space with an open cover $\{X\lalp\}$. Two critical virtual manifold structures
$$\{(X\lalp=\Crit(f\lalp)\sub V\lalp\mapright{f\lalp}\CC, \varphi_\albe:V_\albe\to V_\beal)\},$$  
$$\{(X\lalp=\Crit(g\lalp)\sub W\lalp\mapright{g\lalp}\CC, \psi_\albe:W_\albe\to W_\beal)\}$$
with $\dim V\lalp \le \dim W\lalp$ are \emph{compatible} if there exist closed embeddings
$$\imath\lalp:V\lalp\lra W\lalp,\quad \forall \alpha$$
such that $g\lalp\circ\imath\lalp=f\lalp$ and 
$$\imath_\beta\circ \varphi_\albe = \psi_\albe\circ \imath\lalp, \quad \forall \alpha,\beta.$$
\end{defi}

The two charts $(V\lalp,f\lalp)$ and $(W\lalp,g\lalp)$ are related by the following lemma. 
\begin{lemm}\label{1.60}  Let $W\subset V$ be a complex submanifold in a complex manifold. Let $f:V\to \CC$ be a holomorphic function and $g=f|_W$. Suppose the critical loci $\fX_f=\Crit(f)$ and $\fX_g=\Crit(g)$ are equal. Let $x$ be a closed
point of $\fX_f=\fX_g$. Then %we have the following:\\ (1) T
there is a coordinate system $\{z_1,\cdots, z_r\}$ of $V$ centered at $x$ such that $W$ is defined by the vanishing of $z_1,\cdots,z_m$ and 
$$f=q(z_1,\cdots,z_m)+g(z_{m+1},\cdots,z_r),\quad q(z_1,\cdots,z_m)=\sum_{i=1}^mz_i^2, \quad g=f|_W.
$$ 
%(2) If we let $X\lred$ be the reduced analytic space of $\fX_f=\fX_g$, we have a quasi-isomorphism
%$$\xymatrix{
%[T_W|_{X\lred}\ar[r]^{d(dg)} \ar@{^(->}[d] & T_W^*|_{X\lred}] \ar@{^(->}[d] \\
%[T_V|_{X\lred}\ar[r]^{d(df)} & T_V^*|_{X\lred}]}$$
%so that we have an isomorphism $(\det T_W|_{X\lred})^2\cong (\det T_V|_{X\lred})^2$. \\
%(3) If $\psi:V\to V'$ is an equivalence that restricts to an equivalence $\psi_W:W\to W'$, then we have a commutative diagram
%\[\xymatrix{
%(\det T_W|_{X\lred})^2 \ar[rr]^{\det(d\psi_W)^2}\ar[d]_\cong && (\det T_{W'}|_{X\lred})^2\ar[d]^\cong\\
%(\det T_V|_{X\lred})^2\ar[rr]^{\det(d\psi_W)^2} && (\det T_{V'}|_{X\lred})^2.
%}\]
%(2) Let $\Phi:V\to V$ be an automorphism of $V$  such that $\Phi|_W=\id_W$ and  $f\circ\Phi=f$. Then $\Phi^*:A\ud_f\to A\ud_f$ is $\det(d\Phi|_{X})\,\id_{A\ud_f}$ and $\det(d\Phi|_{X})=\pm 1$.
\end{lemm}
\begin{proof}
%By Proposition \ref{prop1} (2), the lemma is a consequence of the following claim:
We choose coordinates $\{y_1,\cdots,y_r\}$ of $V$ centered at $x$ such that $W$ is defined by the vanishing of $y_1,\cdots,y_m$. Let $I$ be the ideal generated by $y_1,\cdots, y_m$. Since $\fX_f=\fX_g$, i.e. $(df)=(dg)+I$, we have 
$$\frac{\partial f}{\partial y_i}\Big|_W=\sum_{j=m+1}^r a_{ij}\frac{\partial g}{\partial y_j},\quad i=1,\cdots,m$$
for some functions $a_{ij}$ regular at $x$.
By calculus,  we have
$$f=g(y_{m+1},\cdots,y_r)+\sum_{i=1}^m\frac{\partial f}{\partial y_i}\Big|_W\cdot y_i+I^2 $$
$$=g(y_{m+1},\cdots,y_r)+\sum_{j=m+1}^r\frac{\partial g}{\partial y_j}\left( \sum_{i=1}^ma_{ij}y_i\right)+I^2 $$
$$=g(z_{m+1},\cdots, z_r)+\sum_{i,k=1}^m b_{ik}y_iy_k$$ where  $z_j=y_{j}+\sum_{i=1}^ma_{ij}y_i$ for $j\ge m+1$
and $b_{ik}$ are some functions holomorphic near $x$. Since the kernel of the Hessian of $f$ at $x$ is the tangent space of $\fX_f=\fX_g\subset W$ at $x$, the quadratic form $q=\sum_{i,k=1}^m b_{ik}y_iy_k$ is nondegenerate near $x$. Hence we can diagonalize $q=\sum_{i=1}^m z_i^2$ by changing the coordinates $y_1,\cdots, y_m$ to new coordinates $z_1,\cdots,z_m$.  It follows that $z_1,\cdots,z_r$ is the desired coordinate system. 
%The remaining (2) and (3) are straightforward. 
\end{proof}

\bigskip

\section{Orientability of a critical virtual manifold}\label{S1.3}
For a real (Riemannian) manifold $M$ of dimension $n$, its orientation bundle is the top exterior power $\det_\RR T_M=\wedge^n_\RR T_M$ of the tangent bundle $T_M$ and its associated principal $O(1)$-bundle has $\{\pm 1\}$-valued locally constant transition functions since $O(1)=\{\pm 1\}=\ZZ_2$. So the orientation bundle is determined by a class $\xi\in H^1(M,\ZZ_2)$ whose vanishing is equivalent to the orientability of $M$, i.e. the triviality of $\det_\RR T_M$. There is an analogous story for critical virtual manifolds. 

Let $\fX$ be a critical virtual manifold and $(\fX\lalp\mapright{\varphi\lalp} V\lalp \mapright{f\lalp} \CC)$ be a chart. Let 
$$K\lalp^\vee=\varphi_\alpha^*\det T_{V\lalp}|_{X\ured\lalp}$$
be the dual of the canonical bundle of $V\lalp$ restricted to the reduced space $X\ured\lalp$ of $X\lalp$. 
These line bundles are related by the isomorphisms $\xi_\albe:K\lalp^\vee|_{\fX_\albe\ured}\to K_\beta^\vee|_{\fX_\albe\ured}$ defined by
\beq\label{1.11} 
\xi_\albe:=(\varphi\lalp|_{\fX_\albe\ured})^*\det (d\varphi_\albe):\varphi_\alpha^{*}\det T_{V_\alpha}|_{\fX_\albe\ured} \lra \varphi_\beta^{*}\det T_{V_\beta}|_{\fX_\albe\ured}
\eeq

Recall that $\varphi_\albega$ denotes the composition $\varphi_\gaal\circ\varphi_\bega\circ\varphi_\albe:V_\albega\lra V_\alpha.$
\begin{prop}\label{1.12}
For $\fX_\albega=\fX_\alpha\cap \fX_\beta\cap \fX_\gamma\ne \emptyset$, 
\beq\label{1.13}
\xi_\albega=\xi_\gaal\circ\xi_\bega\circ\xi_\albe=(\varphi\lalp|_{\fX_\albega\ured})^*\det(d\varphi_\albega)
\eeq
is locally constant, taking values in $\{\pm 1\}$.
\end{prop}
\begin{proof}
Apply Lemma \ref{1.14}, letting $V=V\lalp$, $U=V_\albega$, $f=f\lalp$, $\varphi=\varphi_\albega$.\end{proof}

\begin{lemm}\label{1.14}
Let $(V,f)$ be an LG pair (Definition \ref{1.4}). Let $U\sub V$ be open and $\varphi: U\to V$ be a biholomorphic map onto its image, such that $f\circ \varphi=f|_U$ and $\varphi|_{\fX_f\cap U}=\mathrm{id}_{\fX_f\cap U}$ where $\fX_f=\Crit(f)$. 
Then $\det(d\varphi)|_{X_f\ured\cap U}$ is locally constant, taking values $\pm1$. 
\end{lemm}
\begin{proof}
Let $H_x=d(df)$ denote the Hessian of $f$ at $x\in {\fX_f}\cap U$. Then $H_x$ is a 
symmetric bilinear form whose kernel is $T_{\fX_f}|_x$. After diagonalizing $H_x$, 
we can choose a subspace $W_x$, complementary to $T_{\fX_f}|_x$, such that the Hessian $H_x$ is nondegenerate on $W_x$.
Since $\varphi|_{\fX_f}=\id_{\fX_f}$, we can write 
\beq\label{1.14d}
d\varphi|_x=\left( \begin{matrix}\mathrm{id} & A\\ 0& B\end{matrix} \right):T_V|_x=T_{\fX_f}|_x\oplus W_x\to T_{\fX_f}|_x\oplus W_x =T_V|_x
\eeq 
where $A:W_x\to T_{\fX_f}|_x$ and $B:W_x\to W_x$ are homomorphisms of vector spaces. 
Since $f\circ\varphi=f|_U$, $B$ preserves the nondegenerate symmetric bilinear form $H_x$ 
and hence $B$ lies in the orthogonal group $O_{H_x}(W_x)$ with respect to $H_x$. 
Since $\det (d\varphi|_x)=\det B$ and an orthogonal matrix has determinant $\pm 1$, 
we find that $\det(d\varphi|_x)\in \{\pm 1\}$ as desired.
\end{proof}

\def\albegade{{\alpha\beta\gamma\delta}}
Since $\xi_\albega=\xi_\gaal\circ \xi_\bega\circ \xi_\albe$ and $\xi_\albe=\xi_\beal^{-1}$ for each pair $(\alpha,\beta)$,
$$
(d\xi)_\albegade=\xi_\albega^{-1}\xi_{\alpha\gamma\delta}^{-1}\xi_{\alpha\beta\delta}\xi_{\beta\gamma\delta} $$
$$
=(\xi_\beal\xi_\gabe\xi_\alga)(\xi_\gaal\xi_{\delta\gamma}\xi_{\alpha\delta})(\xi_{\delta\alpha}\xi_{\beta\delta}\xi_\albe)\xi_{\beta\gamma\delta} $$
$$
=(\xi_\beal\xi_{\gamma\beta}\xi_{\delta\gamma}\xi_{\beta\delta}\xi_\albe)\xi_{\beta\gamma\delta} =(\xi_\albe^{-1}\xi_{\beta\gamma\delta}^{-1}\xi_\albe)\xi_{\beta\gamma\delta}  $$
$$
=\xi_\albe^{-1}\circ\xi_\albe = 1
$$
where the second to the last equality comes from $\xi_{\beta\gamma\delta}=\pm 1$. 
Therefore, $\{\xi_\albega\}$ is a \v{C}ech 2-cocycle with values in $\{\pm 1\}$ and defines a class
\beq\label{1.15} 
\xi:=[\xi_\albega]\in H^2(X,\ZZ_2)\eeq
whose vanishing amounts to, possibly after a refinement of the open cover $\{\fX_\alpha\}$,
the existence of a 1-cochain $\mu=\{\mu_\albe\}$ taking values in $\{\pm 1\}$ such that 
for $\overline{\xi}_\albe=\mu_\albe\xi_\albe$, we have 
\beq\label{1.16}
\overline{\xi}_\albega:= \overline{\xi}_\gaal\circ\overline{\xi}_\bega\circ\overline{\xi}_\albe=1
\eeq
whenever $\fX_\albega\ne \emptyset$. In particular, the 1-cocycle $\{\overline{\xi}_\albe\}$ with values in $\sO^*$ glues the local line bundles $\{K_\alpha^\vee=\varphi_\alpha^*\det T_{V_\alpha}|_{X\lalp\ured}\}$ to a globally defined line bundle $K^\vee_{X}$ on $\fX\ured$.

\begin{defi}\label{1.17}
We say a critical virtual manifold $\fX$ is \emph{orientable} if the class $\xi\in H^2(X,\ZZ_2)$ defined in \eqref{1.15} is zero. 
When $\fX$ is orientable, a line bundle $K^\vee_{\fX}$ on $X\ured$ obtained by gluing $\{K^\vee_\alpha\}$ by $\{\overline{\xi}_\albe=\mu_\albe\xi_\albe\}$ above is called an \emph{orientation bundle} of $\fX$.
\end{defi}

The issue of orientability will be crucial when we construct a natural perverse sheaf and a mixed Hodge module on a critical virtual manifold $\fX$ by gluing locally defined perverse sheaves and mixed Hodge modules of vanishing cycles on local charts $\fX_\alpha$.

\begin{rema}\label{1.18}
An orientation bundle of a critical virtual manifold is not unique. One can always twist $\{\overline{\xi}_\albe\}$ by a $\ZZ_2$-valued 1-cocycle $\{\lambda_\albe\}$ to obtain a new orientation bundle defined by the gluing isomorphisms $\{\lambda_\albe\overline{\xi}_\albe\}$. In fact, since $\xi_\albega\in \{\pm 1\}$, $\{\xi_\albe^2\}$ is a 1-cocycle and glues the locally defined line bundles $\{(K_\alpha^\vee)^2\}$ on $\fX_\alpha\ured$ to a globally defined line bundle $(K_{\fX}^\vee)^2$ on $X\ured$. Hence the square of an orientation bundle is canonically defined. So we find that when a critical virtual manifold $\fX$ is orientable, the set of orientation bundles is an $H^1(X,\ZZ_2)$-orbit by the exact sequence
$$\cdots \lra H^1(\fX,\ZZ_2)\lra H^1(\fX,\sO_\fX^*)\mapright{2} H^1(\fX,\sO_\fX^*)\lra \cdots$$
from the short exact sequence 
$$1\lra \{\pm 1\} \lra \sO_\fX^*\mapright{2} \sO_\fX^*\lra 1.$$
\end{rema}

\begin{exam} Every complex manifold is an orientable critical virtual manifold. 
The critical virtual manifold in Example \ref{1.9e} is not orientable because
$$\det(d\varphi_{\infty,+,-})=1 \and \det(d\varphi_{0,+,-})=\left\{ \begin{matrix}1 & \mathrm{Im}(y)>0\\
-1 & \mathrm{Im}(y)<0\end{matrix}\right. $$
define a nontrivial class in $H^2(\PP^1,\ZZ_2)$.  
\end{exam}

\bigskip
\def\fc{\mathfrak{c} }
\section{Semi-perfect obstruction theory and Behrend's theorem}\label{S1.4}
In this section, we prove that a critical virtual manifold is equipped with a natural semi-perfect obstruction theory (cf. \cite{ChLiS}) which gives us a Donaldson-Thomas type invariant of a critical virtual manifold $\fX$.
We will also prove an analogue of Behrend's theorem for critical virtual manifolds which says that the Donaldson-Thomas type invariant is the weighted Euler number of the analytic space $\fX$ with the weight function given by the Milnor numbers on the local charts. 

\subsection{Behrend function}\label{S1.4.1}
We first recall the theory of Behrend functions from \cite{Beh}. We assume that $\fX$ admits a \emph{global embedding into a complex manifold $\PP$}. Let $C_{\fX/\PP}$ be the normal cone of $\fX$ in $\PP$. Then $C_{\fX/\PP}$ defines the \emph{characteristic cycle}
\beq\label{1.19}
\fc_\fX=\sum_{C'} (-1)^{\dim \pi(C')} \mathrm{mult}(C') \pi(C')
\eeq
where $\pi:C_{\fX/\PP}\to \fX$ is the obvious projection and the sum is over all irreducible components $C'$ of $C_{\fX/\PP}$. Hence $\pi(C')$ denotes the image of $C'$ by $\pi$ which is a prime cycle of $X$. Also, $\mathrm{mult}(C')$ is the length of $C'$ at the generic point of $C'$. 

For a chart $(\fX\lalp\mapright{\varphi\lalp} V\lalp\mapright{f\lalp} \CC)$, we can use the normal cone $C_{\fX\lalp/V\lalp}$ to define the characteristic cycle $\fc_{\fX\lalp}$ by the recipe of \eqref{1.19}. 
By \cite[Proposition 1.1]{Beh}, we have
$$\fc_\fX|_{\fX_\alpha}= \fc_{\fX\lalp}$$
because the characteristic cycle is independent of the embedding into a complex manifold. 

For a prime cycle $W$ of dimension $p$ on $X$, the Nash blow-up $\mu:\tilde{W}\to W$ gives us a constructible function
\beq\label{1.20} Eu(W)=\int_{\mu^{-1}(p)} c(\tilde{T})\cap s(\mu^{-1}(p),\tilde{W})\eeq
where $\tilde{T}$ is the dual of the universal quotient bundle on $\tilde{W}$ in the Grassmannian of rank $p$ quotients of $\Omega_W$ (or $\Omega_\PP$). Here $c(\tilde{T})$ is the total Chern class of $\tilde{T}$ and $s(\mu^{-1}(p),\tilde{W})$ is the Segre class of the normal cone to $\mu^{-1}(p)$ in $\tilde{W}$. Now the Behrend function is defined as the $\ZZ$-valued constructible function (\cite[Definition 1.4]{Beh})
\beq\label{1.21} \nu_\fX=Eu(\fc_\fX)\eeq
When $\fX$ is a critical virtual manifold so that $\fX$ is locally the critical locus $\fX_\alpha\cong \Crit(f_\alpha)$ of a holomorphic function $f_\alpha:V_\alpha\to \CC$, by \cite[Corollary 2.4 (iii)]{PaPr} or \cite[(4)]{Beh}, 
\beq\label{1.22}
\nu_\fX(x)=(-1)^{\dim V_\alpha}(1-\chi(MF_x))
\eeq
where $\chi(MF_x)$ is the Euler characteristic of the Milnor fiber $MF_x$ of $f_\alpha$ at $x$ (cf. \eqref{2.4a}).

The constructible function $\nu_\fX$ defines the weighted Euler characteristic
\beq\label{1.23} \chi_\nu(\fX)=\sum_{n\in \ZZ}n\cdot\chi_c(\nu_\fX^{-1}(n))\eeq
where $\chi_c(\cdot)=\sum (-1)^j \dim H^j_c(\cdot)$ is the Euler characteristic with compact support.

The cotangent bundle $\Omega_\PP$ of a complex manifold $\PP$ admits a canonical symplectic structure: for any system of local coordinates $x_1,\cdots, x_n$ and the coordinates $p_1,\cdots, p_n$ for the basis $dx_1,\cdots,dx_n$ for the fibers of $\pi:\Omega_\PP\to \PP$, $\sum_{i=1}^ndp_i\wedge dx_i$ is a symplectic form which is independent of the choice of the coordinates $\{x_1,\cdots,x_n\}$. 

We say a closed analytic subset $W\subset \Omega_\PP$ is \emph{conic Lagrangian} if the Euler vector field $\theta=\sum p_i\frac{\partial}{\partial p_i}$ is tangent to $W$, $\dim W=\dim \PP$ and the symplectic form vanishes identically on the smooth locus of $W$. A conic Lagrangian subset in the cotangent bundle $\Omega_\PP$ is completely determined by its intersection with the zero section or the characteristic cycle.
\begin{lemm}\label{1.24} \cite[Lemma 4.2]{Beh}
Let $W\subset \Omega_\PP$ be a closed irreducible analytic subset. Let $\bar W=\pi(W)\subset \PP$ be its image and let $\ell(\bar W)\subset \Omega_\PP$ denote the closure of the conormal bundle of any smooth dense open subset of $\bar W$. If $W$ is conic Lagrangian, then $\ell(\bar W)=W$.
\end{lemm}
Here the conormal bundle of a submanifold $\imath:S\hookrightarrow \PP$ means the kernel bundle $\mathrm{ker}(\Omega_\PP|_S\mapright{\imath^*}\Omega_S)$ of the pullback homomorphism.

Note that conic Lagrangian is a local property. The conormal bundle construction $\ell(\bar W)$ for a prime cycle $\bar W$ extends to a homomorphism 
\beq\label{1.27} \ell:Z_*(\PP)\lra Z_{\dim \PP}(\Omega_\PP). \eeq
\begin{prop}\label{1.25} \cite[Propositions 1.12 and 4.6]{Beh} 
If we apply the homomorphism $\ell$ to the characteristic cycle $\fc_\fX$ in \eqref{1.19} together with the Gysin homomorphism
\[ 0^!_{\Omega_\PP}:Z_{\dim \PP}(\Omega_\PP)\lra A_0(\PP) \]
and the degree map $\deg:A_0(\PP)\mapright{p_*}A_0(pt)=\ZZ$ where $p:\PP\to \Spec\, \CC$ is the constant morphism, then we get the Euler characteristic $\chi_\nu(\fX)$ weighted by the Behrend function, i.e.
\beq\label{1.26} \chi_\nu(\fX)=\deg 0^!_{\Omega_\PP}[\ell(\fc_\fX)].\eeq
\end{prop}

\begin{rema}\label{1.58}
By \cite[Proposition 1.1]{Beh}, the Behrend function $\nu_\fX$ depends only on the analytic space structure $(\fX,\sO_\fX)$. Therefore, the weighted Euler characteristic $\chi_\nu(\fX)$ is an invariant of the analytic space $(\fX,\sO_\fX)$.
\end{rema}

In the subsequent subsection, we will see that $\ell(\fc_\fX)$ is actually the obstruction cone in $\Omega_\PP$ for a natural semi-perfect obstruction theory when $\fX$ is a critical virtual manifold, so that
\beq\label{1.28} 0^!_{\Omega_\PP}[\ell(\fc_\fX)]=[\fX]\virt.\eeq
Therefore, the weighted Euler characteristic $\chi_\nu(\fX)$ is the virtual invariant of $\fX$ which satisfies  expected properties such as deformation invariance.

\def\vareps{\varepsilon}

\subsection{Semi-perfect obstruction theories on critical virtual manifolds}\label{S1.4.2}
In this subsection, we show that a critical virtual manifold admits a natural semi-perfect obstruction theory and a virtual fundamental class whose degree gives us a Donaldson-Thomas type invariant.

Recall the following definition from \cite{BeFa}.
\begin{defi}\label{1.28a} A \emph{perfect obstruction theory} on an analytic space $\fX$ refers to a morphism $\phi:E\to \bbL_\fX$ in the derived category $D(\sO_\fX)$ such that \begin{enumerate}
\item $E$ is locally isomorphic to a 2-term complex $[E^{-1}\to E^0]$ of locally free sheaves;
\item $H^{-1}(\phi)$ is surjective and $H^0(\phi)$ is an isomorphism.
\end{enumerate}
Here $\bbL_\fX=\bbL_X^{\ge -1}$ denotes the cotangent complex of $\fX$ in \cite{Illusie}, truncated at $\ge -1$. The perfect obstruction theory $\phi:E\to \bbL_\fX$ is called \emph{symmetric} if there is an isomorphism $\theta:E\to E^\vee[1]$ satisfying $\theta^\vee[1]=\theta$.
\end{defi}
By \cite[Remark 3.7]{Beh}, the obstruction sheaf $Ob_X=H^1(E^\vee)$ of a symmetric obstruction theory is canonically isomorphic to the cotangent sheaf $\Omega_X$. 

\begin{defi}\label{1.35}\cite{ChLiS} Let $\fX$ be an analytic space. A \emph{semi-perfect obstruction theory} of $\fX$ consists of an open covering $\{\fX_\alpha\}$ and perfect obstruction theories
\[ \phi_\alpha:E_\alpha\lra \bbL_{\fX\lalp}\]
for each $\alpha$ such that \begin{enumerate}
\item for each pair $\alpha,\beta$ of indices, there is an isomorphism 
\beq\label{1.29} \psi_\albe:H^1(E_\alpha^\vee)|_{\fX_\albe}\to H^1(E_\beta^\vee)|_{\fX_\albe} \eeq
satisfying $\psi_{\alpha\alpha}=\mathrm{id}$, $\psi_{\beta\alpha}=\psi_\albe^{-1}$ and $\psi_\bega\circ\psi_\albe=\psi_\alga$ for each triple $\alpha,\beta,\gamma$;
\item for each pair $\alpha, \beta$, the perfect obstruction theories $E_\alpha|_{\fX_\albe}$ and $E_\beta|_{\fX_\albe}$ on $\fX_\albe$ give the same obstruction assignment via $\psi_\albe$.
\end{enumerate}
By (1) above, we obtain a sheaf $Ob_X$ which is the gluing of $\{Ob_{X\lalp}=H^1(E\lalp^\vee)\}$ via $\psi_\albe$ and is called the \emph{obstruction sheaf} of $X$ with respect to the semi-perfect obstruction theory.
\end{defi}

The second condition requires further explanation.
\begin{defi}\label{1.30}
Let $x\in \fX$ be a point in an analytic space. An \emph{infinitesimal lifting problem} of $\fX$ at $x$ consists of 
\begin{enumerate}
\item an extension $0\to I\to B\to \bar B\to 0$ of Artin local rings by an ideal $I$ with $I\cdot m_B=0$;
\item a morphism $\bar{g}:\Spec \bar B\to \fX$ sending the unique closed point to $x$. 
\end{enumerate}
\end{defi}
Let $\Delta=\Spec B$ and $\bar \Delta=\Spec \bar B$.
By \cite[\S4]{BeFa}, for any infinitesimal lifting problem, there is a canonical obstruction
\beq\label{1.31}
\omega(\bar g, B,\bar B):=(\bar g^*\bbL_\fX\mapright{\bar g} \bbL_{\bar \Delta}\to \bbL_{\bar\Delta/\Delta}\mapright{\tau^{\ge -1}} I[1])\in Ext^1(\bar g^*\bbL_\fX,I)
\eeq
whose vanishing is necessary and sufficient for the existence of a lifting 
$$g:\Delta=\Spec B\lra \fX$$ 
such that $g|_{\bar \Delta}=\bar g$.  
Here the first morphism $\bar g^*\bbL_\fX\mapright{\bar g} \bbL_{\bar \Delta}$ is the pullback by $\bar g$; the second $\bbL_{\bar \Delta}\to \bbL_{\bar\Delta/\Delta}$ is the natural morphism from the embedding $\bar\Delta\hookrightarrow \Delta$; 
the third morphism $\tau^{\ge -1}$ is the truncation to terms of degree $\ge -1$.

If $\phi:E\to \bbL_\fX$ is a perfect obstruction theory on $\fX$, its composition with $\omega(\bar g, B,\bar B):\bar g^*\bbL_\fX\to I[1]$ gives us 
\beq\label{1.33} 
ob_\fX(\phi, \bar g, B,\bar B):\bar g^* E\lra \bar g^*\bbL_\fX \lra I[1]
\eeq
which is an element of $Ext^1(\bar g^*E,I)=I\otimes_\CC H^1(E^\vee)|_x$. 

\begin{defi}\label{1.34} Let $\phi:E\to \bbL_\fX$ and $\phi':E'\to \bbL_\fX$ be two perfect obstruction theories and $\psi:H^1(E^\vee)\to H^1({E'}^\vee)$ be an isomorphism. We say the two perfect obstruction theories $\phi$ and $\phi'$ give the \emph{same obstruction assignment} via $\psi$ if 
\[
\bar g^*(\psi)(ob_\fX(\phi,\bar g, B,\bar B))=ob_\fX(\phi',\bar g, B,\bar B)
\in I\otimes_\CC H^1({E'}^\vee)|_x
\]
\end{defi}
This explains the second condition of Definition \ref{1.35}. 
\begin{defi}\label{1.36}
A semi-perfect obstruction theory $(\phi\lalp:E\lalp\to \bbL_{X\lalp})$ is called \emph{symmetric} if $\phi\lalp$ are all symmetric (cf. Definition \ref{1.28a}) and the gluing isomorphisms $\psi_\albe$ for the obstruction sheaf $Ob_X$ are the identity maps of $\Omega_{X_\albe}$ via the canonical isomorphism $Ob_{X\lalp}\cong \Omega_{X\lalp}$.
\end{defi}

\medskip

Let us suppose $\fX$ is now a critical virtual manifold equipped with charts $(\fX_\alpha\mapright{\varphi_\alpha} V\lalp \mapright{f_\alpha}\CC)$. Since $\fX_\alpha=\zero(df_\alpha)$ in $V_\alpha$ and $df_\alpha$ is a section of the cotangent bundle $\Omega_{V_\alpha}$, we have a perfect obstruction theory (cf. \cite[\S3]{Beh})
\beq\label{1.37} \phi_\alpha:E_\alpha=[T_{V_\alpha}|_{\fX\lalp}\mapright{d(df\lalp)} \Omega_{V_\alpha}|_{\fX\lalp}]\lra \bbL_{\fX_\alpha}\eeq
where $T_{V_\alpha}=\Omega_{V_\alpha}^\vee$ is the tangent bundle of $V_\alpha$. This is symmetric because the Hessian $d(df\lalp)$ is symmetric. Moreover, $H^1(E^\vee_\alpha)=\Omega_{\fX_\alpha}$ from the exact sequence
\[ J_\alpha/J_\alpha^2\mapright{d} \Omega_{V_\alpha}|_{\fX_\alpha}\lra \Omega_{\fX_\alpha}\lra 0\]
where $J_\alpha=(df_\alpha)$ is the ideal generated by the partial derivatives of $f_\alpha$. 

The obstruction assignment for an infinitesimal lifting problem at $x\in \fX$
\beq\label{1.38} (0\to I\to B\to \bar B\to 0, \bar g:\Spec \bar B\to \fX_\alpha)\eeq
can be described as follows. We extend $\bar g:\bar \Delta=\Spec \bar B\to \fX_\alpha\hookrightarrow V_\alpha$ to a morphism $g':\Delta=\Spec B\to V_\alpha$. Since $\fX_\alpha$ is the vanishing locus of the section $df_\alpha$ of the cotangent bundle $\Omega_{V_\alpha}\to V_\alpha$, $g'$ factors through $\fX_\alpha$ if and only if $df_\alpha\circ g':\Delta\to {g'}^*\Omega_{V_\alpha}$ is zero. 
As $\bar g=g'|_{\bar \Delta}$ factors through $\fX_\alpha$, $df_\alpha\circ g'\in I\otimes \Omega_{V_\alpha}|_{x}$.

Let $\rho:I\otimes_\CC \Omega_{V_\alpha}|_x\to I\otimes_\CC \Omega_{\fX_\alpha}|_x$  be the tautological projection. 
\begin{lemm}\label{1.39}
$ob_\fX(\phi\lalp,\bar g, B,\bar B)=\rho(df_\alpha\circ g')\in I\otimes_\CC \Omega_\fX|_x$.
\end{lemm}
\begin{proof}
By \eqref{1.31}, \eqref{1.33} and \eqref{1.37}, $ob_\fX(\phi_\alpha,\bar g, B,\bar B)$ is the composition
\[ \bar g^*E\lalp \lra \bar g^*\bbL_{\fX\lalp}\lra \bbL_{\bar \Delta/\Delta}\lra I[1] \]
which fits into a commutative diagram
\beq\label{1.40}
\xymatrix{
\bar g^*\Omega_{V\lalp}\ar@{=}[r] \ar[d] & \bar g^*\Omega_{V\lalp} \ar[d] \ar[dr]^{ob_{V\lalp}}\\
\bar g^*E\lalp\ar[r]^{\phi\lalp}\ar[d] & \bar g^*\bbL_{\fX\lalp}\ar[r] & I[1]\\
\bar g^*T_{V\lalp}[1]\ar@{.>}[urr]\ar[d]\\
\bar g^*\Omega_{V\lalp}[1]
}\eeq
The first column is the distinguished triangle
\[ T_{V\lalp}|_{\fX\lalp}\lra \Omega_{V\lalp}|_{\fX\lalp} \lra [T_{V\lalp}\to \Omega_{V\lalp}]|_{\fX_\alpha}\lra T_{V\lalp}|_{\fX\lalp}[1].
\] We have an exact sequence
\beq\label{1.41}
Hom(\bar g^*\Omega_{V_\alpha},I)=I\otimes_\CC T_{V\lalp}|_x\mapright{d(df_\alpha)} Hom(\bar g^*T_{V_\alpha},I)=I\otimes_\CC \Omega_{V_\alpha}|_x\eeq
\[\lra Ext^1(\bar g^*E_\alpha,I)\lra Ext^1(\bar g^*\Omega_{V_\alpha},I).
\]

Since $V_\alpha$ is smooth, the morphism
\[ \bar\Delta=\Spec \bar B\lra \fX_\alpha\hookrightarrow V_\alpha
\]
extends to a morphism $g':\Delta\to V_\alpha$ and hence the homomorphism $ob_{V_\alpha}$ above is zero. Thus 
$ob_\fX(\phi_\alpha,\bar g, B,\bar B)$ lives in the cokernel $I\otimes_\CC \Omega_\fX|_x$ of the first homomorphism in \eqref{1.41}. The distinguished triangles
\[ \Omega_{V\lalp}|_{\fX\lalp}\lra E\lalp \lra T_{V\lalp}|_{\fX\lalp}[1],\]
\[ \bbL_{V\lalp}\lra \bbL_{\fX\lalp}\lra \bbL_{\fX\lalp/V\lalp},\]
\[ \bbL_\Delta|_{\bar \Delta}\lra \bbL_{\bar \Delta}\lra \bbL_{\bar \Delta/\Delta} \]
fit into the diagram
\beq\label{1.42}
\xymatrix{
\bar g^*\Omega_{V\lalp}|_{\fX\lalp}\ar[r]\ar@{=}[d] & \bar g^*E\lalp \ar[r]\ar[d]^{\phi\lalp} & \bar g^*T_{V\lalp}|_{\fX\lalp}[1]
\ar[dr]^{df\lalp\circ g'}\ar[d]\\
\bar g^*\bbL_{V\lalp}|_{\fX\lalp}\ar[r]\ar[d]_{g'} & \bar g^*\bbL_{\fX\lalp} \ar[r]\ar[d]^{\bar g} & \bar g^*\bbL_{\fX\lalp/V\lalp}\ar[d]\ar[r]^{\tau^{\ge -1}} & \bar g^*J\lalp/J\lalp^2[1]\ar[d]\\
\bbL_\Delta|_{\bar \Delta}\ar[r] & \bbL_{\bar \Delta} \ar[r] & \bbL_{\bar\Delta/\Delta}\ar[r]^{\tau^{\ge -1}} & I[1]
}\eeq
From \eqref{1.41} and \eqref{1.42}, we find that $ob_\fX(\phi\lalp,\bar g, B,\bar B)$ equals $\rho(df_\alpha\circ g')$.
\end{proof}

Since $E\lalp=[T_{V\lalp}|_{\fX\lalp}\mapright{d(df\lalp)}\Omega_{V\lalp}|_{\fX\lalp}]$ and $\fX\lalp=\zero(df\lalp)\subset V\lalp$, $H^1(E\lalp^\vee)=\Omega_{\fX\lalp}$ and the transition maps $\varphi_\albe:V_\albe\mapright{\cong} V_\beal$ send $f_\alpha$ to $f_\beta$ and $\Omega_{\fX_\albe}$ to $\Omega_{\fX_\beal}$ respectively. Hence we have isomorphisms
$$\psi_\albe:H^1(E_\alpha^\vee)|_{\fX_\albe}=\Omega_{\fX_\albe}\lra H^1(E_\beta^\vee)|_{X_\beal}=\Omega_{\fX_\beal}$$
induced from $\varphi_\albe$. Certainly these $\{\Omega_{\fX\lalp}, \psi_\albe\}$ glue to the cotangent sheaf $\Omega_\fX$. Moreover these local perfect obstruction theories $\{\phi_\alpha\}$ give the same obstruction assignment.
\begin{lemm}\label{1.43}
Let $f:V\to \CC$ and $g:W\to \CC$ be two equivalent LG pairs, i.e. $f$ and $g$ have only one critical value $0$ and we have a biholomorphic $\Phi:V\to W$ satisfying $g\circ \Phi=f$. 
Let $X_f=\Crit(f)$ and $X_g=\Crit(g)$ be the critical loci of $f$ and $g$ respectively.
Then under the induced isomorphism $\hat{\Phi}=\Phi|_{\fX_f}:\fX_f\mapright{\cong} \fX_g$, the perfect obstruction theories $E=[T_V|_{\fX_f}\to \Omega_V|_{\fX_f}]$ and $E'=[T_W|_{\fX_g}\to \Omega_W|_{\fX_g}]$ give the same obstruction assignment via the canonical isomorphism 
$$H^1(E^\vee)=\Omega_{X_f}\cong \hat{\Phi}^*\Omega_{X_g}=\hat{\Phi}^*H^1({E'}^\vee).$$
\end{lemm}
\begin{proof}
Because $(V,f)$ and $(W,g)$ are equivalent under $\Phi$, $\Phi$ induces an isomorphism of the critical loci 
$\hat{\Phi}: \fX_f\cong \fX_g$. Hence, we have the induced $\hat{\Phi}\sta\Omega_{\fX_g}\cong \Omega_{\fX_f}$
making the following square commutative
\beq\label{1.43d}
\begin{CD}
\Omega_V|_{\fX_f} @>{\cong}>> \hat{\Phi}\sta (\Omega_W|_{\fX_g})\\
@VV{\rho}V @VV{\varrho}V\\
\Omega_{\fX_f} @>{\cong}>> \hat{\Phi}\sta \Omega_{\fX_g},
\end{CD}
\eeq
where the vertical arrows are natural quotient homomorphisms.

We now compare the obstruction assignments. Let $\bar B=B/I$ and $\bar{\lambda}:\spec \bar B\to \fX_f$
supported at $x\in \fX_f$, as in \eqref{1.38}. We let $\bar\mu=\hat{\Phi}\circ\bar{\lambda}: \spec \bar B\to \fX_g$. We need to show that
\beq
ob_{\fX_{f}}(\phi,\bar{\lambda},B,\bar B)=ob_{\fX_{g}}(\phi',\bar\mu, B,\bar B)\in I\otimes_\CC \Omega_{\fX_f}|_x=I\otimes_\CC \Omega_{\fX_g}|_{\hat{\Phi}(x)}.
\eeq

We extend $\bar{\lambda}$ to $\lambda':\spec B\to V$, and let $\mu'=\Phi\circ \lambda': \spec B\to W$.
Using $df\circ \bar{\lambda}=0$ and $I\cdot\bm_B=0$, $df\circ \lambda'\in I\otimes_\CC\Omega_V|_x$.
By a similar reason, $dg\circ\mu'\in I\otimes_\CC\Omega_W|_{\hat{\Phi}(x)}$. From $\Phi\sta(dg)=df$ and $\mu'=\Phi\circ \lambda'$, we 
conclude
$$df\circ \lambda'=\Phi\sta(dg\circ\mu')\in I\otimes_\CC \Omega_V|_x\cong I\otimes_\CC \Omega_W|_{\hat{\Phi}(x)}.
$$
From \eqref{1.43d}, we obtain 
$$\rho(df\circ \lambda')=\varrho(d g\circ\mu')\in I\otimes_\CC\Omega_{\fX_f}|_x\cong I\otimes_\CC\Omega_{\fX_g}|_{\hat{\Phi}(x)}.
$$
Thus $ob_{\fX_{f}}(\phi,\bar{\lambda},B,\bar B)=ob_{\fX_{g}}(\phi',\bar \mu,B,\bar B)$, which proves the lemma.
\end{proof}

So we proved the following.
\begin{prop}\label{1.44}
A critical virtual manifold $\fX$ with charts $(\fX_\alpha\mapright{\varphi\lalp} V\lalp \mapright{f\lalp} \CC)$ admits a symmetric semi-perfect obstruction theory $\{E_\alpha=[T_{V\lalp}|_{\fX\lalp}\to \Omega_{V\lalp}|_{\fX\lalp}]\}$ whose obstruction sheaf $Ob_\fX$ is the cotangent sheaf $\Omega_\fX$ of $\fX$. 
\end{prop}
\begin{rema}\label{1.45h}
In \cite[Definition-Theorem 4.7]{ChLiS}, moduli spaces of simple gluable objects with fixed determinant in the derived category $D^b(Y)$ of coherent sheaves on a Calabi-Yau 3-fold $Y$ are proven to have semi-perfect obstruction theories. Here an object $E\in D^b(Y)$ is simple if $Hom_{D^b(Y)}(E,E)=\CC$ and gluable if $Ext^{<0}(E,E)=0$. It is expected that in fact these moduli spaces are critical virtual manifolds. In Part 2, we will see that this is indeed true for moduli spaces of simple sheaves. 
\end{rema}

\begin{rema}\label{1.45}
The isomorphism $\varphi_\albe$ in Definition \ref{1.5} induces the isomorphism
\beq\label{1.46} 
\xymatrix{
E_\alpha|_{\fX_\albe} \ar[d] \ar@{=}[r] & [T_{V\lalp}|_{\fX_\albe}\ar[d]\ar[r]^{d(df\lalp)} &\Omega_{V\lalp}|_{\fX_\albe}]\ar[d]\\
E_\beta|_{\fX_\albe}  \ar@{=}[r] & [T_{V_\beta}|_{\fX_\albe}\ar[r]^{d(df_\beta)} &\Omega_{V_\beta}|_{\fX_\albe}]
}\eeq
whose vertical maps are isomorphisms induced from $\varphi_\albe$ since $f_\beta\circ\varphi_\albe=f\lalp$. However, because $\varphi_\albega$ in \eqref{mapcoc} need \emph{not} be the identity morphism, there is no obvious reason for the local perfect obstruction theories $\{E\lalp\}$ to glue to a global perfect obstruction theory. 
\end{rema}

\bigskip

\subsection{Virtual cycle and Donaldson-Thomas type invariant}\label{S1.4.3}
In \cite{Beh}, Behrend proved that when a scheme admits a symmetric perfect obstruction theory, the Donaldson-Thomas type invariant, which is the degree of the virtual cycle, is the weighted Euler characteristic (cf. \S\ref{S1.4.1}). We extend this result to the case of symmetric semi-perfect obstruction theory and apply it to critical virtual manifolds.

\medskip

Let $X$ be an analytic space equipped with a symmetric semi-perfect obstruction theory $\{\phi\lalp:E\lalp\to \bbL_{X\lalp}\}$. By \cite[(3.5)]{ChLiS}, we then have the obstruction cone cycle $C_\fX\in Z_*(Ob_\fX)$ where $Ob_X$ denotes the obstruction sheaf of $X$ defined in Definition \ref{1.35}. 
By applying the Gysin map $s^!$ for the sheaf stack $Ob_\fX$ when $\fX$ is proper (cf. \cite[Proposition 3.4]{ChLiS}), we obtain the virtual cycle of $\fX$ defined by
\beq\label{1.47}
[\fX]\virt:=s^![C_\fX]\in A_0(\fX).\eeq

In case $\fX$ admits a global embedding $\imath:\fX\hookrightarrow \PP$ into a complex manifold $\PP$, we can give an alternative description of $[\fX]\virt$ by using the following proposition which is an immediate consequence of \cite[Proposition 2.2]{Beh}. 
\begin{prop}\label{1.48} 
Let $X$ be an analytic space equipped with a symmetric semi-perfect obstruction theory $\{\phi\lalp:E\lalp\to \bbL_{X\lalp}\}$. 
%There exists a unique closed subcone $C_\PP\subset \Omega_\PP$ such that for each $\alpha$ and every lift $\phi$
%\[\xymatrix{
%&\Omega_{V\lalp}|_{\fX\lalp}\ar[d]^{\varphi_\alpha^*}\\
%\Omega_\PP|_{\fX_\alpha}\ar[r]^{\imath^*}\ar[ur]^\phi & \Omega_{\fX\lalp}
%}\]
There exists a unique closed subcone $C_\PP\subset \Omega_\PP$ such that for each local resolution 
$F\to E\lalp^\vee[1]|_{U}$ for open $U\subset X\lalp$ by a locally free sheaf $F$ on $U$ and every lift $\phi$
\[\xymatrix{
&F\ar[d]\\
\Omega_\PP|_{U}\ar[r]^{\imath^*}\ar[ur]^\phi & \Omega_{U}
}\]
%we have $C_\PP|_{\fX\lalp}=\phi^{-1}(C_\alpha)$ where $C\lalp$ is the fiber product of the intrinsic normal cone 
%$$\fC_{\fX\lalp}\hookrightarrow h^1/h^0(E\lalp^\vee)=(\Omega_{V\lalp}|_{\fX\lalp})/(T_{V\lalp}|_{\fX\lalp})$$
%and the quotient morphism $\Omega_{V\lalp}|_{\fX\lalp}\to h^1/h^0(E^\vee_\alpha)$. 
we have $C_\PP|_{U}=\phi^{-1}(C_F)$ where $C_F$ is the fiber product of the intrinsic normal cone 
$$\fC_{\fX\lalp}\hookrightarrow h^1/h^0(E\lalp^\vee)$$
and the quotient morphism $F\to h^1/h^0(E^\vee_\alpha)|_U$. 
\end{prop}

Indeed, by \cite[Proposition 2.2]{Beh}, there are unique closed subcones $C_\PP|_{X\lalp}\subset \Omega_\PP|_{X\lalp}$ satisfying the above lifting property. By the uniqueness over the intersections $X_\albe$, these subcones uniquely glue to a cone $C_\PP\subset \Omega_\PP$. See \cite{BeFa} for the theory of intrinsic normal cones. 

Note that $$C_\PP\in Z_{\dim \PP}(\Omega_\PP)$$ 
since $\fC_{\fX\lalp}$ is 0-dimensional.

\begin{defi}\label{1.49} 
Suppose $\fX$ is a compact analytic space equipped with a symmetric semi-perfect obstruction theory and
%compact critical virtual manifold 
with an embedding into a complex manifold $\PP$. The \emph{virtual fundamental class} of $\fX$ is defined to be
\beq\label{1.50} 
[\fX]\virt=0^!_{\Omega_\PP}[C_\PP]\in A_0(\fX)
\eeq
using the cone in Proposition \ref{1.48}. We define the Donaldson-Thomas type invariant of $X$ to be 
\beq\label{1.51}
DT(\fX)=\deg [\fX]\virt
\eeq
where $\deg:A_0(\fX)\to A_0(\mathrm{pt})=\ZZ$ is the pushforward by the constant map $X\to \mathrm{pt}=\Spec\, \CC$.
\end{defi}

\begin{lemm}\label{1.50l}
$[\fX]\virt$ is independent of the embedding $\fX\hookrightarrow \PP$ into a complex manifold.
\end{lemm}
\begin{proof}
Given two embeddings $\imath, \imath'$ of $\fX$ into $\PP$ and $\PP'$, we have an embedding 
$$(\imath,\imath'):\fX\hookrightarrow \PP\times\PP'$$
and exact sequences 
\[ 0\lra p_2^*\Omega_{\PP'}|_\fX \lra C_{\PP\times \PP'} \lra C_\PP\lra 0 \]
\[ 0\lra p_1^*\Omega_{\PP}|_\fX \lra C_{\PP\times \PP'} \lra C_\PP'\lra 0 \]
where $p_1, p_2$ are the projections from $\PP\times \PP'$ to $\PP$ and $\PP'$ respectively. 
Hence $$0^!_{\Omega_\PP}[C_\PP]=0^!_{\Omega_{\PP\times\PP'}}[C_{\PP\times\PP'}]=0^!_{\Omega_{\PP'}}[C_{\PP'}]$$ as desired.
\end{proof}

By the construction of $s^!$ in \cite[Proposition 3.4]{ChLiS}, it is obvious that the virtual fundamental class in Definition \ref{1.49} coincides with \eqref{1.47}.

\medskip

We next show that the Donaldson-Thomas type invariant \eqref{1.51} is the Euler characteristic \eqref{1.23} weighted by the Behrend function $\nu_\fX$.
\begin{lemm}\label{1.52}
Let $X$ be a compact analytic space equipped with a symmetric semi-perfect obstruction theory.
Fix an embedding $\imath:\fX\hookrightarrow \PP$ of $\fX$ into a complex manifold $\PP$. Then the cone $C_\PP$ from Proposition \ref{1.48} is conic Lagrangian, and $C_\PP=\ell(\fc_\fX)$ where $\fc_\fX$ is the characteristic cycle of $\fX$ in \eqref{1.19}. 
\end{lemm}
\begin{proof}
Being conic Lagrangian is a local property. Hence $C_\PP$ is conic Lagrangian by \cite[Theorem 4.14]{Beh}. For $C_\PP=\ell(\fc_\fX)$, use Lemma \ref{1.24} and the fact that $C_\PP|_{\fX\lalp}$ is the pullback of intrinsic normal cone $\fC_{\fX\lalp}$ for each $\alpha$.
\end{proof}

Combining Lemma \ref{1.52} and Proposition \ref{1.25}, we obtain the following generalization of \cite[Theorem 4.18]{Beh} to the setting of semi-perfect obstruction theory.
\begin{theo}\label{1.53}
Let $\fX$ be a compact analytic space equipped with a symmetric semi-perfect obstruction theory and a closed immersion into a complex manifold $\PP$. Then 
\beq\label{1.54} DT(\fX)=\chi_\nu(\fX)\eeq
where $DT(\fX)=\deg [\fX]\virt$ and $\chi_\nu(\fX)=\sum n\cdot \chi_c(\nu^{-1}(n))$ is the Euler characteristic of $\fX$ weighted by the Behrend function defined by \eqref{1.21}.
\end{theo}

\begin{coro}\label{1.57}
The Donaldson-Thomas type invariant $DT(\fX)=\deg [\fX]\virt$ depends only on the analytic space $(X,\sO_X)$, i.e. $DT(\fX)$ is independent of the symmetric semi-perfect obstruction theory or the embedding into a complex manifold $\PP$. 
\end{coro}
\begin{proof}
By Remark \ref{1.58}, the weighted Euler characteristic $\chi_\nu(\fX)$ depends only on the analytic space $(X,\sO_X)$. The corollary follows from Theorem \ref{1.53}. 
\end{proof}

Let $\fX$ be a critical virtual manifold with charts $(\fX_\alpha\mapright{\varphi\lalp} V\lalp \mapright{f\lalp} \CC)$. Let 
$$\{E_\alpha=[T_{V\lalp}|_{\fX\lalp}\mapright{d(df\lalp)}\Omega_{V\lalp}|_{\fX\lalp}]\}$$ 
be the symmetric semi-perfect obstruction theory from Proposition \ref{1.44} so that the cotangent sheaf $\Omega_\fX$ is the gluing of local obstruction sheaves $\{H^1(E\lalp^\vee)\}$. By Corollary \ref{1.57}, we therefore obtain the following.
\begin{coro}\label{1.59}
Let $\fX$ be a compact critical virtual manifold equipped with an embedding into a complex manifold $\PP$. Then $DT(\fX)=\chi_\nu(\fX)$ and the invariant depends only on the analytic space underlying the critical virtual manifold. 
\end{coro}

\begin{rema}\label{1.56}
In \cite{ChLiS}, the theory of semi-perfect obstruction was developed in the more general setting of a separated proper Deligne-Mumford stack $\fX$ over a smooth Artin stack $\cM$ of pure dimension. The Donaldson-Thomas type invariant $DT(\fX)=\deg [\fX]\virt$ is deformation invariant (cf. \cite[Proposition 3.8]{ChLiS}): If
\[\xymatrix{
\fX'\ar[r]^u\ar[d] & \fX\ar[d]\\
\cM'\ar[r]^v & \cM
}\]
is a fiber product and $v$ is a regular embedding of smooth Artin stacks of pure dimensions,
\[ v^![\fX]\virt=[\fX']\virt .\]
\end{rema}

\bigskip

\subsection{A categorification problem}\label{S1.4.4}
Categorification roughly means finding a deeper structure underlying a known invariant. 
The categorification problem for critical virtual manifolds is the following.
\begin{prob}\label{1.55} 
Find a cohomology theory $\cH^*(\fX)$ for a compact critical virtual manifold $\fX$ whose Euler characteristic 
\[ \sum_m (-1)^m \dim \cH^m(\fX) \]
is equal to the Donaldson-Thomas type invariant
\[ DT(\fX)=\deg [\fX]\virt =\chi_\nu(\fX) .\]
\end{prob}
Originally this problem was proposed by Joyce and Song in \cite[Question 5.7]{JoSo} for the case where $X$ is a moduli space of simple coherent sheaves on a Calabi-Yau 3-fold $Y$. 

In subsequent chapters, we will construct such a cohomology theory as the hypercohomology of a perverse sheaf on $\fX$ which is locally the perverse sheaf $P_\alpha$ of vanishing cycles for $f_\alpha:V\lalp\to \CC$ whenever $\fX$ is an \emph{orientable} critical virtual manifold.

\bigskip

\section{Critical virtual manifolds and $d$-critical loci}\label{S1.5}
\def\cS{\mathcal{S} }
In \cite{Joyce}, Joyce introduced the notion of  a $d$-critical locus as a classical model for a $(-1)$-shifted symplectic derived scheme. In this section, we prove that critical virtual manifolds are $d$-critical loci and vice versa.

\medskip

According to \cite[Definition 2.3]{Joyce}, a $d$-critical locus is a pair $(X,s)$ of an analytic space $X$ and a section $s$ of a sheaf $\cS_X^0$ satisfying the condition that for any $x\in X$ there exists an open neighborhood $U$ of $x$ which is isomorphic to the critcal locus $\Crit(f)$ of a holomorphic function $f$ on a complex manifold $V$ such that $f$ defines the section $s|_U$. Roughly speaking, $\cS_X^0$ sheafifies the choices of charts $(V,f)$ and their compatibility (cf. Definition \ref{1.63}) while allowing the dimensions of the charts to vary.  The definition of the sheaf $\cS_X^0$ is rather involved (cf. \cite[Theorem 2.1]{Joyce}). We will not reproduce the precise definition here but instead we will use the following.
\begin{prop}\label{1.70}
An analytic space $X$ is a $d$-critical locus if and only if there is an open cover $X=\cup\lalp X\lalp$ and LG pairs $\{(V\lalp,f\lalp)\}$ such that
\begin{enumerate}
\item for each $\alpha$, $\Crit(f\lalp)\cong X\lalp$ as analytic spaces;
\item for each pair $(\alpha,\beta)$ of indices, there exist open neighborhoods $V_\albe$ (resp. $V_\beal$) of $X_\albe=X\lalp\cap X_\beta$ in $V\lalp$ (resp. $V_\beta$) and closed embeddings
\beq\label{1.71}
\xymatrix{
V_\albe\ar@{^(->}[r]^{\imath_\albe} & W_\albe & V_\beal\ar@{_(->}[l]_{\imath_\beal}
}\eeq
into a complex manifold $W_\albe$ satisfying
\beq\label{1.72}
g_\albe\circ\imath_\albe=f\lalp|_{V_\albe},\quad g_\albe\circ \imath_\beal=f_\beta|_{V_\beal}\eeq
for some holomorphic function $g_\albe$ on $W_\albe$ whose critical locus $\Crit(g_\albe)$ is $X_\albe$.
\end{enumerate}
\end{prop}
\begin{proof}
If we have an open cover $X=\cup\lalp X\lalp$ and LG pairs $\{(V\lalp,f\lalp)\}$ satisfying the stated conditions, then we certainly have a $d$-critical locus structure on $X$ because the local sections $s\lalp$ defined by $f\lalp$ glue to a globally defined section $s$ by \cite[Theorem 2.1(ii)]{Joyce}.

Conversely, suppose we have a $d$-critical locus $(X,s)$ and at every point $x\in X$ we can find a neighborhood $U$ of $x$ in $X$ which is the critical locus of an LG pair $(V,f)$ such that $f$ defines the section $s|_U$. Hence there is an open cover $X=\cup_{\lambda\in \Lambda} X_\lambda$ and LG pairs $(V_\lambda,f_\lambda)$ such that $\Crit(f_\lambda)\cong X_\lambda$. 
Since $X$ is paracompact, we may assume the covering is locally finite. For each $x\in X$, 
$$\Lambda_x=\{\lambda\in \Lambda\,|\,x\in X_\lambda\}$$
is a finite set. Since a second countable paracompact space is a metric space, we can choose a metric $d$ on $X$. By \cite[Theorem 2.20]{Joyce}, we can pick $\vareps_x>0$ such that 
$$B(x,3\vareps_x)\sub \cap_{\lambda\in \Lambda_x}X_\lambda$$
and that for any $\lambda_1,\lambda_2\in \Lambda_x$, there is an LG pair $(W_{\lambda_1\lambda_2},g_{\lambda_1\lambda_2})$ satisfying 
$$\Crit(g_{\lambda_1\lambda_2})=B(x,3\vareps_x)$$
together with closed embeddings
$$\xymatrix{
V_{\lambda_1}^\circ\ar@{^(->}[r]^{\imath_1} & W_{\lambda_1\lambda_2} & V_{\lambda_2}^\circ \ar@{_(->}[l]_{\imath_2}
}$$
of open neighborhoods $V_{\lambda_i}^\circ$ of $B(x,3\varepsilon_x)$ in $V_{\lambda_i}$ for $i=1,2$, satisfying 
$$g_{\lambda_1\lambda_2}\circ\imath_1=f_{\lambda_1}|_{V_{\lambda_1}^\circ},\quad g_{\lambda_1\lambda_2}\circ\imath_2=f_{\lambda_2}|_{V_{\lambda_2}^\circ}.$$

Fix a map $\lambda:X\to \Lambda$ such that $\lambda(x)\in \Lambda_x$ for all $x$. Then we can find an open neighborhood $V_x$ of 
$$U_x=B(x,\vareps_x)$$
in $V_{\lambda(x)}$ such that the critical locus of $f_x=f_{\lambda(x)}|_{V_x}$ is $U_x$. We claim that 
$$X=\cup_{x\in X}U_x, \quad (V_x,f_x)$$
satisfy the conditions in Proposition \ref{1.70}. We have to check a diagram similar to \eqref{1.71} for $\{V_x\}$. Indeed, if $U_x\cap U_{x'}\ne \emptyset$ with $\vareps_x\ge \vareps_{x'}$, then 
$$U_x\cup U_{x'}\sub B(x,3\vareps_x)$$
and thus $\Lambda_x\sub\Lambda_{x'}$. In particular, both $\lambda(x)$ and $\lambda(x')$ lie in $\Lambda_{x'}$. Now by assumption, we have closed embeddings
$$\xymatrix{V_{\lambda(x)}^\circ\ar@{^(->}[r]^{\imath} & W_{\lambda(x)\lambda(x')} & V_{\lambda(x')}^\circ \ar@{_(->}[l]_{\imath'}
}$$ 
of open neighborhoods $V_{\lambda(x)}^\circ$, $V_{\lambda(x')}^\circ$ of $U_{x'}$  in $V_{\lambda(x)}$, $V_{\lambda(x')}$ respectively, for an LG pair $(W_{\lambda(x)\lambda(x')},g)$ with $\Crit(g)=U_{x'}$ and
$$g\circ\imath=f_{x}|_{V_{\lambda(x)}^\circ}\and 
g\circ\imath'=f_{x'}|_{V_{\lambda(x')}^\circ}.$$
%over $U_{x'}=B(x,\vareps_{x'})$, which 
The diagram certainly restricts to such a diagram for the open subset $U_x\cap U_{x'}$. This proves the proposition.
\end{proof}

\begin{defi}\label{1.76}
A $d$-critical locus $(X,s)$ is called \emph{finite dimensional} if there are an open cover $X=\cup\lalp X\lalp$ and LG pairs $\{(V\lalp,f\lalp)\}$ such that $\{\dim V\lalp, \dim W_\albe\}$ is bounded, in the notation of Proposition \ref{1.70}.
\end{defi}
\begin{prop}\label{1.75}
A critical virtual manifold is a finite dimensional $d$-critical locus. Conversely a finite dimensional $d$-critical locus  is a critical virtual manifold.
\end{prop}
%\begin{rema}\label{1.76}
%The compactness in the second statement is used only to guarantee the existence of a uniform $r$ which is greater than or equal to the dimension of a local chart everywhere.
%\end{rema}
\begin{proof}%[Proof of Proposition \ref{1.75}]
The first statement is obvious from Proposition \ref{1.70} because if
$$X=\cup\lalp X\lalp, \quad X\lalp=\Crit(f\lalp)\sub V\lalp\mapright{f\lalp}\CC$$
is a critical virtual manifold structure, then the biholomorphic map
$\varphi_\albe:V_\albe\to V_\beal$ gives a diagram
$$\xymatrix{
V_\albe\ar@{^(->}[r]^{\imath_\albe} & W_\albe & V_\beal \ar@{_(->}[l]_{\imath_\beal}
}$$
with $W_\albe=V_\albe$, $\imath_\albe=\id$ and $\imath_\beal=\varphi_\beal$. 

Conversely, if $X$ is a finite dimensional $d$-critical locus, then we have an open cover $X=\cup X\lalp$, charts
$$\{(X\lalp=\Crit(f\lalp)\sub V\lalp \mapright{f\lalp}\CC)\}$$
and diagrams \eqref{1.71} satisfying \eqref{1.72}.
Moreover, we can find an integer $r$ such that $$r\ge \dim W_\albe,\quad r\ge \dim V\lalp,\quad \forall \alpha,\beta.$$
By Lemma \ref{1.60}, for any $x\in X$, we have a commutative diagram
$$\xymatrix{
V_\albe^\circ\times \CC^{r-\dim V_{\alpha}} \ar[r]^{\cong}\ar[dr]_{\tilde f\lalp} &
W_\albe^\circ\times \CC^{r-\dim W_\albe}\ar[d]^{\tilde{g}_\albe} & V_\beal^\circ\times\CC^{r-\dim V_\beta}\ar[l]_{\cong}\ar[dl]^{\tilde{f}_\beta}\\
& \CC
}$$
for open neighborhoods $V_\albe^\circ$, $V_\beal^\circ$ and $W_\albe^\circ$ of $x$ in $V_\albe$, $V_\beal$ and $W_\albe$ respectively.
Here $\tilde{f}\lalp=f\lalp+\sum_{i=1}^{r-\dim V_\alpha}z_i^2$, 
$\tilde{f}_\beta=f_\beta+\sum_{i=1}^{r-\dim V_\beta}z_i^2$,
$\tilde{g}_\albe=g_\albe+\sum_{i=1}^{r-\dim W_\albe}z_i^2$.
So we have new charts
$$(\tilde{V}\lalp=V\lalp\times\CC^{r-\dim V\lalp},
\tilde{f}\lalp=f\lalp+\sum_{i=1}^{r-\dim V\lalp}z_i^2)$$
for $X$ and for any $x\in X\lalp\cap X_\beta$, we have a biholomorphic map
$$\tilde{\varphi}_\albe:\tilde{V}\lalp^\circ\lra \tilde{V}_\beta^\circ$$
of open neighborhoods of $x$ such that $\tilde{f}_\beta\circ \tilde{\varphi}_\albe=\tilde{f}_\alpha$. Then the theorem follows from the same argument as in the proof of Proposition \ref{1.70} or more explicitly from Proposition \ref{4.44}.
%After choosing a metric on $X$, we can apply the same argument as in the proof of Proposition \ref{1.70} to deduce that $X$ is a critical virtual manifold with $r$-dimensional charts. We omit the detail to avoid repetition.
\end{proof}

%%%%%%%%%%%%%
%%%%%%%%%%%%%
%%%%%%%%%%%%%
%%%%%%%%%%%%%
%%%%%%%%%%%%%
\def\Re{\mathrm{Re}\,}
\chapter{Perverse sheaves on critical virtual manifolds}\label{ch2}

In this chapter, we provide a solution to the categorification problem (Problem \ref{1.55}).  We prove that when $X$ is an orientable critical virtual manifold with charts $(X\lalp\mapright{\varphi\lalp} V\lalp\mapright{f\lalp}\CC)$, there is a perverse sheaf $P\in Perv(\QQ_X)$ which underlies a polarizable mixed Hodge module $\cM\in MHM(X)^p$ such that the hypercohomology $\bbH^*_c(X,P)$ has Euler characteristic equal to the Euler characteristic $\chi_\nu(X)$ of $X$ weighted by the Behrend function (cf. Theorems \ref{2.23} and \ref{2.32}). 
Our perverse sheaf $P$ on $X$ is the gluing of the perverse sheaves $P_\alpha$ of vanishing cycles for $f\lalp$ and so is the mixed Hodge module $\cM$. The obstruction for gluing $\{P\lalp\}$ is shown to coincide with that for the orientability of $X$.
Moreover, we construct a perverse sheaf $\hat{P}$ with $\chi_c(X,\hat{P})=\chi_\nu(X)$, which underlies a direct sum of polarizable Hodge modules (cf. Corollary \ref{2.33}). The hypercohomology $\bbH^*_c(X,\hat{P})$ will be of fundamental use for a mathematical theory of the Gopakumar-Vafa invariant in Part 2. 

\bigskip

\section{Perverse sheaves of vanishing cycles}\label{S2.1}
In this section, we recall the perverse sheaf of vanishing cycles for the constant sheaf on a complex manifold $V$ and a holomorphic function $f$ on $V$. We show that if $X$ is a critical virtual manifold with charts $(X\lalp\mapright{\varphi\lalp}V\lalp\mapright{f\lalp}\CC)$, there are sheaf complexes $P_\alpha=\varphi_\alpha^*P_{f\lalp}\in D^b_c(\QQ_{X_\alpha})$ and isomorphisms $\sigma_\albe:P\lalp|_{X_\albe}\to P_\beta|_{X_\albe}$ (cf. Proposition \ref{2.7} and Corollary \ref{2.11}). 

\medskip

Recall that an LG pair $(V,f)$ refers to a complex manifold $V$ together with a holomorphic function $f$ on $V$ which has only one critical value $0$. 
Let $$V_{>0}=\{x\in V\,|\, \mathrm{Re}\, f(x)>0\} \and 
V_{\le 0}=\{x\in V\,|\, \mathrm{Re}\, f(x)\le 0\}.$$
For any sheaf $\cF$ of $\QQ$-vector spaces on $V$, let
$\Gamma_{V_{\le 0}}\cF$ denote the sheaf 
$$U\mapsto \ker[\Gamma(U,\cF)\to \Gamma(U\cap V_{>0},\cF) ],\quad  \forall \, U\text{ open}.$$
Let $R\Gamma_{V_{\le 0}}:D^b_c(\QQ_V)\to D^b_c(\QQ_V)$ be the right derived functor of the left exact functor $\Gamma_{V_{\le 0}}$ where $D^b_c(\QQ_V)$ denotes the derived category of bounded (cohomologically) constructible complexes of sheaves of $\QQ$-vector spaces on $V$. 

\begin{defi}\label{2.1}
The \emph{perverse sheaf of vanishing cycles} for an LG pair $(V,f)$ is defined as
\beq\label{2.2}
P_f=R\Gamma_{V_{\le 0}}\QQ_V[\dim V]|_{f^{-1}(0)}\in D^b_c(\QQ_{f^{-1}(0)}). \eeq
\end{defi}
Here $[\dim V]$ denotes the translation functor, shifting a complex to the left by $\dim V$. 

\begin{exam}  \label{2.9} Let $q=\sum_{i=1}^ry_i^2:V=\CC^r\to\CC$. The set $V_{>0}\subset V=\CC^r$  is a disk bundle over $\RR^r-\{0\}$ which is homotopic to $S^{r-1}$. From the distinguished triangle 
\beq\label{edtpv} R\Gamma_{V_{\le 0}}\QQ\lra \QQ\lra R\imath_*\imath^*\QQ\mapright{+1},
\eeq 
where $\imath:V_{>0}\hookrightarrow V$ is the inclusion, we find that 
$P_q[1-r]= R\Gamma_{V_{\le 0}}\QQ[1]$ is a sheaf complex supported at the origin $0$ satisfying $P_q[1-r]|_0\cong\QQ[-r+1]$, i.e. $P_q\cong\QQ_0.$
\end{exam}

By \cite[Exercise VIII.13]{KaSh} or \cite[Proposition 4.2.9]{Dim}, \eqref{2.2} is equivalent to the following definition: Let $\widetilde{\CC}^\times\to \CC^\times$ denote the universal cover of $\CC^\times=\CC-\{0\}$ and let $\widetilde{V}^\times$ be the fiber product of $V^\times=V-f^{-1}(0)\mapright{f}\CC^\times$ with the universal cover $\widetilde{\CC}^\times\to \CC^\times$ so that we have the diagram
\beq\label{2.2a} \xymatrix{
&& \widetilde{V}^\times \ar[r]\ar[d]^\pi \ar[dl]_{\bar\pi} & \widetilde{\CC}^\times\ar[d]\\
f^{-1}(0)\ar@{^(->}[r]^\imath & V & V^\times\ar[r]\ar@{_(->}[l]_\jmath & \CC^\times.
}\eeq
Then $P_f[1-\dim V]$ is isomorphic to the cone of the morphism
\beq\label{2.3}
\QQ_{f^{-1}(0)}=\imath^*\QQ_V\lra \imath^*R\jmath_*R\pi_*\pi^*\jmath^*\QQ_V=\imath^*R\jmath_*R\pi_*\QQ_{\widetilde{V}^\times}
=\imath^*R\bar\pi_*\QQ_{\widetilde{V}^\times}
\eeq
in the triangulated category $D^b_c(\QQ_{f^{-1}(0)})$, given by the adjunctions $\id\to R\pi_*\pi^*$ and $\id\to R\jmath_*\jmath^*$. 

By \cite[Example 4.2.6]{Dim}, the stalk cohomology of $P_f$ at $x\in f^{-1}(0)$ is the reduced cohomology
\beq\label{2.4} \cH^k(P_f)_x=\tilde{H}^{k+\dim V-1}(MF_x,\QQ)  \eeq
of the Milnor fiber 
\beq\label{2.4a} 
MF_x=f^{-1}(\delta)\cap B_\epsilon(x) \quad\text{for  }0<\delta <\!< \epsilon<\!< 1.\eeq 
Since $f$ is submersive away from $\fX_f=\Crit(f)$, the Milnor fibers $MF_x$ for $x\in f^{-1}(0)-\fX_f$ are contractible and $P_f$ is trivial on $f^{-1}(0)-\fX_f$. Hence we have 
$$P_f\in D^b_c(\QQ_{\fX_f}).$$
Moreover, by \eqref{2.4}, we have
\beq\label{2.5} \chi(P_f)_x=\sum_m (-1)^m \dim \cH^m(P_f)_x=(-1)^{\dim V}(1-\chi(MF_x)). \eeq
By \eqref{1.22}, we find that the right side of \eqref{2.5} equals the value $\nu_{\fX_f}(x)$ of the Behrend function for the analytic space $\fX_f=\Crit(f)$ defined by the partial derivatives of $f$. Therefore we have
\beq\label{2.5a} \chi(P_f)_x=\nu_{\fX_f}(x).\eeq

Given any bounded constructible complex of sheaves $P$ on a topological space $X$ and open $U\subset X$, we have a distinguished triangle
\[R\jmath_!\jmath^{-1}P\lra P\lra R\imath_*\imath^{-1}P\lra \]
where $\imath:X-U\to X$ and $\jmath:U\to X$ are inclusion maps. Taking the compact support hypercohomology, we obtain an exact sequence
\[\cdots\lra \bbH^k_c(U,P)\lra \bbH^k_c(X,P)\lra \bbH^k_c(X-U,P)\lra \bbH_c^{k+1}(U,P)\lra\cdots\]
So we have the additive property of the Euler characteristic:
\beq\label{2.5b} \chi_c(X,P)=\chi_c(U,P)+\chi_c(X-U,P)\eeq
where $\chi_c(-,P)=\sum_m(-1)^m\dim \bbH^m_c(-,P)$. 

By using the usual spectral sequence from sheaf cohomology to hypercohomology
\[ E_2^{p,q}=H^p_c(X,\cH^q(P))\Rightarrow \bbH^{p+q}_c(X,P) \]
together with \eqref{2.5a} and \eqref{2.5b}, we find that the Euler characteristic of the hypercohomology $\bbH_c^*(\fX_f,P_f)$ equals the Euler characteristic
\[ \chi_\nu(\fX_f) \]
weighted by the Behrend functin $\nu=\nu_{\fX_f}$, since $P_f$ is (cohomologically) constructible.
By the same argument, we obtain the following solution to Problem \ref{1.55}.

\begin{prop}\label{2.6}
Let $\fX$ be a critical virtual manifold with charts $$(\fX_\alpha\mapright{\varphi_\alpha} V\lalp \mapright{f\lalp}\CC)$$ and let $P\lalp=\varphi^*\lalp P_{f\lalp}$ denote the pullback of the perverse sheaf $P_{f\lalp}$ of vanishing cycles for $f\lalp:V\lalp\to \CC$ via $\fX\lalp\cong \varphi\lalp(\fX\lalp)=\Crit(f\lalp)=\fX_{f\lalp}$. Suppose there are an object $P\in D^b_c(\QQ_{\fX})$ and isomorphisms $P|_{\fX\lalp}\cong P\lalp$ in $D^b_c(\QQ_{\fX\lalp})$ for all $\alpha$. Then the Euler characteristic of the compact support hypercohomology $\bbH_c^*(\fX,P)$ equals the Euler characteristic $\chi_\nu(\fX)$ weighted by the Behrend function \eqref{1.23}. 
\end{prop}
In other words, our categorification problem (Problem \ref{1.55}) will be solved if we can glue $\{P\lalp\}$ to a globally defined sheaf complex $P$.  In the subsequent sections we will construct canonical gluing isomorphisms of $\{P\lalp\}$, called of geometric origin, and show that these isomorphisms enable us to glue $\{P\lalp\}$ to a perverse sheaf $P$ on $X$ when
the critical virtual manifold $\fX$ is orientable. Moreover $P$ is unique up to twisting by a $\ZZ_2$-local system in $H^1(\fX,\ZZ_2)$. 

For gluing $\{P\lalp\}$, we will use the following.
\begin{prop}\label{2.7} (cf. \cite[Proposition 4.2.11]{Dim} or \cite[Exercise VIII.15]{KaSh})
Let $(V,f)$ and $(W,g)$ be two LG pairs (Definition \ref{1.4}). Let $\Phi:V\to W$ be a homeomorphism satisfying $g\circ \Phi=f$. Then 
\beq\label{2.8} R\hat{\Phi}_*(P_f)\cong P_g\eeq
where $\hat{\Phi}:f^{-1}(0)\to g^{-1}(0)$ is the restriction of $\Phi$. Since $R\hat{\Phi}_*$ and $\hat{\Phi}^{-1}=\hat{\Phi}^*$ are adjoints, we also have an isomorphism 
\beq\label{2.8a}P_f\cong \hat{\Phi}^*P_g.\eeq
\end{prop}
\begin{proof}
Using the notation of \eqref{2.2a}, consider the diagram
\[\xymatrix{
f^{-1}(0)\ar@{^(->}[r]^{\imath_V} \ar[d]_{\hat{\Phi}} & V\ar[d]^\Phi & \widetilde{V}^\times\ar[d]^{\tilde{\Phi}}\ar[r]\ar[l]_{\bar\pi_V} & \widetilde{\CC}^\times\ar@{=}[d]\\
g^{-1}(0)\ar@{^(->}[r]^{\imath_W} & W & \widetilde{W}^\times\ar[l]_{\bar\pi_W}\ar[r] & \widetilde{\CC}^\times
}\]
where the vertical maps except the last are the homeomorphisms induced from $\Phi$.
Then we have 
\[ R\hat{\Phi}_*[\imath_V^*R(\bar\pi_V)_*\QQ_{\widetilde{V}^\times}] = \imath_W^*R\Phi_*R(\bar\pi_V)_*\QQ_{\widetilde{V}^\times} \]\[
=\imath_W^*R(\bar\pi_W)_*(R\tilde{\Phi}_*\QQ_{\widetilde{V}^\times})
=\imath_W^*R(\bar\pi_W)_*\QQ_{\widetilde{W}^\times}. \]
Since $R\hat{\Phi}_*$ is an exact functor of triangulated categories and $\hat{\Phi}$ is a homeomorphism, we have a commutative diagram of distinguished triangles
\[\xymatrix{
R\hat{\Phi}_*\QQ_{f^{-1}(0)}\ar@{=}[d]\ar[r] & R\hat{\Phi}_*[\imath_V^*R(\bar\pi_V)_*\QQ_{\widetilde{V}^\times}] \ar@{=}[d]\ar[r] &
R\hat{\Phi}_*P_f[1-\dim V]\ar[r]\ar@{.>}[d] &\\
\QQ_{g^{-1}(0)}\ar[r] & \imath_W^*R(\bar\pi_W)_*\QQ_{\widetilde{W}^\times}\ar[r] & P_g[1-\dim W]\ar[r]&
}\]
which gives us the isomorphism \eqref{2.8} because $\dim V=\dim W$ by $\Phi$.
\end{proof}

We can give an alternative proof by using \eqref{2.2}. Indeed, by \eqref{2.2}, $P_f$ fits into the distinguished triangle
\[ P_f[-n]\lra \QQ_V|_{f^{-1}(0)} \lra Rv_*\QQ_{V_{>0}}|_{f^{-1}(0)}\lra \]
where $v:V_{>0}=\{x\in V\,|\, \Re f(x)>0\} \to V$ is the inclusion and $n=\dim V$. Applying $R\hat{\Phi}_*$, we get a distinguished triangle
\[ R\Phi_*P_f[-n] \lra \QQ_W|_{g^{-1}(0)}\lra Rw_*\QQ_{W_{>0}}|_{g^{-1}(0)}\lra \]
where $w:W_{>0}\to W$ is the inclusion, because $\Phi$ sends $V_{>0}$ to $W_{>0}$. By \eqref{2.2}, $P_g$ fits into the distinguished triangle 
\[ P_g[-n]\lra \QQ_W|_{g^{-1}(0)}\lra Rw_*\QQ_{W_{>0}}|_{g^{-1}(0)}\lra .\]
Therefore, $P_g\cong R\hat{\Phi}_*P_f$. 
\begin{coro}\label{2.11}
Let $\fX$ be a critical virtual manifold equipped with charts $(\fX_\alpha\mapright{\varphi_\alpha} V\lalp \mapright{f\lalp}\CC)$ and let $P\lalp=\varphi^*\lalp P_{f\lalp}$. Then we have an isomorphism
\beq\label{2.11a} \sigma_\albe:P\lalp|_{\fX_\albe}\lra P_\beta|_{\fX_\albe}\eeq
in $D^b_c(\QQ_{\fX_\albe})$.
\end{coro}
\begin{proof}
To simplify the notation, we drop the restriction to $X_\albe$ below. We use the notation of Definition \ref{1.5}. We have a biholomorphic map $\varphi_\albe:V_\albe\to V_\beal$ which gives an  
isomorphism $P_{f\lalp}\mapright{\cong} \varphi^*_\albe P_{f_\beta}$
by \eqref{2.8a}. Applying $\varphi\lalp^*$, we get an isomorphism 
$$P\lalp=\varphi_\alpha^*P_{f\lalp}\mapright{\cong} \varphi_\alpha^*\varphi_\albe^*P_{f_\beta}=\varphi_\beta^*P_{f_\beta}=P_\beta$$
because $\varphi_\beta=\varphi_\albe\circ\varphi_\alpha$. 
\end{proof}
The perverse sheaves of vanishing cycles satisfy the following self-duality.
\begin{prop}\label{2.12} \cite[Proposition 4.2.10]{Dim}
Let $(V,f)$ be an LG pair. Then there exists an isomorphism
$$\mathbf{D}P_f\cong P_f$$
where $\mathbf{D}P_f=R\cH om_{\fX_f}(P_f,\omega_{\fX_f})$ is the Verdier dual of $P_f$ with respect to $\omega_{\fX_f}=p^!\QQ_{\mathrm{pt}}$ via the constant map $p:\fX_f\to \mathrm{pt}$. 
\end{prop}

\bigskip

\section{Gluing of perverse sheaves}\label{S2.2}

In this section, we recall the category of perverse sheaves $Perv(\QQ_X)$ and their gluing properties. We then state a theorem (Theorem \ref{2.20}) which implies that the 2-cocycle obstruction for gluing the local perverse sheaves $P\lalp=\varphi\lalp^*P_{f\lalp}\in Perv(\QQ_{X\lalp})$ from \S\ref{S2.1} coincides with the 2-cocycle obstruction for the orientability of the critical virtual manifold $X$. We thus obtain a global perverse sheaf $P\in Perv(\QQ_X)$ for an \emph{orientable} critical virtual manifold $X$ with charts $X=\cup X\lalp$ which is the gluing of the local perverse sheaves $\{P\lalp\}$ (cf. Theorem \ref{2.23}). 

\medskip

The full subcategory of perverse sheaves in $D^b_c(\QQ_X)$ on an analytic space $X$ is an abelian subcategory of the derived category $D^b_c(\QQ_X)$ of bounded constructible complexes of sheaves of $\QQ$-vector spaces on $X$, whose objects are defined as follows. 
\begin{defi}\label{2.13} An object $P\in D^b_c(\QQ_X)$ is called a \emph{perverse sheaf} (with respect to the middle perversity) if 
\begin{enumerate}
\item $\dim \{x\in X\,|\, H^i(\imath_x^*P)=\bbH^i(B_\varepsilon(x);P)\ne 0\} \le -i$ for all $i$;
\item  $\dim \{x\in X\,|\, H^i(\imath_x^!P)=\bbH^i(B_\varepsilon(x),B_\varepsilon(x)-\{x\};P)\ne 0\} \le i$ for all $i$
\end{enumerate}  
where $\imath_x:\{x\}\hookrightarrow X$ is the inclusion and $B_\varepsilon(x)$ is the open ball of radius $\varepsilon$ centered at $x$ for $\varepsilon$ small enough. 
\end{defi}
Perverse sheaves form an abelian category $Perv(\QQ_X)$ which is the core of a t-structure (\cite[\S2]{BBD}).
An example of perverse sheaf is the perverse sheaf $P_f$ of vanishing cycles defined in Definition \ref{2.1} (cf. \cite[Theorem 5.2.21]{Dim}).

Although sheaf complexes do not have the gluing property in general,
it is known that perverse sheaves and their isomorphisms glue.

\begin{prop}\label{2.14}  \cite[Paragraph 2]{BBD}
Let $X$ be a reduced complex analytic space and let $\{X\lalp\}$ be an open covering of $X$.

(1) Suppose that for each $\alpha$ we have $P_\alpha\in Perv(\QQ_{X_\alpha})$ and for each pair $\alpha, \beta$ we have isomorphisms \[ \sigma_{\alpha\beta}:P_\alpha|_{X_\albe}\mapright{\cong}P_\beta|_{X_\albe}\]
satisfying $\sigma_\albe=\sigma_{\beta\alpha}^{-1}$, $\sigma_{\alpha\alpha}=\id_{V_\alpha}$ and the cocycle condition $\sigma_\gaal\circ\sigma_{\beta\gamma}\circ\sigma_{\alpha\beta}=1$. Then $\{P_\alpha\}$ glue to define a perverse sheaf $P$ on $X$ equipped with isomorphisms 
$$\sigma\lalp:P|_{X_\alpha}\mapright{\cong} P_\alpha$$
such that $\sigma_{\alpha\beta}$ equals the composition 
$$\sigma_\beta\circ\sigma\lalp^{-1}:P\lalp|_{X_\albe} \mapleft{\cong} P|_{X_\albe}\mapright{\cong} P_\beta|_{X_\albe}.$$

(2) Suppose $P, Q\in Perv(\QQ_X)$ and $\sigma_\alpha:P|_{X_\alpha}\mapright{\cong}Q|_{X_\alpha}$ are isomorphisms such that $\sigma_\alpha|_{X_\albe}=\sigma_\beta|_{X_\albe}$. Then there exists an isomorphism $\sigma:P\to Q$ such that $\sigma|_{X_\alpha}=\sigma_\alpha$ for all $\alpha$.
\end{prop}
By the Riemann-Hilbert correspondence established by Kashiwara and Mebkhout, perverse sheaves correspond to regular holonomic D-modules. As sheaves of D-modules glue, so do perverse sheaves. 

\medskip

Our gluing isomorphism $\sigma_\albe$ in Corollary \ref{2.11} is not an arbitrary isomorphism but arose from a biholomorphic map. Recall that an LG pair $(V,f)$ is a holomorphic function $f$ on a complex manifold $V$ that has only one critical value $0$. For an LG pair $(V,f)$, $X_f$ denotes the critical locus $\Crit(f)$ defined by the partial derivatives of $f$.
\begin{defi}\label{2.15} (Geometric origin)
Let $(V_1,f_1)$ and $(V_2,f_2)$ be two LG pairs, and let $\zeta: X_{f_1}\to X_{f_2}$ be an isomorphism of analytic spaces.
An isomorphism $\sigma:P_{f_1}\mapright{\cong} \zeta\sta P_{f_2}$ of perverse sheaves 
is said to be \emph{of geometric origin} if there exists an open neighborhood $X_{f_1}\sub U_1\sub V_1$, a holomorphic map
$\varphi: U_1\to V_2$ biholomorphic onto its image such that 
\beq\label{2.16}\varphi|_{X_{f_1}}=\zeta,\quad f_2\circ\varphi=f_1|_{U_1}\eeq and $\sigma$
is the isomorphism defined in Proposition \ref{2.7} induced from $\varphi$.
\end{defi}

\begin{defi}\label{2.17}
Let $\fX$ be a critical virtual manifold equipped with charts $(\fX_\alpha\mapright{\varphi\lalp}V\lalp\mapright{f\lalp}\CC)$. 
A \emph{geometric gluing} of the perverse sheaves $P_\alpha=\varphi\lalp^{*}P_{f\lalp}$ of vanishing cycles on $\fX_\alpha$ is a perverse sheaf $P$ on $X$ together with isomorphisms $\sigma_\alpha:P|_{X_\alpha}\to P_\alpha$  such that 
$$\sigma_\beta\circ \sigma_\alpha^{-1}:P_{\alpha}|_{X_\albe}\lra P_\beta|_{X_\albe}$$
are of geometric origin.
\end{defi}

The following is immediate from Definition \ref{2.17} and Proposition \ref{2.14}.

\begin{coro}\label{2.18}
Let $\fX$ be a critical virtual manifold equipped with charts $(\fX_\alpha\mapright{\varphi\lalp}V\lalp\mapright{f\lalp}\CC)$. 
There exists a geometric gluing of $\{P\lalp\}$ if possibly after a refinement of the covering $X=\cup\lalp X\lalp$ of $X$, there exists a 1-cochain $\mu=\{\mu_\albe\}$ taking values in $\ZZ_2=\{\pm 1\}$ such that for $\bar\sigma_\albe=\mu_\albe\sigma_\albe$, and for $X_\albega\ne \emptyset$, we have
$$\bar\sigma_\albega:=\bar\sigma_\gaal\circ\bar\sigma_\bega\circ\bar\sigma_\albe=1
$$ where $X_\albega=X_\alpha\cap X_\beta\cap X_\gamma$.
\end{coro}

The gluing condition in this corollary demands, in particular, that the 2-cocycle
\beq\label{2.19}
\sigma_\albega:=\sigma_\gaal\circ\sigma_\bega\circ\sigma_\albe:P\lalp|_{X_\albega}\lra P\lalp|_{X_\albega}
\eeq
be locally constant with values $\mu_\albega=\mu_\gaal\mu_\bega\mu_\albe\in \ZZ_2=\{\pm 1\}$. 
This is indeed true by the following.

\begin{theo}\label{2.20}
Let $(V,f)$ be an LG pair and $X_f=\Crit(f)$ (cf. Definition \ref{1.4}). Let $U$ be an open subset such that $X_f\sub U\sub V$ and let $\Phi:U\to V$ be a holomorphic map, biholomorphic onto its image such that $f\circ\Phi=f|_U$ and $\Phi|_{\fX_f}=\mathrm{id}_{\fX_f}$.
Then the isomorphism $\sigma$ from Proposition \ref{2.7} is equal to
$$\det (d\Phi|_{X_f})\cdot \id: P_f\lra P_f, 
$$ 
where $\det (d\Phi|_{X_f})$ is locally constant with values in $\{\pm 1\}$. %The same holds for mixed Hodge modules.
\end{theo}

Our proof of this theorem is rather lengthy and independent of the rest of this chapter. 
So we postpone the proof of Theorem \ref{2.20} to Chapter \ref{ch3}. See \cite[Corollary 3.2]{BBDJS} for a different proof.

\vsp

The following are immediate consequences of Theorem \ref{2.20}.

\begin{coro}\label{2.21}
Let $\fX$ be a critical virtual manifold equipped with charts $(\fX_\alpha\mapright{\varphi\lalp}V\lalp\mapright{f\lalp}\CC)$.
Let $\{\xi_\albega\}$ be the $\{\pm1\}$-valued 2-cocycle in \eqref{1.13} for the gluing of the anticanonical line bundles $\{K_\alpha^\vee=\varphi^*_\alpha\det T_{V\lalp}|_{X\lalp\ured}\}$. Let $\{\sigma_\albega\}$ be the 2-cocycle in \eqref{2.19}. Then $\sigma_\albega=\xi_\albega$ whenever $X_\albega\ne \emptyset$.\end{coro}
\begin{proof}
Simply let $\Phi$ be the composition $\varphi_\albega$ in \eqref{mapcoc} and use the definitions of $\xi_\albega$ and $\sigma_\albega$, together with Theorem \ref{2.20}.
\end{proof}

\begin{coro}\label{2.22}
Two geometric gluings $P$ and $\tilde{P}$ of perverse sheaves $\{P_\alpha\}$ can differ only by a $\ZZ_2$-local system, i.e. there exists a $\ZZ_2$-local system $\rho\in H^1(X,\ZZ_2)$ such that 
$$\tilde{P} \cong P\otimes \rho.$$
\end{coro}
\def\tisi{\tilde{\sigma} }
\begin{proof}
We have isomorphisms $\tilde{\sigma}_\alpha:\tilde{P}|_{\fX_\alpha}\to P_\alpha$ and $\tilde{\sigma}_\albe=\tisi_\beta\circ\tisi^{-1}_\alpha$. Then $\sigma_\albe$ and $\tisi_\albe$ are two isomorphisms from $P_\alpha|_{\fX_\albe}$ to $P_\beta|_{\fX_\albe}$ of geometric origin, arising from biholomorphic $\varphi_\albe$ and $\tilde{\varphi}_\albe$ up to sign. The composition 
$$\tisi_\albe^{-1}\circ\sigma_\albe:P_\alpha|_{\fX_\albe}\to P_\alpha|_{\fX_\albe}$$
is also an isomorphism of geometric origin, i.e. it is the pullback 
(cf. Proposition \ref{2.7}) by the biholomorphic map $\tilde{\varphi}_\beal\circ \varphi_\albe$ from an open neighborhood of $\fX_\albe$ in $V_\alpha$ to itself preserving $f_\alpha$ and $\fX_\albe$ up to sign. By Theorem \ref{2.20}, $\rho_\albe=\tisi_\albe^{-1}\circ\sigma_\albe=\pm 1$ which coincides with the determinant of the tangent map $d(\tilde{\varphi}_\beal\circ \varphi_\albe)|_{X_\albe}$. Since both $\{\sigma_\albe\}$ and $\{\tisi_\albe\}$ are cocycles, $\{\rho_\albe\in \ZZ_2\}$ is also a cocycle and defines a $\ZZ_2$-local system $\rho\in H^1(X,\ZZ_2)$.
\end{proof}

\medskip

If $\fX$ is an orientable critical virtual manifold, then possibly after a refinement of the covering $\{X_\alpha\}$ of $X$, we can find a 1-cochain $\{\mu_\albe\}$ with values in $\{\pm 1\}$ such that $\bar\sigma_\albega=\bar\xi_\albega=1$ 
by using the notation of Definition \ref{1.17} and Corollary \ref{2.18}.
%where $\bar\sigma_\albega=\mu_\albega\sigma_\albega$ and $\bar\xi_\albega=\mu_\albega\xi_\albega$. 
Therefore we obtain a geometric gluing of $\{P\lalp\}$. We summarize the above discussion as follows.
\begin{theo}\label{2.23}
Let $\fX$ be an \emph{orientable} critical virtual manifold with charts $(\fX_\alpha\mapright{\varphi\lalp}V\lalp\mapright{f\lalp}\CC)$. 
Then there exists a geometric gluing $P$ of the local perverse sheaves $P\lalp=\varphi\lalp^*P_{f\lalp}$ which is unique up to twisting by a $\ZZ_2$-local system.
\end{theo}
We thus obtain a solution to the categorification problem (Problem \ref{1.55}). 
\begin{coro}\label{2.24}
Let $\fX$ be an orientable critical virtual manifold and $P$ be the perverse sheaf in Theorem \ref{2.23}. Then the Euler characteristic $\chi_c(\fX,P)$ of the hypercohomology $\bbH_c^*(\fX,P)$ of $P$ is equal to the Euler characteristic $\chi_\nu(\fX)$ of $\fX$ weighted by the Behrend function. When $X$ is compact, this in turn equals the Donaldson-Thomas type invariant $DT(\fX)=\deg [\fX]\virt$. 
\end{coro}

\begin{rema}\label{2.25} The geometric gluing condition is highly nontrivial. For instance, when $f_\alpha=0$ for all $\alpha$ so that $P\lalp=\QQ|_{X\lalp}$, the gluing isomorphisms $\QQ|_{X_\albe}\to \QQ|_{X_\albe}$ can only be either $1$ or $-1$. 
In particular, a geometric gluing $P$ can only be the trivial bundle $\QQ_X$ twisted by a $\ZZ_2$ local system when $\fX$ is smooth.\end{rema}

In Part 2, we will prove that a moduli space $\fX$ of simple sheaves on a Calabi-Yau 3-fold $Y$ admits a structure of critical virtual manifold which is orientable when there is a tautological family. Hence the Donaldson-Thomas invariant of $Y$ along $\fX$ can be categorified by the hypercohomology $\bbH_c^*(\fX,P)$ of the perverse sheaves in Theorem \ref{2.23}. As an application, this cohomological Donaldson-Thomas invariant will provide us with a mathematical theory 
of the Gopakumar-Vafa invariant \cite{GoVa}.

\bigskip

\def\gr{\mathrm{gr} }
\def\bbD{\mathbb{D} }
\def\Sym{\mathrm{Sym} }

\section{Gluing of mixed Hodge modules}\label{S2.3}

We proved that there is a perverse sheaf $P$ on an orientable critical virtual manifold $X$ whose hypercohomology $\bbH^*_c(X,P)$ has Euler characteristic equal to the Euler characteristic $\chi_\nu(X)$ of $X$ weighted by the Behrend function, which in turn coincides with the Donaldson-Thomas type invariant $DT(X)$ when $X$ is compact. In this section, 
we show that there is a Hodge theory on $\bbH^*_c(X,P)$ so that the decomposition theorem and the hard Lefschetz theorem hold. A Hodge theory for a perverse sheaf means a mixed Hodge module defined and studied by Morihiko Saito \cite{Sai88, Sai90}. (See \cite{Schn} for a quick survey.) We prove that there is a mixed Hodge module $\cM$ on $X$ whose underlying perverse sheaf is the perverse sheaf $P$ constructed in the previous section (cf. Theorem \ref{2.32}). 

\medskip

\subsection{Hodge modules}\label{2.3.1}
Let $X$ be an analytic space embedded into a complex manifold $\PP$. The category of mixed Hodge modules is independent of $\PP$ thanks to Kashiwara's equivalence \cite[Theorem 1.6.1]{HTT}.  Considering vector fields on $\PP$ as $\CC_\PP$-derivations of holomorphic functions in $\sO_\PP$, the sheaf $D_\PP$ of differential operators is defined as the subsheaf of $End_{\CC_\PP}(\sO_\PP)$ generated by the sheaf $\sO_\PP$ and the tangent bundle $T_\PP$. 
There is a natural filtration of $D_\PP$ defined inductively by $F_0D_\PP=\sO_\PP$ and 
$$F_lD_\PP=\{ P\in End_{\CC_\PP}(\sO_\PP)\,|\,[P,f]\in F_{l-1}D_\PP, \ \forall f\in \sO_\PP\}.$$
It is straightforward to see that $$\gr D_\PP=\oplus_l F_lD_\PP/F_{l-1}D_\PP\cong \Sym T_\PP.$$

A \emph{Hodge module} $\cM$ on $X\subset \PP$ consists of 
\begin{enumerate}
\item a regular holonomic $D_\PP$-module $M$ whose support lies in $X$;
\item a perverse sheaf $P$ on $X$;
\item an isomorphism $DR(M)\cong \CC\otimes_\QQ P$;
\item a good filtration of $M$ by $\sO_\PP$-coherent subsheaves $\{M_i\}$ such that $M_i\cdot D_j\subset M_{i+j}$ and $\gr M=\oplus M_i/M_{i-1}$ is coherent over $\gr D_\PP\cong \Sym T_\PP$
\end{enumerate}
which satisfy suitable local conditions.  See \cite{Sai88, HTT, Schn} for precise statements of these local conditions. 
For our purpose of gluing (mixed) Hodge modules, these local conditions are always satisfied from the start. We refer to \cite{HTT} for the definitions of regular holonomic $D_\PP$-modules, de Rham functor $DR$ etc. The isomorphism (3) is said to give a rational structure on $M$. The category $HM(X)$ of Hodge modules is a full subcategory of the category $HW(X)$ of tuples $(M,P,M_\bullet)$ as above without local conditions. 
%and a Hodge module comes with the notion of weight which is defined inductively via the nearby cycles functors. 

When $V$ is a complex manifold of dimension $n$ and $\omega_V=\wedge^nT^*_V$, the constant variation of Hodge structures gives a Hodge module which consists of $P=\QQ_V[n]$, $M=\omega_V$, and $M_i=M$ for $i\ge -n$ and $0$ for $i<-n$. We denote this Hodge module by $\QQ^H_V[n]$ and call it the constant Hodge module. 

A polarization of a Hodge module $\cM$ of weight $w$ refers to an isomorphism $\cM(w)\cong \mathbf{D}\cM$ in $HM(X)$ where $\mathbf{D}$ denotes the Verdier dual. A Hodge module $\cM\in HM(X)$ is called \emph{polarizable} if it admits a polarization. We let $HM(X)^p\subset HM(X)$ denote the full subcategory of polarizable Hodge modules. For example, the constant Hodge module $\QQ^H_V[n]$ is a polarizable Hodge module because $R\cH om(\QQ[n],\QQ[2n])\cong \QQ[n]$ and 
this extends to the canonical isomorphism 
$\QQ^H_V[n]\cong \mathbf{D}\QQ^H_V[n]$  of Hodge modules. 
Polarizable Hodge modules satisfy the following useful properties. 
\begin{theo}\label{2.30}\cite[Theorem 5.3.1]{Sai88}
Let $f:X\to Y$ be a projective morphism of analytic spaces (admitting an embedding into complex manifolds).
Let $\ell$ be the first Chern class of a relatively ample line bundle on $X$. 
Then for a polarizable Hodge module $\cM$ on $X$,
\begin{enumerate}
\item $^p\!R^if_*\cM\in HM(Y)^p$ for all $i$;
\item $\ell^i:\, ^p\!R^{-i}f_*\cM\lra\ ^p\!R^if_*\cM$ is an isomorphism of Hodge modules.
\end{enumerate}
\end{theo}
The hard Lefschetz property (2) above gives us an isomorphism $$\ell^i: \,^p\!R ^{-i}f_*P\to \,^p\!R ^if_*P$$ where $P$ is the underlying perverse sheaf of $\cM=(M,P,M_\bullet)$. 
Together with Deligne's argument from \cite{Deli} (cf. \cite[p466]{GrHa}) for degeneration of spectral sequences, this isomorphism gives us the decomposition theorem.
\begin{coro}\label{2.31} Under the hypothesis of Theorem \ref{2.30}, we have a (non-canonical) isomorphism
$$Rf_*P\cong \oplus_i\,^p\!R ^if_*P[-i]$$
and each $^p\!R^if_*P[-i]$ is a perverse sheaf underlying a polarizable Hodge module. 
\end{coro}

\medskip
\def\rat{\mathrm{rat} }
\subsection{Mixed Hodge modules}\label{S2.3.2}

As in \S\ref{2.3.1}, $X$ is an analytic space embeddable into a complex manifold $\PP$. 

A \emph{polarizable mixed Hodge module} on $X$ consists of  \begin{enumerate}
\item a regular holonomic $D_\PP$-module $M$ whose support lies in $X$;
\item a perverse sheaf $P$ on $X$;
\item an isomorphism $DR(M)\cong \CC\otimes_\QQ P$;
\item a good filtration of $M$ by $\sO_\PP$-coherent subsheaves $M_\bullet=\{M_i\}$ such that $M_i\cdot D_j\subset M_{i+j}$ and $\gr M=\oplus M_i/M_{i-1}$ is coherent over $\gr D_\PP\cong \Sym T_\PP$;
\item  a finite increasing filtration $W_\bullet$ of $\cM=(M,P,M_\bullet)$ with $\gr^W_i \cM\in HM(X)^p$ for all $i$
\end{enumerate}
which satisfy suitable local conditions (cf. \cite{Sai90}). As our purpose is gluing mixed Hodge modules, these local conditions are always satisfied from the start (before gluing).  The category $MHM(X)^p$ of polarizable mixed Hodge modules is the full subcategory of the category $MHW(X)$ of tuples $(M,P,M_\bullet, W_\bullet)$ without local conditions.

By the definition, we have a forgetful functor 
\beq\label{2.36}\rat: MHM(X)^p\lra Perv(\QQ_X), \quad (M,P,M_\bullet, W_\bullet)\mapsto P\eeq
which is an exact and faithful functor via the Riemann-Hilbert correspondence. 

The constant Hodge module $\QQ^H_V[n]$ on a complex manifold $V$ of dimension $n$ is a polarizable mixed Hodge module (cf. \cite[Theorem 3.8, (4.5.5)]{Sai90}). If $f:V\to \CC$ is a holomorphic function on a complex manifold $V$, then there is a polarizable mixed Hodge module $\cM_f:=\phi_f\QQ^H_V[n]$ supported on $f^{-1}(0)$ such that $$\rat(\cM_f)=P_f$$
is the perverse sheaf of vanishing cycles defined in 
Definition \ref{2.1} where $\phi_f$ denotes the vanishing cycle functor. This is actually part of the local conditions required for mixed Hodge modules. 

By \cite[Theorem 0.1]{Sai90}, when there is a morphism $\Phi$ of analytic spaces, we have natural functors $\Phi_*$, $\Phi_!$, $\Phi^*$, $\Phi^!$, $\psi_g$, $\phi_{g,1}$, $\mathbf{D}$, $\boxtimes$, $\otimes$, and $\cH om$ between the derived categories of mixed Hodge modules. 
By \cite[Theorem 2.14]{Sai90}, when $\Phi:V\to W$ 
is a biholomorphic map of complex manifolds and $g:W\to \CC$ is a holomorphic function with $f=g\circ\Phi$, 
we have a \emph{canonical} isomorphism 
\beq\label{2.39}\sigma^H:\cM_f=\phi_f\QQ^H_V[n]\mapright{\cong} \hat{\Phi}^*\phi_g\QQ^H_W[n] =\hat{\Phi}^*\cM_g\eeq
which induces \eqref{2.8a} when $\rat$ is applied. Here $\hat{\Phi}:f^{-1}(0)\to g^{-1}(0)$ is the isomorphism of analytic spaces induced from $\Phi$. 
%This follows from the definition \cite[(2.2.6)]{Sai90} of the vanishing cycle functors $\phi_f$ and $\phi_g$. 
By \cite[Proposition 2.6]{Sai90}, if $f:V\to\CC$ is a holomorphic function, we have a \emph{canonical} isomorphism
\beq\label{2.37}\phi_f\mathbf{D}\QQ^H_V[n]\cong \mathbf{D}\phi_f\QQ^H_V[n]=\mathbf{D}\cM_f\eeq
where $n=\dim V$. 
The canonical isomorphism $\QQ^H_V[n]\cong \mathbf{D}\QQ^H_V[n]$ sending 1 to 1 gives us a \emph{canonical} isomorphism 
\beq\label{2.38} \cM_f=\phi_f\QQ^H_V[n]\cong \phi_f\mathbf{D}\QQ^H_V[n].\eeq
Composing \eqref{2.38} with \eqref{2.37}, we obtain a \emph{canonical} isomorphism
\beq\label{2.40} \cM_f\lra \mathbf{D}\cM_f,\eeq
i.e. $\cM_f$ has a \emph{canonical} polarization. Since \eqref{2.40} and \eqref{2.39} are canonical, the isomorphism \eqref{2.39} is an isomorphism of polarizable mixed Hodge modules.

By \cite[\S2]{Sai90} and \cite[\S1.6]{Sai91}, if $X=\cup_\alpha X_\alpha$ is an open cover, the category $MHM(X)^p$ of polarizable mixed Hodge modules is equivalent to the category of the collections $\cM_\alpha\in MHM(X_\alpha)^p$ together with isomorphisms 
\beq\label{2.41}\sigma_\albe^H:\cM_\alpha|_{X_\albe}\to \cM_\beta|_{X_\albe}\eeq 
in $MHM(X_\albe)^p$, satisfying 
$$\sigma_\albega^H=\sigma_\gaal^H\circ \sigma_\bega^H\circ \sigma_\albe^H=1$$ whenever $X_\albega\ne \emptyset$.  Likewise, isomorphisms can be glued.
In other words, Proposition \ref{2.14} holds for polarizable mixed Hodge modules.

Combining all in this subsection, we obtain the following theorem. We use the notation of Definition \ref{1.5}.
\begin{theo}\label{2.32}
Let $X$ be an orientable critical virtual manifold with charts $(X\lalp\mapright{\varphi\lalp}V\lalp\mapright{f\lalp}\CC)$.
Let $\cM\lalp=\varphi\lalp^*\cM_{f\lalp}\in MHM(X\lalp)^p$ such that $\rat(\cM\lalp)=P\lalp$. Let $\sigma_\albe^H:\cM\lalp|_{X_\albe}\to \cM_\beta|_{X_\albe}$ be the isomorphism induced from \eqref{2.39} by the biholomorphic map $\varphi_\albe$ as in Corollary \ref{2.11}.
There is a polarizable mixed Hodge module $\cM\in MHM(X)^p$ together with isomorphisms $\sigma\lalp^H:\cM|_{X\lalp}\mapright{\cong} \cM\lalp$ such that $\sigma_\albe^H$ equals the composition $\sigma_\beta^H\circ{\sigma\lalp^H}^{-1}:\cM\lalp|_{X_\albe}\mapleft{\cong}\cM|_{X_\albe}\mapright{\cong} \cM_\beta|_{X_\albe}$. Moreover, $\rat(M)$ is the perverse sheaf $P$ in Theorem \ref{2.23}.
\end{theo}
\begin{proof} Since \eqref{2.39} is an isomorphism of polarizable mixed Hodge modules, it gives an isomorphism $\sigma_\albe$ of polarizable mixed Hodge modules for each pair $(\alpha,\beta)$, as in the proof of Corollary \ref{2.11}.
Since the functor $\rat$ is faithful and $\rat(\sigma_\albega^H)=\sigma_\albega$ where $\sigma_\albega$ is from \eqref{2.19}, the gluing condition $\sigma_\albega=1$ for the perverse sheaves $\{P\lalp\}$ implies the gluing condition $\sigma_\albega^H=1$ for the mixed Hodge modules $\{\cM\lalp\}$. As the perverse sheaves $\{P\lalp\}$ glue, so do the polarizable mixed Hodge modules $\{\cM\lalp\}$.
\end{proof}

Let $P$ be the perverse sheaf which underlies a polarizable mixed Hodge module. 
Let $$\hat{P}=\gr^W P$$ denote the gradation with respect to the weight filtration. Then $\hat{P}$ is a direct sum of polarizable Hodge modules. If $\psi:X\to Y$ is a projective morphism of analytic spaces, then the hard Lefschetz theorem (Theorem \ref{2.30}) and the decomposition theorem (Corollary \ref{2.31}) hold for $\hat{P}$. On the other hand, by the long exact sequence for hypercohomology groups and Corollary \ref{1.59}, we have
$$\chi_c(X,\hat{P})=\chi_c(X,{P})=\chi_\nu(X)=DT(X).$$
 We therefore have another solution to Problem \ref{1.55}. 
 \begin{coro}\label{2.33}
 Let $X$ be an orientable critical virtual manifold. Then there is a perverse sheaf $\hat{P}$ underlying a direct sum of polarizable Hodge modules such that the Euler characteristic $\chi_c(X,\hat{P})$ of the hypercohomology $\bbH^*_c(X,\hat{P})$ equals the Euler characteristic $\chi_\nu(X)$ of $X$ weighted by the Behrend function. 
 \end{coro}
We will use $\hat{P}$ for a mathematical theory of the Gopakumar-Vafa invariant.

%%%%%%%%%%%%%%%%%
%%%%%%%%%%%%%%%%%
%%%%%%%%%%%%%%%%%
%%%%%%%%%%%%%%%%%
%%%%%%%%%%%%%%%%%
\chapter{Rigidity of perverse sheaves}\label{ch3}
 
In this chapter, we provide a proof of Theorem \ref{2.20}.
The main points are the following: \begin{enumerate}
\item Perverse sheaves and their isomorphisms are rigid under continuous deformations.
\item We can use vector fields to produce isotopies of 
equivalences of LG pairs which generate continuous deformations of isomorphisms of perverse sheaves.
\item We can find isotopies from a self-equivalence $\Phi$ of an LG pair to linear transformations. 
\end{enumerate}
Our proof reveals the dependence of the perverse sheaf of vanishing cycles 
on the obstruction theory of the critical locus.

\bigskip

\section{Sebastiani-Thom isomorphism and rigidity}
\label{S3.1}
We begin with a few properties of perverse sheaves that we will use. 

\begin{prop}\label{3.1}  \cite[\S2]{Massey} (Sebastiani-Thom isomorphism)\\
Let $g:V\to \CC$ and $h:W\to \CC$ be holomorphic functions on complex manifolds. Let $f=g+h:V\times W\to \CC$ be defined by $f(z,y)=g(z)+h(y)$. Then we have an isomorphism of perverse sheaves 
$$P_f\cong pr_1^{-1}P_g\otimes pr_2^{-1} P_h $$
where $pr_1:V\times W\to V$ and $pr_2:V\times W\to W$ denote the projections.
\end{prop}
\begin{coro}\label{3.2}
When $h$ is the quadratic function $\sum_{i=1}^ry_i^2$ on $\CC^r$ in Example \ref{2.9}, $P_h=\QQ_0$ 
and hence $P_f\cong P_g$ as perverse sheaves defined on the critical locus $\fX_g=\fX_f$.
\end{coro}

We recall the following ``elementary construction'' of perverse sheaves by MacPherson and Vilonen.
\begin{theo}\label{3.3} \cite[Theorem 4.5]{MV}
Let $X$ be an analytic space. 
Let $S\subset X$ be a closed stratum of complex codimension $c$. The category $Perv(\QQ_X)$ is equivalent to the category of objects $(B,C)\in Perv(\QQ_{X-S})\times Sh_{\QQ}(S)$ together with a commutative triangle
\beq\label{3.4}\xymatrix{ R^{-c-1}\pi_*\kappa_*\kappa^*B\ar[rr]\ar[dr]_m && R^{-c}\pi_*\gamma_!\gamma^*B\\ 
& C\ar[ur]_n  }\eeq
such that $\ker (n)$ and $\coker (m)$ are local systems on $S$, where $\kappa:K\hookrightarrow L$ and $\gamma:L-K\hookrightarrow L$ are inclusions of the perverse link bundle $K$ and its complement $L-K$ in the link bundle $\pi:L\to S$. The equivalence of categories is explicitly given by sending $P\in Perv(\QQ_X)$ to $B=P|_{X-S}$ together with the natural morphisms
\[\xymatrix{
 R^{-c-1}\pi_*\kappa_*\kappa^*B\ar[rr]\ar[dr]_m && R^{-c}\pi_*\gamma_!\gamma^*B\\ 
& R^{-c}\pi_*\varphi_!\varphi^*P\ar[ur]_n }\]
where $\varphi:D-K\hookrightarrow D$ is the inclusion into the normal slice bundle.
\end{theo}
See \cite[\S4]{MV} for precise definitions of links $K$, $L$ and $D$. Morally the above theorem says that an extension of a perverse sheaf on $X-S$ to $X$ is obtained by adding a sheaf on $S$. 

An application of Theorem \ref{3.3} is the following \emph{rigidity property of perverse sheaves}.
\begin{lemm}\label{3.5} Let $P\in Perv(\QQ_U)$ be a perverse sheaf on an analytic space $U$. Let $\pi:T\to U$ be a continuous map from a topological space $T$ with connected fibers and let $T'$ be a subspace of $T$ such that $\pi|_{T'}$ is surjective. Suppose that an isomorphism 
$\mu:\pi^{-1}P\mapright{\cong} \pi^{-1}P$ satisfies $\mu|_{T'}=\id_{(\pi^{-1}P)|_{T'}}$. Then $\mu=\id_{\pi^{-1}P}$.
\end{lemm}
\begin{proof}
We first prove the simple case. If we let $C$ be a locally constant sheaf over $\QQ$ of finite rank supported on a subset $Z\subset U$ and $\bar\mu:\pi^{-1}C\to \pi^{-1}C$ be a homomorphism such that $\bar\mu|_{T'\cap \pi^{-1}(Z)}=\id$, then $\bar\mu$ is the identity morphism. Indeed, since the issue is local, we may assume that $Z$ is connected and that $C\cong \QQ^r$ so that $\bar\mu:\QQ^r\to\QQ^r$ is given by a continuous map $\pi^{-1}(Z)\to GL(r,\QQ)$. By fiber connectedness, this obviously is a constant map which is $1$ along $T'\cap\pi^{-1}(Z)$. We thus proved the lemma in the sheaf case.

For the general case, we use Theorem \ref{3.3}. By replacing $U$ by the support of $P$ if necessary, we may assume that the support of $P$ is $U$.  We stratify $U$ and let $U^{(i)}$ be the union of strata of codimension $\le i$. Since $P$ is a perverse sheaf, $P|_{U^{(0)}}[-\dim U]$ is isomorphic to a locally constant sheaf and hence $\mu|_{U^{(0)}}$ is the identity map. For $U^{(1)}=U^{(0)}\cup S$, using the notation of Theorem \ref{3.3}, $C=R^{-1}\pi_*\varphi_!\varphi^*P$ is a locally constant sheaf and $\mu$ induces a homomorphism $\pi^{-1}C\to \pi^{-1}C$ which is identity on $T'\cap \pi^{-1}(S)$. Therefore $\mu$ induces the identity morphism of the pullback of \eqref{3.4} by $\pi$ to itself and hence $\mu$ is the identity on $U^{(1)}$. Continuing in this fashion, we obtain Lemma \ref{3.5}.
\end{proof}

Typically we will consider the case where $\pi:T=U\times I\to U$ is the projection and $I$ is the interval $[0,1]$ or a disc in $\CC$. If there is a continuous family $\mu$ of isomorphisms $P\to P$ parameterized by $I$ which is $\id$ over $U\times \{0\}$, all the isomorphisms are the identity. More precisely, let $\Phi:V\times I\to V$ be a continuous family of homeomorphisms of a complex manifold $V$, i.e. $\Phi$ is continuous and $\Phi_t:=\Phi|_{V\times\{t\}}$ are homeomorphisms for all $t\in I$. 
Let $f:V\to \CC$ be a holomorphic function, satisfying $f(\Phi(x,t))=f(x)$, i.e. $f\circ \Phi_t=f$ for all $t$. 
As always, we let $X_f=\Crit(f)$ denote the critical locus of $f$. 
Suppose $\Phi_t|_{X_f}$ is the identity map for all $t$. Then the homeomorphism $\Phi_t$ gives us the isomorphism
\beq\label{3.8} \sigma_t:P_f\lra P_f\eeq
by Proposition \ref{2.7}.
 
Let $P_f=R\Gamma_{\{x\in V\,|\,\mathrm{Re} f(x)\le 0\}}\QQ_V[\dim V]|_{f^{-1}(0)}$ be the perverse sheaf of vanishing cycles from Definition \ref{2.1}. Since $\pi^{-1}=\pi^*$ is an exact functor (on $\QQ$-sheaves), we have $\pi^{-1}P_f=R\Gamma_{\{x\in V\,|\,\mathrm{Re} f(x)\le 0\}\times I}\QQ_{V\times I}[\dim V]|_{f^{-1}(0)\times I}$. 
By the proof of Proposition \ref{2.7}, we then have an isomorphism 
\beq\label{3.7}\mu:\pi^{-1}P_f\lra \pi^{-1}P_f\eeq 
whose restriction to $X_f\times\{t\}$ is the isomorphism $\sigma_t:P_f\to P_f$ in \eqref{3.8}. 
In this situation, Lemma \ref{3.5} gives the following.
\begin{prop}\label{3.6} Let $(V,f)$ be an LG pair (Definition \ref{1.4}). Let $\Phi_t$, $P_f$, $\sigma_t$ be as above. 
Suppose $\Phi_0$ is the identity map of $V$ and $\Phi_t|_{X_f}:X_f\to X_f$ is the identity map for all $t\in I$.
Then $\sigma_t=\id_{P_f}$ for all $t\in I$.
\end{prop}
\begin{proof}
The proposition is immediate from Lemma \ref{3.5} by letting $U=X_f$, $T=U\times I\mapright{\pi} U$ and $\mu$ given by \eqref{3.7}, because $\mu|_{U\times\{0\}}=\id$ by $\Phi_0=\id_V$.
\end{proof}

\bigskip 

\section{Vector fields and isotopies} \label{S3.2}
As the issue of Theorem \ref{2.20} is local, we may restrict our concern to an open submanifold $V\sub\CC^n$. In this section, we use vector fields to generate isotopies. Together with the rigidity (Proposition \ref{3.6}), this will give us a 2-cocycle property of perverse sheaves of vanishing cycles (Lemma \ref{3.12}).

\begin{lemm}\label{3.9}
Let $V\sub\CC^n$ be an open submanifold and $f:V\to\CC$ be an LG pair (Definition \ref{1.4}).   
Let $(df)$ be the ideal generated by the partial derivatives of $f$ and let $X_f=\Crit(f)$ denote the analytic subspace defined by the ideal $(df)$. Then there is an open neighborhood $V'$ of $X_f$ in $V$ such that
$$f|_{V'}\in (df)|_{V'} \and \frac{f}{|\!|df|\!|}\to 0\text{ as } df\to 0.$$
\end{lemm}

\begin{proof}
Let $\pi:\tilde V\to V$ be a resolution of $X_f$ so that $\pi^{-1}(\fX_f)$ is a normal crossing divisor. 
Near every $\tilde x\in \pi\upmo(\fX_f)$, the pullback $\pi^*(df)$ of the ideal $(df)$ is a principal ideal sheaf generated by some monomial $\varphi=z_1^{k_1}\cdots z_r^{k_r}$ with $k_i>0$ for a system $z_1,\cdots,z_n$ of local coordinates of $\tilde V$ centered at $\tilde x$. Let $x=\pi(\tilde x)$, and let $w_1,\cdots, w_n$ be the coordinate functions of $\CC^n$. Then $\varphi$ divides $\pi\sta \frac{\partial f}{\partial w_i}$ for all $i$. 

We claim that $\pi^*f$ is divisible by $z_1^{k_1+1}\cdots z_r^{k_r+1}$. We first write $\pi^*f=c_{m_1}z_1^{m_1}+c_{m_1+1}z_1^{m_1+1}+\cdots$, where $c_ k$ are holomorphic functions of $\{z_2,\cdots, z_n\}$. Because near $\tilde x$ and away from
$\pi\upmo(X_f)$, $\pi$ is biholomorphic, 
\beq\label{3.10}
\frac{\partial (\pi\sta f)}{\partial z_1}=\sum_i \pi\sta \left(\frac{\partial f}{\partial w_i}\right)\cdot \frac{\partial (\pi\sta w_i)}{\partial z_1}
\eeq
holds away from $\pi\upmo(X_f)$. As all terms in this identity are holomorphic, it holds in a neighborhood of $\tilde x$. 
Because the right hand side is divisible by $\varphi$, so is the left hand side. Since
$$\frac{\partial (\pi\sta f)}{\partial z_1}=m_1c_{m_1}z_1^{m_1-1}+(m_1+1)c_{m_1+1}z_1^{m_1}+\cdots ,
$$
we must have $m_1-1\ge k_1$. Therefore $z_1^{k_1+1}$ divides $\pi^*f$. Likewise $z_i^{k_i+1}$ divides $\pi^*f$ for each $i$. 
Therefore by \eqref{3.10}, 
\beq\label{3.11}
\pi^*(f)\subset \pi^*(df)\,\sqrt{\pi^*(df)}\subset \pi^*(df).
\eeq
Since $\pi^*:\cO_V\to \pi_*\cO_{\tilde V}$ is injective, we have $f\in (df)$. 

Finally, using \eqref{3.11}, near $\tilde x$ we write $\pi\sta f=\sum b_i\pi\sta\bl \frac{\partial f}{\partial w_i}\br$ for some holomorphic functions $b_i$ such that $b_i$ vanishes along $\pi^{-1}(X_f)$. Since $\pi: \ti V\to V$ is proper, we have $\frac{f}{|\!|df|\!|}\to 0$ as $df\to 0$.
This proves the lemma.
\end{proof}

\begin{lemm}\label{3.12}
Let $V\sub \CC^n$ be open. Let $(V,f_0)$ and $(V,f_1)$ be two LG pairs. 
Let $f_t=(1-t)f_0+tf_1$ for $0\le t\le 1$. Suppose the critical locus $\fX_{f_t}=\Crit(f_t)$ is independent of $t$ as an analytic subspace of $V$. Then there is an isomorphism 
$$\tau_{01}:P_{f_0}\mapright{\cong} P_{f_1}.
$$
If we have a third LG pair $(V,f_2)$ such that the critical locus of $(1-t-s)f_0+tf_1+sf_2$ is independent of $t,s\in [0,1]$. Then $$\tau_{12}\circ\tau_{01}=\tau_{02}:P_{f_0}\lra P_{f_2}.$$
\end{lemm}

\begin{proof}
We use the standard hermitian inner product on $\CC^n$. Let $\nabla f_t$ be the gradient vector field of $f_t$ defined by $df_t(v)=\nabla f_t\cdot v$ for tangent vectors $v$. By Lemma \ref{3.9}, we have 
$\frac{f}{|\!|df|\!|}\to 0$ as $df\to 0$ in a neighborhood of $X_f$.
We define a time dependent vector field 
\[ \xi_t=\frac{f_0-f_1}{|\!|\nabla f_t|\!|^2}\ \overline{\nabla f_t}.
\]
We claim that this is a well defined vector field on $V$. It suffices to show that 
\[ |\!|\xi_t|\!|=\frac{|f_0-f_1|}{|\!|\nabla f_t|\!|}=\frac{|f_0-f_1|}{|\!|d f_t|\!|}\]
approaches zero as the point approaches $\fX_{f_t}=\fX_{f_0}=\fX_{f_1}$. Since $(df_t)\supset (df_0)=(df_1)$ by assumption, we can express the partial derivatives of $f_0$ and $f_1$ as linear combinations of the partial derivatives of $f_t$. Hence,
\[ |\!|d f_0|\!|\le C|\!|d f_t|\!|\and |\!|d f_1|\!|\le C |\!|d f_t|\!|\]
for some $C>0$. Thus by Lemma \ref{3.9},
\[ |\!|\xi_t|\!|=\frac{|f_0-f_1|}{|\!|d f_t|\!|}\le C^{-1}\left( \frac{|f_0|}{|\!|d f_0|\!|}+\frac{|f_1|}{|\!|d f_1|\!|} \right)\to 0 \text{ as } df_0, df_1\to 0. \]
So we proved the claim. 

Let $x_t$ for $t\in [0,1]$ be an integral curve of the vector field $\xi_t$, so that $$\frac{dx_t}{dt}=\dot{x}_t=\xi_t(x_t).$$
Then $f_t(x_t)$ is constant because $\frac{d}{dt}f_t(x_t)$ is equal to
\[ df_t(\dot{x}_t)+f_1-f_0=\nabla f_t\cdot\dot{x}_t+f_1-f_0=\nabla f_t\cdot \frac{f_0-f_1}{|\!|\nabla f_t|\!|^2}\ \overline{\nabla f_t}+f_1-f_0=0.\]
Therefore the flow of the vector field $\xi_t$ from $t=0$ to $t=1$ gives a homeomorphism $\Phi_{01}:U\to U'$ of  neighborhoods of $X_{f_0}$ such that $f_1(\Phi_{01}(x))=f_0(x)$ for $x\in U$.  If $x\in X_{f_0}=X_{f_1}$, $f_0(x)=f_1(x)=0$ and hence $\xi_t(x)=0$ for all $t$. So $\Phi_{01}|_{X_{f_0}}$ is the identity map of $X_{f_0}$. By Proposition \ref{2.7}, we obtain an isomorphism $\tau_{01}:P_{f_0}\cong P_{f_1}$. 

Suppose that we have three holomorphic functions $f_0, f_1,f_2$ on $V$ as stated in Lemma \ref{3.12}.  The composition $\tau_{12}\circ\tau_{01}$ is obtained from the composition of the diffeomorphism $\Phi_{01}:U\to U'$ with $f_1(\Phi_{01}(x))=f_0(x)$ and the diffeormorphism $\Phi_{12}:U'\to U''$ with $f_2(\Phi_{12}(x))=f_1(x)$. By replacing $f_1$ by $(1-s)f_1+sf_2$ for $0\le s\le 1$, 
we obtain an isotopy from $\Phi_{12}\circ\Phi_{01}$ to $\Phi_{02}$, so that $\Phi_{02}^{-1}\circ\Phi_{12}\circ\Phi_{01}$ is isotopic to the identity map. By Proposition \ref{3.6}, we have $\tau_{02}^{-1}\circ\tau_{12}\circ\tau_{01}=\mathrm{id}$.
\end{proof}

The following is a consequence of Lemmas \ref{3.9} and \ref{3.12}.
\begin{lemm}\label{3.14}
Let $V\sub \CC^n$ be open and let $(V,g)$ be an LG pair (Definition \ref{1.4}).
Let $u$ be a nowhere vanishing holomorphic function on $V$and let $f=ug$.
Then there is an open neighborhood
$X_g\sub U\sub V$ so that $\fX_f\cap U=\fX_g$.
Suppose further $u|_{\fX_f\cap U}=1$. Then $g=uf$ induces a canonical isomorphism $P_f|_U\cong P_g$.
\end{lemm}
\begin{proof} By Lemma \ref{3.9}, $g\in (dg)$ on an open neighborhood of $X_g$. Since $f=ug$ and $df=u\,dg+g\,du$, the ideal $(df)$ is contained in $(dg)$, thus $\fX_g\sub \fX_f$. Let $\mathfrak A=\fX_f\cap (f=0)$. By Lemma \ref{3.9}, $\mathfrak A\sub \fX_f$ is both open and closed. We let $U\sub V$ be an open neighborhood of $\mathfrak A\sub V$ so that $U\cap \fX_f=\mathfrak A$.
Then using $g=u^{-1}f$, the same argument shows that $(dg)$ is contained in $(df|_U)$. Hence $\fX_f\cap U= \fX_g$. 

Suppose $u|_{\fX_f}=1$. We 
let $u_t=(1-t)+tu: V\to\CC$. Then $u_t$ are invertible in a neighborhood of $X_f$ for any $t$. Let $f_t=u_tg$;
then $f_0=g$ and $f_1=f$. By the same argument, and using that $[0,1]$ is compact, 
we can find an open neighborhood $\fX_g\sub U'\sub V$ so that $\fX_{f_t}\cap {U'}=\fX_{f_0}$ for all $t\in [0,1]$. Applying 
Lemma \ref{3.12}, we obtain an induced isomorphism $P_f|_{U}=P_{f_1}|_{U'}\cong P_{f_0}|_{U'}=P_g$,
thus proving the lemma.
\end{proof}

\bigskip

\section{Obstruction theory and isotopy}\label{S3.3}

In order to use Lemma \ref{3.12} for a proof of Theorem \ref{2.20}, we have to make sure that the critical loci $X_{f_t}$ is independent of $t$. In this section, we will use the obstruction assignments (Definition \ref{1.34}) to give a criterion for the constancy of the critical loci. 

Recall that an LG pair $(V,f)$ gives us a symmetric obstruction theory 
$$E_V=[T_V|_{X_f}\mapright{d(df)}\Omega_V|_{X_f}] \lra \bbL_{X_f}
$$ 
and there is an obstruction class $ob_{X_f}(\phi,\bar{g},B,\bar{B})\in I\otimes_\CC\Omega_{X_f}|_x$ for an infinitesimal lifting problem (Definition \ref{1.8}).  We proved that two equivalent LG pairs give the same obstruction assignment (Lemma \ref{1.43}). A natural case where non-equivalent LG pairs give the same obstruction assignment is the following. 

\begin{lemm}\label{3.15}
Let $V\sub \CC^n$ be open and let $(V,f)$ be an LG pair with a critical point $x\in X_f=\Crit(f)$. 
Then there exist a system of coordinates $y_1,\cdots,y_n$ centered at $x$ in a 
neighborhood $U$ of $x$ in $V$, a holomorphic function $h$ on $U$ in the form $h=h(y_{r+1},\cdots,y_n)$, and an invertible function $u$ on $U$ with $u|_{\fX_f\cap U}=1$ such that 
\[ f|_U=u\cdot (y_1^2+\cdots +y_r^2+h)\and T_x\fX_f=T_xU'
\]
where $U'=(y_1=\cdots=y_r=0)\cap U$. If we let $h'=h|_{U'}$, then $x\in \fX_{h'}=\fX_f\cap U$. Moreover, we have a canonical isomorphism $P_f|_{U}\cong P_{h'}$. 
\end{lemm}

\begin{proof}
If $T_x\fX_f=T_xV$, then there is nothing to prove. So we assume $T_x\fX_f\ne T_xV$. Choose a local coordinate system $w_1,\cdots, w_{n}$ of $V$ centered at $x$, and let 
$$H=\left(\frac{\partial^2f}{\partial w_i\partial w_j}\Big\vert_x\right)
$$
be the Hessian matrix of $f$ at $x$. Since $\fX_f$ is the vanishing locus of $\frac{\partial f}{\partial w_i}$, we have  $T_x\fX_f=\ker H\ne T_xV$. Since $H$ is symmetric, there exists an element $v\in T_xV$ such that $v^THv\ne 0$. After a coordinate change, we may assume that $v=(1,0,0,\cdots,0)$, i.e. $\frac{\partial^2f}{\partial w_1^2}|_x\ne 0$. By the Weierstrass preparation lemma \cite[p7]{GrHa},  we can write $$f=u_1\cdot (w_1^2+w_1p_1+q_1)$$ for some invertible holomorphic function $u_1$ and 
holomorphic functions $p_1,q_1$ in $w_2,\cdots,w_{n}$. By completing the square, we can write $f=u_1\cdot ({\bar w}_1^2+\tilde{q}_1(w_2,\cdots,w_{n}))$ with  $\bar w_1=w_1+p_1/2$. We can repeat the same argument with $\tilde q_1$ in place of $f$ to obtain $f=u_1\cdot({\bar w}_1^2+u_2({\bar w}_2^2+\tilde{q}_2(w_3,\cdots,w_{n})))$ for some invertible $u_2$. We continue this way until we reach 
$$f=u_1{\bar w}_1^2+u_1u_2{\bar w}_2^2+\cdots +\left(\prod_{i=1}^r u_i\right){\bar w}_r^2+\left(\prod_{i=1}^{r}u_i\right)\tilde h(w_{r+1},\cdots,w_{n})
$$
for invertible functions $u_1,\cdots, u_{r}$. Let $\tilde u=\prod_{i=1}^{r}u_i$, $\tilde y_j=(\prod_{i=j+1}^{r} u_i)^{-\frac12}\bar w_j$ for $j=1,\cdots,r$ and $y_{r+k}=w_{r+k}$ for $k\ge 1$. Then 
$$f=\tilde u\cdot (\tilde h(y_{r+1},\cdots,y_n)+\tilde y_1^2+\cdots +\tilde y_r^2) 
$$ 
in a neighborhood $U$ of $x$, expressed in terms of the coordinate variables $(\ti y,y)=(\ti y_1,\cdots, \ti y_r,y_{r+1},\cdots, y_n)$.
Let $u=\tilde u(\tilde y,y)/\tilde u(0,y)$,  $y_j=\sqrt{\tilde u(0,y)}\tilde y_j$ for $j\le r$ and $h={\tilde u(0,y)}\tilde h$ so that we have
$$f=u\cdot (y_1^2+\cdots + y_r^2+h(y_{r+1},\cdots,y_n)) ,\quad u(0,z)=1.
$$
We let 
$$g(y_1,\cdots,y_n)=y_1^2+\cdots + y_r^2+h(y_{r+1},\cdots,y_n).
$$
By Lemma \ref{3.14}, possibly after shrinking $x\in U$, we have $\fX_f\cap U=\fX_g\cap U$, and canonical isomorphism $P_f|_U\cong P_g|_U$. 
Applying the Sebastiani-Thom isomorphism, we get $P_g|_U\cong P_{h'}|_{U'}$, and thus an isomorphism $P_g|_U\cong P_{h'}|_{U'}$. This proves the lemma. 
\end{proof}

\begin{lemm}\label{3.16}
Let the situation be as in Lemma \ref{3.15}. Then the symmetric obstruction theories defined by the LG pairs $(V,f)$ and $(V',h')$ give the same obstruction assignment of $\fX_f\cap U=\fX_{h'}$ at $x$.
\end{lemm}

\begin{proof}The proof is parallel to that of Lemma \ref{1.43}, using that $u(x)=1$. 
We will omit the details here. 
\end{proof}

\begin{prop}\label{3.17} Let $V\subset \CC^n$ be open. 
Let $(V,f_0)$ and $(V,f_1)$ be two LG pairs such that $\fX_{f_0}=\fX_{f_1}$ with $x\in \fX_{f_0}$.
Let $f_t=(1-t)f_0+tf_1$. 
Suppose that the symmetric obstruction theories defined by $(V,f_0)$ and $(V,f_1)$ give the same obstruction assignment at $x$,
and that $f_1-f_0\in \fm_x^3$ where $\fm_x$ is the ideal of holomorphic functions on $V$ vanishing at $x$. 
Then there is an open neighborhood $x\in U\sub V$ so that
$\fX_{f_t}\cap U$ is independent of $t\in [0,1]$. 
\end{prop}

\begin{proof} We consider $F: \CC\times V\to \CC$ defined by $F(t,z)=(1-t)f_0(z)+tf_1(z)$,
and define the relative critical locus $\fX_{F/\CC}$ to be defined by the ideal $(d_V F)$, where $d_V$ is the
differential along the $V$ directions. Because $\fX_{f_0}=\fX_{f_1}$, 
\beq\label{3.18}
\CC\times \fX_{f_0}\sub \fX_{F/\CC}\sub \CC\times V
\eeq
are closed analytic subspaces. We prove that there is an open subset $x\in U\sub V$ and $[0,1]\sub W\sub \CC$
so that 
\beq\label{3.19}
W\times\fX_{f_0}=\fX_{F/\CC}\cap (W\times U).
\eeq

Let $\pi: \CC\times V\to V$ be the projection. Because $\pi\sta ( df_0)\sub (d_V F)$ is a coherent subsheaf of $\sO_V$-modules,
the quotient $(d_V F)/\pi\sta (df_0)$ is a coherent sheaf of $\sO_V$-modules. Thus
$$\Sigma=\{z\in \CC\times V\mid (d_V F)/\pi\sta (df_0)\big\vert_z\ne 0\}
$$
is a closed analytic subset of $V$. So if $\Sigma\cap (\CC\times x)=\emptyset$, we can find open neighborhoods $x\in U\sub V$ and 
$[0,1]\sub W\sub \CC$ such that \eqref{3.19} holds.

Let $t\in\CC$ and let $\hat \fX_{f_t}$ be the formal completion of $\fX_{f_t}$ at $x$. Since $t\times_\CC \fX_{F/\CC}=\fX_{f_t}$,
$(t,x)\not\in\Sigma$ if and only if $\hat\fX_{f_t}=\hat \fX_{f_0}$. We now prove this identity.
Because $\hat\fX_{f_0}\sub \hat\fX_{f_t}$, we have a surjective homomorphism $\sO_{\hat \fX_{f_t}}\twoheadrightarrow \sO_{\hat \fX_{f_0}}$.
We will show $\sO_{\hat \fX_{f_t}}=\sO_{ \hat \fX_{f_0}}$
by showing that for any $k\ge 1$, 
\beq\label{3.20}
\sO_{\fX_{f_t}}\big/ \bl \sO_{ \fX_{f_t}}\cap \fm_x^{k}\br =\sO_{\fX_{f_0}}
\big/ \bl \sO_{ \fX_{f_0}}\cap \fm_x^{k}\br. %\sub \sO_{ V}\big/ \bl \sO_{ V_x}\cap \fm_x^{k}\br.
\eeq
The identity for $k=2$ follows from $f_1-f_0\in \fm_x^3$, thus $T_x\fX_{f_t}=T_x\fX_{f_0}$. 
Suppose $\sO_{ \hat\fX_{f_t}}\neq \sO_{ \hat \fX_{f_0}}$. Then there is a $k_0\ge 2$ so that \eqref{3.20} is true 
for $k=k_0$ but not for $k=k_0+1$.

We let $B=\sO_{\fX_{f_t}}/(\sO_{ \fX_{f_t}}\cap\fm_x^{k_0+1})$, and let $I\sub B$ be the kernel of the 
restriction homomorphism 
$$B\to \bar B=\sO_{\fX_{f_t}}/(\sO_{ \fX_{f_t}}\cap\fm_x^{k_0}).
$$
We let $\bar g: \spec \bar B\to \fX_{f_t}$ be the tautological morphism. Because 
the identity \eqref{3.20} holds for $k_0$, the composite $\spec \bar B\to \fX_{f_t}\to V$ factors through
$\bar{g}': \spec \bar B\to \fX_{f_0}$.

By the definition of $(I,B, \bar g)$, $\bar g$ extends to $\spec B\to \fX_{f_t}$. Thus
the obstruction class $ob_{\fX_{f_t}}(\bar g,B,\bar B)$ is zero. 
We claim that the obstruction class $ob_{\fX_{f_0}}(\bar g',B,\bar B)$ to extending $\bar{g}'$
to $\spec B\to \fX_{f_0}$ is trivial too. Indeed, 
Using the identity $\hat\fX_{f_1}=\hat\fX_{f_0}\sub \hat V$, we can view $\bar{g}'$ as a
morphism from $\spec \bar B$ to $\fX_{f_1}$. 
Because $f_0-f_1\in\fm_x^3$, we have
\beq\label{3.21}
\Omega_{\fX_{f_0}}|_x=\Omega_{\fX_{f_1}}|_x=\Omega_{\fX_{f_t}}|_x,
\eeq
as quotient spaces of $\Omega_V|_x$. Adding that $f_t=(1-t)f_0+t f_1$, we obtain
\beq\label{3.22}
ob_{\fX_{f_t}}(\bar g,B,\bar B)=(1-t)\cdot ob_{\fX_{f_0}}(\bar{g}',B,\bar B)+t\cdot ob_{\fX_{f_1}}(\bar{g}',B,\bar B).
\eeq
Here we can equate and add because they are elements in $I\otimes_\CC\Omega_{\fX_{f_0}}|_x$.

Because $\fX_{f_0}$ and $\fX_{f_1}$ have identical
obstruction assignments at $x$,  we have
$$ob_{\fX_{f_0}}(\bar{g}',B,\bar B)=ob_{\fX_{f_1}}(\bar{g}',B,\bar B).
$$
Using \eqref{3.22} and $ob_{\fX_{f_t}}(\bar{g},B,\bar B)=0$, we get
$ob_{\fX_{f_0}}(\bar{g}',B,\bar B)=0$. Hence
$\bar\varphi'$ also extends to $\varphi': \spec B\to \fX_{f_0}$. 
Composing with the ring homomorphism induced by $\hat\fX_{f_0}\sub\hat\fX_{f_t}$, we obtain a composite
ring homomorphism
\beq\label{3.23}
\sO_{ \fX_{f_t}}\big/ \bl \sO_{ \fX_{f_t}}\cap \fm_x^{k_0+1}\br
\lra \sO_{ \fX_{f_0}}\big/ \bl \sO_{ \fX_{f_0}}\cap \fm_x^{k_0+1}\br\mapright{\varphi^{\prime\ast}}
\sO_{ \fX_{f_t}}\big/ \bl \sO_{ \fX_{f_t}}\cap \fm_x^{k_0+1}\br
\eeq
whose restriction to $\sO_{ \fX_{f_t}}/ ( \sO_{ \fX_{f_t}}\cap \fm_x^{2})$ is the identity map. Thus the composite 
\eqref{3.23} is an isomorphism. Because $\hat\fX_{f_0}\sub \hat \fX_{f_t}$, we conclude that
\eqref{3.20} holds for $k=k_0+1$, a contradiction. This proves $\hat \fX_{f_t}=\hat \fX_{f_0}$,
and the proposition.
\end{proof}

\begin{rema}
The proposition may fail when $f_0-f_1\not\in\fm_x^3$. For example, let $V=\CC$, $f_0=z^2$ and $f_1=-z^2$.
Then $\fX_{f_0}=\fX_{f_1}$ while $\fX_{f_{1/2}}\ne \fX_{f_0}$. 
\end{rema}

\begin{coro}\label{3.24}
Let the situation be as in Proposition \ref{3.17}. Then there is an open neighborhood $x\in U\sub V$
such that the family of holomorphic functions $f_t$ induces an isomorphism
$$\tau_{01}: P_{f_0}|_{X_{f_0}\cap U}\cong P_{f_1}|_{X_{f_1}\cap U}.
$$
\end{coro}

\begin{proof}
Applying Lemma \ref{3.12} to the family $f_t$ gives the corollary. 
\end{proof}

\bigskip

\section{Proof of Theorem \ref{2.20}}\label{S3.4}

In this section, we complete our proof of Theorem \ref{2.20}. 
\begin{lemm}\label{3.27}
Let $\Psi:\CC^r\to \CC^r$ be the linear isomorphism defined by the diagonal matrix $\mathrm{diag}(\pm 1,1,\cdots,1)$. Let $q(z_1,\cdots,z_r)=z_1^2+\cdots+z_r^2$. Then the pullback isomorphism
$P_q\lra P_q$
induced from $\Psi$ by Proposition \ref{2.7} is $\pm \mathrm{id}$. Moreover if $\Phi:V\times\CC^r\to V\times \CC^r$ is $\mathrm{id}\times \Psi$ and $f(z,y)=g(y)+q(z)$ for $y\in V$ and $z\in \CC^r$, then the pullback isomorphism induced from $\Phi$ by Proposition \ref{2.7} is $\pm \mathrm{id}:P_f\to P_f$. 
\end{lemm}
\begin{proof}
The Milnor fiber at $0$ is homotopic to the sphere $S^{r-1}$ whose reduced cohomology is $\tilde{H}^{r-1}(S^{r-1})=\QQ$. Obviously, $\Psi$ acts $\tilde{H}^{r-1}(S^{r-1})$ on as $\pm 1$. For the second statement, use Proposition \ref{3.1}.
\end{proof}

\begin{theo}[Theorem \ref{2.20}]\label{3.25}
Let $(V,f)$ be an LG pair. 
Let $U$ be an open neighborhood of $X_f$ in $V$ and $\Phi:U\to V$ be a holomorphic map, biholomorphic
onto its image, such that $f\circ\Phi=f|_U$ and $\Phi|_{\fX_f}=\mathrm{id}_{\fX_f}$. Then the isomorphism $\sigma:P_f\to P_f$ induced from $\Phi$ by Proposition \ref{2.7}  is $\det (d\Phi|_{\fX_f})\cdot \mathrm{id}$. \end{theo}
Note that by Lemma \ref{1.14}, $\det (d\Phi|_{\fX_f})$ is locally constant with values in $\{\pm 1\}$. 

\begin{proof}
Since this is a local problem, for any $x\in\fX_f$, by shrinking $V$, we can 
assume that $V$ embeds into $\CC^n$ as an open subset, so that and $x\in V$ corresponds to $0\in \CC^n$. 
Applying Lemma \ref{3.15}, we can assume that for the standard complex coordinate variables
$(y_1,\cdots, y_n)$ of $\CC^n$, % let $f: V\to\CC$ take the form
\beq\label{3.26}
f=u\cdot \bl y_1^2+\cdots+y^2_r+h(y_{r+1},\cdots,y_n)\br,
\eeq
where $u:V\to\CC$ is a holomorphic function such that $T_0\fX_f$ is the 
linear space $V_0\defeq \{y_1=\cdots=y_r=0\}\cap V$ and $u|_{V_0}=1$.

Using $V\sub\CC^n$ and the canonical isomorphism $T_x\CC^n\cong \CC^n$, 
we form the linear transformation
$$\Psi\defeq d\Phi|_0:\CC^n\lra\CC^n.
$$
Because $\Phi|_{\fX_f}=\id_{\fX_f}$, and because $T_x\fX_f$ is the linear subspace
$\{y_1=\cdots=y_r=0\}$, we see that $\Psi(V_0)\sub V_0$ and $\Psi|_{V_0}=\id_{V_0}$. Because of \eqref{3.26}, possibly after shrinking $x\in V_0$, we have $\fX_f\sub V_0$. Thus $\Psi|_{\fX_f}=\id_{\fX_f}$.

For $s\in [0,1]$, we define
$$\varphi_s=(1-s)\cdot \Phi+s\cdot \Psi: U\lra \CC^n.
$$
Since $d\varphi_s|_0=\Psi$ are invertible, by shrinking $U$ if necessary,
$\varphi_s$ map into $V\sub \CC^n$, and are biholomorphic onto their images.
Note that $\varphi_0=\Phi|_U$ and $\varphi_1=\Psi|_U$ is linear.
Because both $\Phi|_{\fX_f}=\Psi|_{\fX_f}=\id_{\fX_f}$, we have $\varphi_s|_{\fX_f\cap U}=\id_{\fX_f\cap U}$.

We now construct various isomorphisms. 
We let $g=u\upmo f$, and let $f_t$, $t\in [0,1]$, be $u_t\cdot g|_U$, as in the proof of Lemma \ref{3.14}. 
Thus $f_t$ interpolates between $f|_U$ and $g|_U$, with $f_0=f|_U$ and $f_1=g|_U$.
We then form the composite isomorphism
$$\eta_{s,t}: P_{f_t}|_{X_{f_{t}}\cap U}\mapright{\cong} P_{f_t\circ\varphi_s}|_{X_{f_t\circ\varphi_s}}\mapright{\cong} 
P_{f_t}|_{X_{f_t}\cap U},
$$
where the first isomorphism is obtained by applying Proposition \ref{2.7} to the map $\varphi_s$ from
$(U,f_t)$ to $(U, f_t\circ\varphi_s)$, and
the second isomorphism is obtained by applying Corollary \ref{3.24} to $f_t\circ\varphi_s-f_t\in\fm_x^3$. 
Applying Corollary \ref{3.24} to $f_t$, possibly after shrinking $x\in U$,
we have a family of isomorphisms $\tau_{t}$ as shown: %=P_{f}|_{X_f\cap U}.
$$\begin{CD}
P_{f_t}|_{X_{f_{t}}\cap U}@>{\tau_t}>> P_{f_0}|_{X_{f_{0}}\cap U}\\
@V{\eta_{s,t}}VV @V{\eta_{0,0}}VV\\
P_{f_t}|_{X_{f_{t}}\cap U}@>{\tau_t}>> P_{f_0}|_{X_{f_{0}}\cap U}.
\end{CD}
$$
Since for $(s,t)=(0,0)$, the square is commutative, and since $\tau_t$ and $\eta_{s,t}$ are continuous
families of isomorphisms, applying Lemma \ref{1.14}, the above square is commutative at $(s,t)=(1,1)$.

Since $\Phi|_{V_0}=\mathrm{id}_{V_0}$, $\Psi$ is of the form
$$\left(\begin{matrix}
A&0\\
B&I
\end{matrix}\right)$$
with $A\in O(r)$. Since $O(r)$ has only two connected components, we can find a continuous path $A_t$  in $O(r)$ from $A$ to the diagonal $r\times r$ matrix $\diag(e,1\cdots,1)$ where $e=\det A=\det\, d\Phi|_0$. Let $B_t=(1-t)B$. So that we have a path $\Psi_t$ from $\Psi$ above to the $n\times n$ diagonal matrix $D:=\diag(e,1,\cdots,1)$. Then $g\circ \Psi_t-g\in \fm_x^3$ for all $t$ because $h\in \fm_x^3$. Applying the same argument as above, we obtain a commutative diagram of perverse sheaves
\[\xymatrix{
P_g\ar[r]^{\tau'} \ar[d]_{D} & P_g\ar[r]^{\tau_1} \ar[d]^{\eta_{1,1}} & P_f\ar[d]^{\eta_{0,0}}\\
P_g\ar[r]^{\tau'} & P_g\ar[r]^{\tau_1} & P_f
,}\]
all restricted to $U$. By Lemma \ref{3.27}, the left vertical is $e\cdot\mathrm{id}$. The right vertical is the isomorphism $\sigma$ in Theorem \ref{3.25}. Since the horizontal maps are the same, we find that $\sigma=e\cdot\mathrm{id}$.
%By construction, we have $\eta_{0,0}=\Phi_\ast$ and $\eta_{1,1}=\Psi_\ast$. 
%Applying the Sebastiani-Thom isomorphism (Proposition \ref{3.1}), we conclude that $\eta_{1,1}=\Psi_\ast=e\cdot\id$, where $e$ is the determinant of
%$d\Phi|_x$ in the statement of the theorem. Thus the commutativity of the square at $(1,1)$ implies that
%the isomorphism $\sigma$ equals
%$\eta_{0,0}=\tau_1\circ e\cdot \id\circ\tau_1\upmo=e\cdot\id$ with $e=\det\, d\Phi|_0$. 
%This proves the theorem.
\end{proof}

%%%%%%%%%%%%

\def\SH{\mathfrak{Sh} _{si}}
\def\End{\mathrm{End}\, }
\part{Critical virtual manifolds and moduli of sheaves}

In Part 2, we prove that moduli spaces $X$ of simple sheaves on a Calabi-Yau 3-fold $Y$ are critical virtual manifolds. Here a simple sheaf means a coherent sheaf whose automorphism group is $\CC^*\cdot \id$. By the Seidel-Thomas twists, we only
need to consider moduli spaces $X$ of simple holomorphic vector bundles on $Y$ (cf. Proposition \ref{6.4} and Theorem \ref{6.8}).  

Via gauge theory, a moduli space $X$ of simple holomorphic vector bundles on $Y$ can be thought of as an open analytic space in the quotient $\fV_{si}=\cA^{int}_{si}/\cG_k$ of the space $\cA^{int}_{si}$ of simple integrable semiconnections on a fixed complex vector bundle $E$ by the action of the gauge group $\cG_k$. Here a semiconnection means a $\dbar$-operator $$\dbar:\Omega^{0,0}(\End E)\lra \Omega^{0,1}(\End E)$$ which is integrable if its $(0,2)$-curvature $F^{0,2}_\dbar=\dbar\circ\dbar$ is zero. The space $\cA_{si}$ of simple semiconnections comes equipped with a function $CS$, called a holomorphic Chern-Simons functional, whose critical locus $\zero(d\,CS)$ coincides with the space $\cA^{int}_{si}$ of simple integrable semiconnections. Using the theory of elliptic partial differential operators, Joyce and Song proved that Miyajima's local section $V_0$ 
of $\cA_{si}\to \cA_{si}/\cG_k= \cB_{si}$ is a finite dimensional complex manifold such that the critical locus of the holomorphic function $f_0=CS|_{V_0}$ is an open set in $\fV_{si}$. Unfortunately, the dimension of $V_0$ varies from point to point and hence the Miyajima-Joyce-Song charts $V_0$ do not give us a critical virtual manifold structure on an open subspace $X\sub\fV_{si}$.

The idea to overcome this issue is to enlarge the Miyajima-Joyce-Song chart $V_0$, in a controlled fashion. We introduce the notion of a CS framing, which will be an infinitesimal guide in increasing $V_0$ to a larger complex manifold $V$ in $\cA_{si}$ such that the critical locus of $f=CS|_V$ stays the same as that of $f_0$. We modify Miyajima's elliptic differential operator so that the CS framing $\Xi$ can enlarge the solution space. We prove that in this way we obtain a complex submanifold $V$ of $\cA_{si}$ containing $V_0$, whose dimension $r$ is the fixed number
we choose from the beginning. We call such a complex manifold a \emph{CS chart}. Therefore we have an open cover $X=\cup_\lambda X_\lambda$ such that each $X_\lambda$ is the critical locus of a holomorphic function $f_\lambda$ on an $r$-dimensional complex manifold $V_\lambda$ (cf. Theorem \ref{4.28} and Corollary \ref{4.39}).

To get a critical virtual manifold structure on $X$, we further need to compare these local charts $\{V_\lambda\}$.
This is attained by a topological argument (cf. Theorem \ref{4.45}) together with the observations that
\begin{enumerate}
\item the CS chart $V$ centered at $x\in \cA^{int}_{si}$ and determined by a hermitian metric $h$ on $E$ and a CS framing $\Xi$ is actually independent of the choice of $(x,h,\Xi)$ up to a local equivalence (cf. Proposition \ref{4.41});
\item for any point $y$ sufficiently close to the center $x$ of a CS chart $V$, the CS chart $V'$ centered at $y$ is locally equivalent to $V$ (cf. Proposition \ref{4.42}). 
\end{enumerate}
We thus find that $X$ is a critical virtual manifold (Theorem \ref{4.45}).

The orientability of the critical virtual manifold $X$ is proved to be equivalent to the existence of a square root of the determinant line bundle 
$\det Rp_*R\cH om(\cE,\cE)$ on $X\ured$ where $\cE$ is the (local) tautological family of simple sheaves on $Y$ (cf. Theorem \ref{6.20}).
Then by an argument of Okounkov's (Proposition \ref{6.40}), such a square root exists if there is a global tautological family on $X\ured$. In particular, when $X$ is the moduli space of stable sheaves $\cF$ on $Y$ whose Hilbert polynomial is $\chi(\cF(m))=dm+e$ with $d,e$ coprime, then $X$ has an orientable critical virtual manifold structure and hence we have a perverse sheaf $P$ (resp. a polarizable mixed Hodge module $\cM$) on $X$, which is the gluing of the local perverse sheaves (resp. mixed Hodge modules) of vanishing cycles (cf. Theorem \ref{6.55}).

Finally we use the perverse sheaf $P$ to obtain a mathematical theory of the Gopakumar-Vafa invariant. When $X$ is a moduli space of 1-dimensional stable sheaves, we have morphisms 
$$X\lra S\lra \mathrm{pt}$$
where the first arrow comes from the Hilbert-Chow morphism, taking the supports of stable sheaves. By the relative hard Lefschetz property and the decomposition theorem for the hypercohomology of perverse sheaves underlying polarizable Hodge modules, we find that the hypercohomology $\bbH^*(X,\hat{P})$ of the perverse sheaf $\hat{P}$ is an $sl_2\times sl_2$-representation space where $\hat{P}=\mathrm{gr}^WP$ is the gradation of $P$ by the weight filtration. Hence we can decompose $\bbH^*(X,\hat{P})$ as the direct sum of copies of the cohomology $H^*(Jac_g)$ of $g$-dimensional complex tori for various $g$. If we imagine that the space of curves in $X$ is discrete, it is natural to define the number $n_g$ of copies of $H^*(Jac_g)$ in $\bbH^*(X,\hat{P})$ as a curve counting invariant. This is Gopakumar-Vafa's proposal which is now rigorous by our results. We conjecture that the Gopakumar-Vafa invariant carries all the information about the Gromov-Witten invariant of $Y$ (cf. Conjecture \ref{7.5}). We verify this conjecture for the genus 0 case (Proposition \ref{7.4}) and for K3-fibered Calabi-Yau 3-folds (Proposition \ref{7.10}). 

\bigskip

The layout of Part 2 is as follows. In Chapter \ref{ch4}, we use gauge theory to prove the existence of CS charts after proving the existence of CS framings. In Chapter \ref{ch5}, we prove two key results, Propositions \ref{4.41} and \ref{4.42}, that enable us to prove the existence of a  critical virtual manifold structure. In Chapter \ref{ch6}, we show that the Seidel-Thomas twists reduce the general case to the vector bundle case and prove the orientability criterion.
In Chapter \ref{ch7}, we apply the theory of critical virtual manifold to give a mathematical theory of the Gopakumar-Vafa invariant.

\chapter{Critical loci and moduli of simple vector bundles}\label{ch4}

Let $Y$ be a smooth projective Calabi-Yau 3-fold over $\CC$. In particular its canonical line bundle 
$K_Y=\Omega^3_Y$ is trivial.
Let $\Omega\in H^0(Y,\Omega^3_Y)$ be a fixed nowhere vanishing holomorphic 3-form on $Y$. 
Simple sheaves of fixed total Chern class $c$ on $Y$ form a stack $\mathfrak{Sh} _{si}^c$ 
and simple locally free sheaves form an open substack $\mathfrak{V} _{si}^c$. 
In this chapter, we prove that an open analytic subspace $X$ of $\mathfrak{V} _{si}^c$ 
has an open covering $\cup X\lalp$ such that $X\lalp$ is the critical locus $\Crit(f\lalp)$ 
for an LG pair $f\lalp:V\lalp\to \CC$ (cf. Theorem \ref{4.28}).  Then we prove that these charts give us a critical virtual manifold structure on $X$
assuming two propositions (Propositions \ref{4.41} and \ref{4.42}) that will be proved in Chapter \ref{ch5} (cf. Theorem \ref{4.45}). 

In Chapter \ref{ch6}, we will prove that any open bounded analytic space $X\subset \mathfrak{Sh} _{si}^c$ is isomorphic to an open $X'\subset \mathfrak{V} _{si}^{c'}$ for some $c'$ 
and hence admits a critical virtual manifold structure.

\section{Semiconnections and Chern-Simons functional}\label{S4.1}

In this section we recall necessary gauge theoretic background. 
Our presentation follows \cite{Miya, Koba}, and \cite[Chapter 9]{JoSo}.

\subsection{Gauge theory of semiconnections}\label{S4.1.1}

Let $E$ be a complex vector bundle on $Y$ of total Chern class $c$. We fix a K\"ahler metric on the Calabi-Yau 3-fold $Y$ once and for all. We let $h$ be a hermitian metric on $E$. Let 
\[ \Omega^{0,q}(E)=C^\infty(E\otimes_\CC \wedge^{0,q}T_Y^*) \]
denote the sheaf of $C^\infty$-sections of the vector bundle $E\otimes \wedge^{0,q}T_Y^*$. 
\begin{defi}\label{4.1} A smooth \emph{semiconnection} is a $\CC$-linear operator 
\beq\label{4.2} \dbar:\Omega^{0,0}(E)\lra \Omega^{0,1}(E)\eeq
satisfying the Leibniz rule: $\dbar(fs)=f\dbar s +s\otimes \dbar f$ for all smooth sections $s$ of $E$ and all smooth functions $f$ on $Y$, where $\dbar f$ is the usual $(0,1)$-part of the differential $df$. 
\end{defi}

A semiconnection $\dbar:\Omega^{0,0}(E)\to \Omega^{0,1}(E)$ defines an operator
\beq\label{4.3} \dbar: \Omega^{0,q}(E)\lra \Omega^{0,q+1}(E),\quad 
\dbar(s\otimes\alpha)=s\otimes \dbar \alpha +  \dbar s\wedge \alpha
\eeq
for all smooth sections $s$ of $E$ and all $(0,q)$-forms $\alpha$. The \emph{curvature} of a semiconnection $\dbar$ is defined by
\beq\label{4.4}
F^{0,2}_\dbar=\dbar^2:\Omega^{0,0}(E)\mapright{\dbar} \Omega^{0,1}(E)\mapright{\dbar} \Omega^{0,2}(E)
\eeq 
which is an element of $\Omega^{0,2}(\End E)$. 

\begin{defi}\label{4.5} 
A smooth semiconnection $\dbar$ is called \emph{integrable} or a \emph{holomorphic structure} on $E$ if $F^{0,2}_\dbar=\dbar^2=0$. \end{defi}

By the Nirenberg-Newlander theorem, an integrable semiconnection $\dbar$ defines the locally free sheaf 
(holomorphic vector bundle) $\cE$ of $\dbar$-flat sections by assigning to each open $U$, the vector space $$\cE(U)=\{s\in C^\infty(E)\,|\,\dbar s=0\}.
$$
 
Two such smooth semiconnections define the same holomorphic vector bundle $\cE$ if and only if they lie in the same orbit by the action of the gauge group 
$$\cG=C^\infty(\mathrm{Aut} (E))$$
of smooth sections of the principal bundle $\mathrm{Aut}(E)$ associated to $E$, where the gauge action is defined by $\dbar\cdot g=g^{-1}\circ \dbar\circ g$.
For an integrable connection $\dbar$ and its associated holomorphic vector bundle $\cE$, the stabilizer of $\dbar$ in $\cG$ is canonically isomorphic to the automorphism group of $\cE$. We say a semiconnection $\dbar$ or a holomorphic vector bundle $\cE$ is \emph{simple} if the stabilizer $\cG_\dbar$ or the automorphism group of $\cE$ is $\CC\cdot \mathrm{id}$. 

To apply the theory of elliptic differential operators, we complete the semiconnection space and the gauge group. We fix a pair of integers $k\ge 4$ and $\ell>6$. Using the hermitian metric $h$ of $E$, we form the completions
\[ \Omega^{0,q}(\End E)_k \and \cG_{k}\]
under the Sobolev norm $L^\ell_k$ and $L^\ell_{k+1}$ respectively where the $L^\ell_k$ norm is the sum of $L^\ell$-norms of the $s$th derivatives for $s\le k$. We let $\cA$ denote the space of $L^\ell_k$-semiconnections on $E$; let $\cA _{si}$ be the set of simple semiconnections in $\cA$; let $\cA _{si} ^{int}$ be the set of all simple integrable smooth semiconnections on $E$.

For a smooth semiconnection $\dbar$, we have a bijection
\beq\label{4.6}
\dbar+\cdot: \Omega^{0,1}(\End E)_k \lra \cA,\quad \fa\mapsto \dbar+\fa
\eeq
because $\dbar+\fa$ is a semiconnection and the difference of two semiconnections lies in $\Omega^{0,1}(\End E)_k$. 
Hence $\cA$ is an affine space based on $\Omega^{0,1}(\End E)_k$. 
Let 
\beq\label{4.7}
\cB _{si}=\cA _{si}/\cG_{k}\, \sub\, \cA/\cG_k=\cB
\eeq
and let 
\beq\label{4.8} \Pii:\cA\lra \cB\eeq
denote the quotient map. 
Then $\cB _{si}$ is a complex Banach manifold and the complex analytic space of gauge equivalence classes of simple integrable semiconnections
\beq\label{4.9} 
\fV _{si}^c=\cA _{si} ^{int}/\cG_k \subset \cB _{si}
\eeq
is the moduli space of simple holomorphic vector bundles on $Y$ of total Chern class $c$. 
%Therefore, 
%$$\fV _{si}^c=\cA _{si} ^{int}/\cG_k\subset \cB _{si}
%$$
%{\red which is isomorphic to the algebraic space of the moduli of simple holomorphic vector bundles whose underlying smooth bundle is isomorphic to $E$.}
We comment that the $L_k^\ell$ and $L^\ell_{k+1}$ structures on $\cA$ and $\cG_k$ do
not depend on the choice of $h$ and the reference smooth semiconnection $\dbar$. 

The operator \eqref{4.2} extends to the completions 
$$\dbar: \Omega^{0,q}(\End E)_{k-q+1}\to  \Omega^{0,q+1}(\End E)_{k-q}
$$
and we can form the Laplacian 
\beq\label{4.10}
\Lap_{\dbar}=\dbar\dbar\stax+\dbar\stax\dbar: \Omega^{0,q}(\End E)_{k-q+1}\to \Omega^{0,q}(\End E)_{k-q-1}.
\eeq
We denote by $\Lap_\dbar\upmo(0)^{0,q}$ the set of harmonic forms in $\Omega^{0,q}(\End E)$.\vsp
In Chapter \ref{ch5}, we will also use the truncated eigenspaces (for $q=1,2$)
\beq\label{4.11} \Theta^{0,q}_\dbar(\eps)\defeq \CC\text{-span}\,\{\fa\in \Omega^{0,q}(\End E)_{k-q+1}\mid \dbar\sta \fa=0,\  
\Lap_{\dbar}\fa=\lambda \fa, \  \lambda< \eps\}\sub \ker(\dbar\sta)^{0,q}
\eeq
for $\eps>0$.

\medskip

\subsection{Chern-Simons functionals}\label{S4.1.2}
For an integrable smooth semiconnection $\dbar\in \cA _{si} ^{int}$, by \eqref{4.6}, an arbitrary element of $\cA$ is
$\dbar +\fa$ for $\fa\in \Omega^{0,1}(End\,E)_k$.  %We write $\dbar_\fa=\dbar +\fa$.
Then the action of an element $g\in \cG_k$ of the gauge group is 
\[ (\dbar+\fa)\cdot g= \dbar + g^{-1}dg+g^{-1}\fa g \]
and the curvature of $\dbar +\fa$ is 
\beq\label{4.12} F^{0,2}_{\dbar+\fa} = (\dbar+\fa)^2=\bpar \fa+\fa\wedge\fa.\eeq
By the Sobolev inequality, $F^{0,2}_{\dbar+\cdot}$
is a continuous operator from $\Omega^{0,1}(\End E)_k$ to $\Omega^{0,2}(\End E)_{k-1}$, analytic in $\fa$. 

Recall that $\Omega$ is a non-trivial holomorphic (3,0)-form on $Y$ that we fixed once and for all.
The holomorphic Chern-Simons functional with the reference point $\dbar\in\cA ^{int} _{si}$ is 
\beq\label{4.13}
CS(\dbar+\fa)=\frac1{4\pi^2}\int_Y \tr \left(\frac12 (\bpar\fa)\wedge\fa + \frac13 \fa\wedge\fa\wedge \fa\right)\wedge \Omega,
\quad \dbar+\fa\in \cA.
\eeq 
It is a cubic polynomial in $\fa\in\Omega^{0,1}(\End E)_k$ whose quadratic part is 
\beq\label{4.14}
CS_{2}(\dbar+\fa)=\frac1{4\pi^2}\int_Y \tr \bigl(\frac12 (\bpar\fa)\wedge\fa\br \wedge \Omega,
\quad \fa\in\Omega^{0,1}(\End E)_k.
\eeq
Since the directional derivative of $CS$ at $\fa$ in the direction of $\fb$ is 
\beq\label{4.15}
\delta\,CS(\dbar+\fa)(\fb)=\frac1{4\pi^2}\int_Y \tr (\fb\wedge F^{0,2}_{\bpar+\fa})\wedge\Omega,
\eeq
$\delta\,CS(\dbar+\fa)=0$ if and only if $\dbar+\fa$ is integrable. Thus the complex analytic subspace
$\cA _{si} ^{int}$  of simple integrable smooth semiconnections in $\cA_{si}$ is the critical locus  
of $CS$. 

\begin{rema}\label{4.16} 
The Chern-Simons functional $CS$ depends on the reference semiconnection $\dbar$. Using \eqref{4.15}, we see that
for different reference points $\dbar$ and $\dbar'\in\cA ^{int} _{si}$, 
the Chern-Simons functionals \eqref{4.13} may differ only by a constant, 
and are equal when $\dbar$ and $\dbar'$ lie
in the same connected component of $\cA ^{int} _{si}$ or in the same orbit under the action of the gauge group $\cG_k$. 

For our purpose of constructing a critical virtual manifold structure on 
an open analytic subspace $X\sub \fV^c_{si}=\cA^{int}_{si}/\cG_k$, 
we may assume that $X$ is connected because we can construct a 
critical virtual manifold structure component by component. 
%Also adding a constant to CS won't make any difference.
Hence any two reference points $\dbar\in  \Pii^{-1}(X)\sub \cA^{int}_{si}$ via the quotient map $\Pii:\cA\to \cA/\cG_k=\cB$ give 
the same Chern-Simons functional. So in what follows, we won't mention the reference point for the Chern-Simons functionals. 
% and use only the Chern-Simons functional $CS$ with the reference point $\dbar_0$. 
%In the following, we will call the $CS$ referenced at any $\dbar\in\cA ^{int} _{si}$ a Chern-Simons functional. 
\end{rema}

\bigskip

\section{CS framings and CS charts} \label{S4.2}

We use the notation introduced in \S\ref{S4.1}: $E$ is a fixed complex vector bundle of total Chern class $c$ and we fix a hermitian metric $h$ on $E$. $\cA$ is the $L^\ell_k$-completion of the space of smooth semiconnections and $\cG_k$ is the $L^\ell_{k+1}$-completion of the gauge group, etc. 

For psychological reasons, we often write $x\in \cA _{si} ^{int}$ for a simple integrable semiconnection $\dbar$. 
Let $x=\dbar\in \cA _{si} ^{int}$  and denote by $\cE_x$ (with subscript $x$) the holomorphic vector bundle over the Calabi-Yau 3-fold $Y$ defined by $\dbar$. We pick an integer
$$r\ge \dim Ext^1_Y(\cE_x,\cE_x).
$$

\begin{defi}\label{4.19}
A complex submanifold $V$ of $\cA _{si}$ of dimension $r$ 
is called a  \emph{CS chart of dimension $r$ at $x$} if 
$x\in V$ and the critical locus $\Crit(CS|_V)$ of the restriction $CS|_V$ of the Chern-Simons functional to $V$ is isomorphic to an open neighborhood of $\bar{x}=\Pii(x)\in \fV _{si}^c$, by the restriction $\Pii|_{\Crit(CS|_V)}$ of the quotient map $\Pii:\cA\to \cB$ (cf. \eqref{4.8} and \eqref{4.9}).
\end{defi}

Our goal in this section is to prove that there is a CS chart at every $x\in \cA _{si} ^{int}$ of any given dimension $r$. 
\begin{theo}\label{4.22}
For each point $x=\dbar\in \cA _{si} ^{int}$ and any integer $r\ge \dim Ext^1_Y(\cE_x,\cE_x)$ where $\cE_x$ is the holomorphic vector bundle defined by $x$, there is a CS chart $V$ of dimension $r$ at $x$. 
\end{theo}
The rest of this section is devoted to a proof of Theorem \ref{4.22}.

\medskip

\subsection{The Miyajima-Joyce-Song chart}

A prototype of a CS chart was constructed by Miyajima and Joyce-Song 
(cf. \cite[\S1]{Miya} and \cite[Theorem 5.5]{JoSo}).
\begin{prop}\label{4.20} \cite{Miya, JoSo}
For a sufficiently small $\eps>0$, 
\beq\label{4.21}
V_{0}=\{\dbar+\fa\,|\, |\!|\fa|\!|<\varepsilon,\ \dbar^*\fa=0,\ \dbar^*(\dbar \fa+\fa\wedge \fa)=0\}
\eeq
is a finite dimensional complex submanifold of $\cA$ with $x=\dbar\in V_{0} $ and $T_xV_{0} =\Lap^{-1}_\dbar(0)^{0,1}\cong Ext^1(\cE_x,\cE_x)$  
where $\cE_x$ is the holomorphic vector bundle defined by $x=\dbar\in \cA _{si} ^{int}$. 
Moreover, $V_{0} $ is a CS chart of dimension $\dim Ext^1_Y(\cE_x,\cE_x)$ at $x$.
\end{prop}
We call $V_{0} $ the \emph{Miyajima-Joyce-Song chart} at $x$. 
The two equations $\dbar^*\fa=0, \dbar^*(\dbar \fa+\fa\wedge \fa)=0$ in \eqref{4.21} combine to give us the single equation
\beq\label{4.18}
\Lap_\dbar(\fa)+\bpar^*(\fa\wedge\fa)=\bpar(\bpar^*\fa)+\bpar^*(\bpar \fa+ \fa\wedge\fa)=0.
\eeq
In fact, this single equation is equivalent to the system of two equations in \eqref{4.21}: 
By the Hodge decomposition, \eqref{4.18} is equivalent to the system of equations 
$\bpar(\bpar^*\fa)=0$ and $\bpar^*(\bpar \fa+ \fa\wedge\fa)=0$, the former of which is obviously equivalent to
$\bpar^*\fa=0$.

The idea for a proof of Theorem \ref{4.22} is to enlarge the Miyajima-Joyce-Song chart by taking the inverse image of a finite dimensional subspace by the operator  $\Lap_\dbar(\fa)+\bpar^*(\fa\wedge\fa)$ in \eqref{4.18}. 

\medskip

\subsection{CS framings}

Let $x=\dbar\in \cA ^{int} _{si}$. By Remark \ref{4.16}, we may use any reference point for our purpose of constructing an LG pair whose critical locus is an open neighborhood of $\bar{x}=\Pii(x)\in \fV^c _{si}$. So we use the Chern-Simons functional $CS$ with the reference point $x=\dbar$. 
Since $\Lap_x \upmo(0)^{0,1}=\ker(\Lap_x)^{0,1}$ lies in the null subspace of the quadratic form $CS_{2}$
(cf. \eqref{4.14}), $CS_{2}$ induces a bilinear pairing
\beq\label{4.23}
\cQ : T_{x}\cA/\bigtriangleup_x\upmo(0)^{0,1} \times T_{x}\cA/\Lap_x\upmo(0)^{0,1} \lra \CC,\eeq
defined by 
$$\cQ (\fa,\fb)=\frac{1}{8\pi^2}\int_Y\mathrm{tr}(\dbar\fa\wedge\fb)\wedge \Omega$$
which is symmetric by Stokes' theorem. 
Here we use the canonical isomorphism $T_x\cA\cong \OadE$ induced by the isomorphism \eqref{4.6}.

\begin{defi}\label{4.24} An $r$-dimensional \emph{CS framing}  at $x\in \cA ^{int} _{si}$ is an $r$-dimensional subspace 
$\Xi \sub \ker (\dbar^*)_k^{0,1}\sub \OadE= T_x\cA$ 
such that
\begin{enumerate}
\item $\Lap_x\upmo(0)^{0,1}\sub \Xi $;
\item the restriction $\cQ|_{\Xi /\Lap_x\upmo(0)^{0,1}}$ is a nondegenerate quadratic form.
\end{enumerate}
\end{defi}

\begin{lemm}\label{4.29}
For any $x=\dbar\in \cA ^{int} _{si}$ and $r\ge \dim Ext^1_Y(\cE_x,\cE_x)=\dim \Lap_x^{-1}(0)^{0,1}$, there is an $r$-dimensional CS framing at $x$. % where $\cE$ is the holomorphic vector bundle on $Y$ defined by $x=\dbar$. 
\end{lemm}
\begin{proof}
Note that 
%the pairing $\cQ_x(\dbar\mathfrak{c},\fa)=0$ for any $\fa$. Hence tangent vectors $\dbar \mathfrak{c}$ to the orbit $\cG x$ of $x=\dbar$  lie in the kernel of the quadratic form $\cQ_x$.On the other hand, 
$\cQ $ is nondegenerate on $\ker (\dbar^*)_{k}^{0,1}/\Lap_x^{-1}(0)^{0,1}$. Indeed, if $\cQ (\fa,\fb)=0$ for all $\fb$, 
then $\dbar\fa=0$. Hence $\fa\in  \Lap_x^{-1}(0)^{0,1}$ since $\dbar^*\fa=0$. By the standard argument in linear algebra on diagonalizing symmetric operators, we can find a subspace $N$ of dimension $r-\dim \Lap_x^{-1}(0)^{0,1}$ of $\ker (\dbar^*)_{k}^{0,1}/\Lap_x^{-1}(0)^{0,1}$ which is nondegenerate with respect to $\cQ $. Then $\Xi =N\oplus \Lap_x^{-1}(0)^{0,1}\subset \ker (\dbar^*)_{k}^{0,1}$ is an $r$-dimensional CS framing at $x$. 
\end{proof}

\medskip

\subsection{CS charts from CS framings}
Let $x=\dbar\in \cA _{si} ^{int}$ and let $\Lap_x=\dbar\dbar^*+\dbar^*\dbar$.
For a CS framing $\Xi \subset \ker (\dbar^*)_{k}^{0,1}\subset \Omega^{0,1}(\End E)_k$, let 
\beq\label{4.25}
\tilde{\Xi} :=\Lap_x^{-1}(0)^{0,1}\oplus \Lap_x(\Xi )
\eeq
and define the quotient homomorphism
\beq\label{4.26a}
P_x: \Omega^{0,1}(\End E)_{k-2} \lra  \Omega^{0,1}(\End E)_{k-2}\big/\tilde{\Xi} ,
\eeq
and the elliptic operator
\beq\label{4.26}
\bL_x: \Omega^{0,1}( \End E)_k
\lra \Omega^{0,1}(\End E)_{k-2}/\tilde{\Xi},\quad  \bL_x(\fa)=P_x \bl \Lap_x \fa+\dbar\sta \fa\wedge\fa\br.
\eeq

\begin{theo}\label{4.28} Let $x=\dbar\in\cA ^{int} _{si}$, and let $\Xi $ be an $r$-dimensional CS framing at $x$.
For a sufficiently small $\veps >0$, 
\beq\label{4.27}
V =\{ \dbar+\fa\mid \fa\in \Omega^{0,1}(\End E)_k,\ \bL_x(\fa)=0, \
\parallel\! \fa\!\parallel_k<\varepsilon \}. 
\eeq
is an $r$-dimensional CS chart at $x$.
\end{theo}
The CS chart $V$ clearly depends on the choice of the point $x\in \cA^{int}_{si}$, 
the hermitian metric $h$ on $E$ and the CS framing $\Xi$ at $x$. 
We will write $V=V(x,h,\Xi)$ if we want to emphasize the choice of $(x,h,\Xi)$.

\begin{lemm}\label{4.34}
$\bL_x(\fa)=0$ if and only if 
\beq\label{4.31}
\dbar\sta\fa=0, \and  P_x\circ \dbar\stax (\dbar \fa+\fa\wedge\fa)=0.
\eeq
\end{lemm}
\begin{proof}
It is immediate to see that \eqref{4.31} implies $\bL_x(\fa)=0$.
For the other direction, suppose $\bL_x(\fa)=0$. Since $\Xi \sub \ker(\dbar\sta)_{k}^{0,1}$,
$\tilde{\Xi} \sub \ker(\dbar\sta)_{k-2}^{0,1}$. Applying $\dbar\sta$ to
$\bL_x(\fa)=0$, we obtain $\dbar^*\dbar\dbar\sta\fa=0$, which forces $\dbar\sta\fa=0$. 
Having this, we obtain $P_x\circ \dbar\stax (\dbar \fa+\fa\wedge\fa)=0$. 
\end{proof}
\begin{lemm}\label{4.35}
The tangent space $T_xV $ to $V $ at $x$ is $\Xi $.
\end{lemm}
\begin{proof}
By Lemma \ref{4.34}, letting $\fa=0$, $T_xV $ is defined by 
$$\dbar^*\fb=0\and \dbar^*\dbar\fb\in \Lap_x(\Xi ).$$
Since $\Xi \subset \ker(\dbar\sta)_{k}^{0,1}$, we find that $T_xV =\Xi $.
\end{proof}

\begin{lemm}\label{4.36} 
For a CS framing $\Xi \sub \ker(\dbar\sta)_{k}^{0,1}$ at $x=\dbar\in\cA _{si} ^{int}$, let 
$$\hat{\Xi}= \Lap_x^{-1}(0)^{0,2}\oplus \dbar(\Xi ).$$
Then the bilinear pairing 
\beq\label{4.37}
\langle \cdot,\cdot\rangle: \Omega^{0,1}(\End E)_k\times \Omega^{0,2}(\End E)_{k-1}\lra \CC,
\eeq
defined by
\[ \langle \fa,\fb\rangle=\frac{1}{8\pi^2}\int_Y \tr(\fa\wedge\fb)\wedge\Omega\]
is nondegenerate on $\Xi \times \hat{\Xi} $.
\end{lemm}
The pairing $\langle \cdot,\cdot\rangle$ relates to the symmetric quadratic form $CS_{2}$ defined by \eqref{4.14}  via
\beq\label{4.38}
CS_{2}(\fa,\fb)=\langle\fa,\dbar \fb\rangle, \quad \fa,\fb\in \OadE. 
\eeq
\begin{proof}[Proof of Lemma \ref{4.36}]
By the Serre duality theorem, the pairing $\langle \cdot,\cdot\rangle$ is nondegenerate on $\Lap_x^{-1}(0)^{0,1}\times \Lap_x^{-1}(0)^{0,2}$. 
Let $e_{1},\cdots, e_{l}\in \Xi $ be such that
$\dbar e_{1},\cdots,$ $\dbar e_{l}$ form a basis for $\dbar (\Xi )$.
By Hodge theory, $e_{1},\cdots, e_{l}$ and $\Lap_x\upmo(0)^{0,1}$ span $\Xi $.
Because $\Lap_x\upmo(0)^{0,1}$ is orthogonal to $\image(\dbar)$ under $\langle \cdot,\cdot\rangle$, by \eqref{4.38}, 
$\cQ $ is nondegenerate on $\Xi /\Lap_x\upmo(0)^{0,1}$ 
if and only if $\langle e_{i},\dbar  e_{j}\rangle$ form an invertible $l\times l$ matrix,
which is equivalent to saying that $\langle\cdot,\cdot\rangle$ on $\Xi \times \hat{\Xi}$ is nondegenerate. 
This proves the lemma. 
\end{proof}

\begin{proof}[Proof of Theorem \ref{4.28}]
By Proposition \ref{4.20},
$(d(CS|_{V_0})=0)\sub V_0$ is an open neighborhoof of $\bar x\in \fV^c_{si}$.
By the proof of \cite[Proposition 9.12]{JoSo},  
$$(F^{0,2}=0)\cap V_0= (d(CS|_{V_0})=0).$$
Further, using the definition of $V$ in \eqref{4.27}, we have the identity 
$ (F^{0,2}=0)\cap V=(F^{0,2}=0)\cap V_0$. 
Thus to prove the theorem, it suffices to
show that
\beq\label{AA}
(d(CS|_V)=0)= (F^{0,2}=0)\cap V.
\eeq

Inspired by the proof of \cite[Proposition 9.12]{JoSo}, we define a subbundle 
\beq\label{4.30}
R\defeq \{(\dbar+\fa,\fb)\in V \times\Omega^{0,2}(\End E)_k \mid P_x\circ\dbar\sta\fb=\dbar\sta(\dbar\fb-\fb\wedge\fa-\fa\wedge\fb)=0\}.
\eeq
Because the two equations in the bracket are holomorphic in $\fa$, 
for $\veps$ small enough, applying the implicit function theorem, we see
that $R$ is a holomorphic subbundle of $V \times\Omega^{0,2}(\End E)_k$ over $V $. 
We think of the curvature $F^{0,2}_{\dbar+\fa}=\dbar\fa+\fa\wedge\fa$ as a section $F^{0,2}|_{V }$ of the trivial bundle 
$V \times\Omega^{0,2}(\End E)_k\to V $.
Then the Bianchi identity coupled with the equations \eqref{4.31} 
ensure that the restriction of the curvature section $F^{0,2}$ to $V $ is a section of $R $. 

We then define a bundle map
\beq\label{4.32} \varphi : R \to T^* V , \quad
(\dbar+\fa,\fb)\mapsto (\dbar+\fa, \alpha_\fb),
\eeq
where
$\alpha_\fb\in T_{\dbar +\fa}^* V $ is defined as
$$\alpha_\fb(\cdot)=\frac{1}{8\pi^2}\int_Y \tr(\cdot\wedge\fb)\wedge\Omega=\langle\cdot,\fb\rangle$$ using the bilinear pairing \eqref{4.37}.
Clearly, $\varphi $ is holomorphic and $\varphi \circ F^{0,2}|_{V }=df $ where $f =CS|_{V }$.
Therefore, if $\varphi $ is an isomorphism of vector bundles, 
$$(df=0)=(\varphi\circ F^{0,2}|_{V}=0)=(F^{0,2}|_{V}=0).
$$
This will prove the theorem.

We show that by choosing $\veps$ small, we can make $\varphi $ an isomorphism of vector bundles over
$V $. We claim that restricting to $x=\dbar \in V $, we have
\beq\label{4.33}
R |_{x}=\{\fb\in\Omega^{0,2}(\End E)_{k} \mid P_x\circ\dbar\sta\fb=\dbar\sta\dbar\fb=0\}=\dbar(\Xi )\oplus \Lap_x\upmo(0)^{0,2}.
\eeq
Indeed the first equality follows from the definition of $R $.
We prove the second equality. The inclusion $\supset$ is obvious. To prove $\subset$, we let $\fb\in R |_{x}$.  
Since $\dbar\sta\dbar\fb=0$, we have $\dbar\fb=0$. 
Thus we can write $\fb=\fb_0+\dbar\fc$
with $\Lap_x\fb_0=0$. Since $P_x\circ\dbar\sta\fb=0$ and $\Xi \sub \ker (\dbar^*)_k^{0,1}$, $\dbar^*\fb=\dbar^*\dbar\fc\in  \Lap_x(\Xi )=\dbar^*\dbar(\Xi )$. 
Hence we may take $\fc\in \Xi $.
This proves \eqref{4.33}.

By Lemmas \ref{4.35} and \ref{4.36}, $R |_x=\hat{\Xi} $ is dual to $T_xV =\Xi $ by the bilinear form \eqref{4.37}. Therefore, $\varphi |_x$ is an isomorphism. Since being an isomorphism is an open condition, after shrinking $\varepsilon$ if necessary, we find that $\varphi $ is an isomorphism as desired. This proves the theorem.
\end{proof}

\begin{coro}\label{4.39}
Let $X\subset \fV^c _{si}=\cA^{int}_{si}/\cG_k$ be an analytic space which is open in $\fV^c _{si}$. Suppose there is an integer $r$ such that  $r\ge \dim Ext^1_Y(\cE_x,\cE_x)$ for all $x\in X$.
%holomorphic vector bundles $\cE$ on $Y$ parameterized by $X$. 
We have an open cover $X=\cup_{\lambda\in \Lambda} X_\lambda$ and $r$-dimensional CS charts 
$$X_\lambda\mapleft{\cong} \Crit(f_\lambda)\subset V_\lambda\mapright{f_\lambda}\CC$$
where each $V_\lambda\subset \cA$ is a CS chart and $f_\lambda=CS|_{V_\lambda}$ is the restriction of the Chern-Simons functional.
\end{coro}  
\begin{proof}
Let $\bar{x}\in X$ and ${x}\in \cA _{si} ^{int}$ be a point over $\bar{x}$, i.e. $\Pii(x)=\bar{x}$ via the quotient map $\Pii:\cA\to \cB=\cA/\cG_k$. By Lemma \ref{4.29}, we have an $r$-dimensional CS framing $\Xi$ at $x$. The CS chart $V$ constructed in Theorem \ref{4.28} has the critical locus $\Crit(f)$ of $f=CS|_V$ which is isomorphic to an open neighborhood $\bar{x}$ in $\fV^c _{si}$ by $\Pii$. Varying $\bar{x}\in X$, we get the desired open cover by the critical loci of CS charts.
\end{proof}

\bigskip

\section{Critical virtual manifold structure} \label{S4.3}

What else do we need to get a critical virtual manifold structure on an open analytic subspace $X$ of $\fV^c _{si}=\cA^{int}_{si}/\cG$, on top of Theorem \ref{4.28} and Corollary \ref{4.39}? As we have an open covering by critical loci of CS charts, we further need transition maps comparing CS charts.

First we note that the definition of the CS chart $V$ over a point $\bar{x}\in X\subset \fV^c _{si}$  in Theorem \ref{4.28} relied on the choice of a lift $x\in \cA _{si} ^{int}$, a CS framing $\Xi$ and a hermitian metric $h$ on $E$. 
We first prove that our CS chart $V$ constructed in Theorem \ref{4.28} is independent of the choice of $(x,h,\Xi)$, up to a local equivalence.

\begin{defi}\label{4.40} 
Given two LG pairs $(V,f)$ and $(V',f')$ and $x\in \Crit(f)$, 
we say $(V,f)$ is \emph{locally equivalent at} $x$ to $(V',f')$
if there is an open neighborhood $V^\circ$ of $x$  in $V$ and a holomorphic map $\Phi:V^\circ\to V'$, biholomorphic onto an open set ${V'}^\circ$ in $V'$, such that $f'\circ\Phi=f|_{V\circ}$. 
\end{defi}
Recall that an LG pair $(V,f)$ is a holomorphic function $f$ on a complex manifold $V$ which has only one critical value $0$.

\begin{prop}\label{4.41} Let $x$ and $x'$ be points in $\cA _{si} ^{int}$ with $\Pii(x)=\Pii(x')$ where $\Pii:\cA\to \cB=\cA/\cG_k$ is the quotient map. Let $h$ and $h'$ be two hermitian metrics on $E$. Let $\Xi$ and $\Xi'$ be CS framings at $x$ and $x'$ respectively. Let $CS$ be the Chern-Simons functional.
Let $V=V(x,h,\Xi)$ and $V'=V(x',h',\Xi')$ be the CS charts constructed in Theorem \ref{4.28} with data $(x,h,\Xi)$ and $(x',h',\Xi')$ respectively. If we let $f=CS|_V$ and $f'=CS|_{V'}$, then $(V,f)$ is locally equivalent to $(V',f')$ at $x$.
\end{prop}
We will prove Proposition \ref{4.41} in Chapter \ref{ch5}. See Remark \ref{4.16} for the reference point of the Chern-Simons functional.  

Moreover, we have the following comparison result.

\begin{prop}\label{4.42} 
Let $x=\dbar_x\in \cA _{si} ^{int}$. Let $h$ be a hermitian metric on $E$ and $\Xi$ 
be a CS framing at $x$.
%$r\ge \dim Ext^1_Y(\cE,\cE)$ where $\cE$ is the holomorphic vector bundle on the Calabi-Yau 3-fold $Y$ defined by $\dbar_x$. 
Let $V=V(x,h,\Xi)$ be the CS chart constructed in Theorem \ref{4.28} with sufficiently small $\varepsilon$. Let $f=CS|_V$.
Then for each $y\in \Crit(f)$, there is a CS framing $\Xi'$ at $y$ such that for the CS chart $V'=V(y,h,\Xi')$ with sufficiently small $\varepsilon$ and $f'=CS|_{V'}$, $(V,f)$ is locally equivalent at $y$ to $(V',f')$.
\end{prop}

We will prove Proposition \ref{4.42} in Chapter \ref{ch5}.

For psychological reasons, when we think of a point $x\in \cA _{si}$ as an operator $\Omega^{0,0}(\End E)_k\to \Omega^{0,1}(\End E)_{k-1}$,  we will denote this operator by $\dbar_x$ as in Proposition \ref{4.42}.

From Propositions \ref{4.41} and \ref{4.42}, we obtain the following.
\begin{prop}\label{4.43} 
Let $x, x'\in \cA^{int}_{si}$ such that $\Pii(x)=\bar{x}, \Pii(x')=\bar{x}'\in \fV_{si}^c=\cA^{int}_{si}/\cG$. 
Let $\Xi$ and $\Xi'$ be CS framings at $x$ and $x'$ respectively of the same dimension $r$. Let $h$ and $h'$ be hermitian metrics on $E$.
Let $V=V(x,h,\Xi)$ and $V'=V(x',h',\Xi')$ be the CS charts constructed in Theorem \ref{4.28} 
(for sufficiently small $\varepsilon$ so that Proposition \ref{4.42} holds for $V$ and $V'$) 
with respect to $(x,h,\Xi)$ and $(x',h',\Xi')$ respectively. 
Let $f=CS|_V$ and $f'=CS|_{V'}$ denote the restrictions of the Chern-Simons functional, so that $\Crit(f)$ and $\Crit(f')$ give open neighborhoods 
$U=\Pii(\Crit(f))$ and $U'=\Pii(\Crit(f'))$ of $\bar{x}$ and $\bar{x}'$ in $\fV^c_{si}$ respectively. 
Suppose we have $y\in \Crit(f)$ and $y'\in \Crit(f')$ such that 
$\Pii(y)=\Pii(y')=\bar{y}\in U\cap U'\sub \fV^c_{si}$. Then $(V,f)$ is locally equivalent at $y$ to $(V',f')$. 
\end{prop}
\begin{proof}
By Proposition \ref{4.42}, $(V,f)$ is locally equivalent to $(V_y,CS|_{V_y})$ at $y$ where 
$V_y=V(y,h,\Xi_y)$ for some CS framing $\Xi_y$ at $y$. 
Likewise, $(V',f')$ is locally equivalent to $(V_{y'},CS|_{V_{y'}})$ at $y'$
where $V_{y'}=V(y',h',\Xi_{y'})$ for some CS framing $\Xi_{y'}$ at $y'$. By Proposition \ref{4.41}, $(V_y,CS|_{V_y})$ is locally equivalent to $(V_{y'},CS|_{V_{y'}})$ at $y$. Therefore $(V,f)$ is locally equivalent to $(V',f')$ at $y$.  
\end{proof}

Proposition \ref{4.43} is sufficient for the existence of a critical virtual manifold structure because of the following.
\begin{prop}\label{4.44}
Let $X$ be an analytic space and $X=\cup_{\lambda\in \Lambda} U_\lambda$ be a locally finite open cover. Let $(V_\lambda,f_\lambda)$ be an LG pair whose critical locus $\Crit(f_\lambda)$ is isomorphic to $U_\lambda$. Suppose for any $y\in U_\lambda\cap U_{\lambda'}$, $(V_\lambda,f_\lambda)$ is locally equivalent to $(V_{\lambda'},f_{\lambda'})$ at $y$. Then there is a critical virtual manifold structure on $X$, given by restricting the charts $(V_\lambda,f_\lambda)$.
\end{prop}
We will often use the following simple observation: Let $(V,f)$ be an LG pair and $U\sub \Crit(f)$ be an open subset. Then there is an open subset $V^\circ\sub V$ such that $\Crit(f|_{V^\circ})=U$. For instance, we may let $V^\circ=V-(\Crit(f)-U)$. 
\begin{proof}
The analytic space $X$ is a metric space by Urysohn's metrization theorem 
because $X$ is a Hausdorff paracompact second countable topological space. 
We fix a metric $d$ on $X$. 

For $x\in X$, let 
$$\Lambda_x=\{\lambda\in \Lambda\,|\, U_\lambda\ni x\}.$$
Since the open cover is locally finite, $\Lambda_x$ is finite. 
By assumption, there exists $\varepsilon_x>0$ such that for any pair $\lambda,\lambda'\in \Lambda_x$, we have an equivalence 
$$\Phi_{\lambda\lambda'}:V^\circ_\lambda\lra V_{\lambda'}^\circ$$
of LG pairs $(V_\lambda^\circ,f_\lambda|_{V_\lambda^\circ})$ and $(V_{\lambda'}^\circ,f_{\lambda'}|_{V_{\lambda'}^\circ})$
for some open neighborhoods $V^\circ_\lambda\sub V_\lambda$, ${V}^\circ_{\lambda'}\sub V_{\lambda'}$ of the ball 
$$B(x,3\varepsilon_x)=\{x'\in X\,|\,d(x,x')<3\varepsilon_x\}.$$
In particular, $B(x,3\varepsilon_x)\subset \cap_{\lambda\in \Lambda_x}U_\lambda$.

Fix a map $\lambda:X\to \Lambda$ with $x\in U_{\lambda(x)}$ or $\lambda(x)\in \Lambda_x$. Then restricting $(V_{\lambda(x)},f_{\lambda(x)})$ to an open neighborhood of 
$$U_x:=B(x,\varepsilon_x)$$
gives us a chart $f_x:V_x\to \CC$ with $\Crit(f_x)=U_x$. 

We claim that $X=\cup_{x\in X} U_x$ together with charts
$$U_x=\Crit(f_x)\hookrightarrow V_x\mapright{f_x}\CC$$
form a critical virtual manifold. We have to check that there is an equivalence 
$$\varphi_{xx'}:V_x^\circ\lra V_{x'}^\circ,\quad \text{with } f_{x'}\circ \varphi_{xx'}=f_x|_{V_x^\circ}$$ 
of open neighborhoods $V_x^\circ\subset V_x\sub V_{\lambda(x)}$ 
and $V_{x'}^\circ\sub V_{x'}\sub V_{\lambda(x')}$ of $U_x\cap U_{x'}$.
Indeed, if $U_x\cap U_{x'}\ne \emptyset$ with $\varepsilon_x\ge \epsilon_{x'}$, then 
$$U_x\cup U_{x'}\sub B(x,3\varepsilon_x)\sub \cap_{\lambda\in\Lambda_x}U_\lambda$$
and hence $\Lambda_x\subset \Lambda_{x'}$. By construction, for $\lambda=\lambda(x),\lambda'=\lambda(x')\in \Lambda_{x'}$, we have an equivalence
$$\Phi_{\lambda\lambda'}:V_{\lambda}^\circ\lra V_{\lambda'}^\circ$$
of open neighborhoods $V_\lambda^\circ\sub V_\lambda$, $V_{\lambda'}^\circ\sub V_{\lambda'}$ of $U_{x'}$, which certainly restricts to an equivalence
$$\varphi_{xx'}=\Phi_{\lambda\lambda'}|_{V_{xx'}}:V_{xx'}\lra V_{x'x}$$
of open neighborhoods $V_{xx'}\sub V_\lambda^\circ\sub V_{\lambda(x)}$, $V_{x'x}\sub V_{\lambda'}^\circ\sub V_{\lambda(x')}$ of $U_x\cap U_{x'}$. We let $\varphi_{x'x}=\varphi_{xx'}^{-1}$. This proves the proposition.
\end{proof}

By Proposition \ref{4.43}, we have an open cover $\{\Crit(f)\subset V(x,h,\Xi)\}_{x\in X}$ satisfying the assumption of Proposition \ref{4.44}. Since $X$ is paracompact by assumption, we can find a locally finite subcover. 
Therefore from Propositions \ref{4.43} and \ref{4.44}, we obtain the following, because $T_xX=Ext^1_Y(\cE_x,\cE_x)$.
% where $\cE$ is the holomorphic vector bundle defined by $x\in \cA^{int}_{si}$.
\begin{theo}\label{4.45}
Let $X\subset \fV^c_{si}=\cA^{int}_{si}/\cG_k$ be an open analytic space. Suppose there is an integer $r$ such that $r\ge \dim 
T_xX$ for all $x\in X$.
%Ext^1_Y(\cE,\cE)$ for all holomorphic vector bundles $\cE$ on $Y$ parameterized by $X$. 
Then there is a critical virtual manifold structure on $X$ given by the CS charts from Theorem \ref{4.28}.
\end{theo}

\begin{rema}\label{4.46} 
The existence of such an integer $r$ in Theorem \ref{4.45} is guaranteed if for instance $X$ is (quasi-)compact or quasi-projective or $X$ is a Noetherian scheme. % by the upper semicontinuity of $\dim Ext^1_Y(\cE,\cE)$ . 
\end{rema}

So it only remains to prove Propositions \ref{4.41} and \ref{4.42}.
Chapter \ref{ch5} is devoted to their proofs.

%%%%%%%%%%%%%
%%%%%%%%%%%%%
%%%%%%%%%%%%%
%%%%%%%%%%%%%
%%%%%%%%%%%%%

\chapter{Critical virtual manifolds and moduli of vector bundles}\label{ch5}

In this chapter we prove Propositions \ref{4.41} and \ref{4.42} to complete our proof of Theorem \ref{4.45}. 
The common nature of the two theorems is that the CS charts we constructed in Theorem \ref{4.28} 
are locally equivalent when we perturb the data $(x,h,\Xi)$. The main point of our proof is that we can construct
analytic families of CS charts which are locally trivial. By using vector field arguments, 
we can find local trivializations which preserve the Chern-Simons functional.

\bigskip

\section{Proof of Proposition \ref{4.42}}\label{S5.1}

In this section, we prove Proposition \ref{4.42}.

Let $x=\dbar_x\in \cA^{int}_{si}$ be an integrable simple semiconnection on the complex vector bundle $E$ over the Calabi-Yau 3-fold $Y$. 
Let $\cE_x$ denote the holomorphic vector bundle on $Y$ determined by $x=\dbar_x$ and let
\beq\label{5.3} d_x=\dim Ext^1_Y(\cE_x,\cE_x)=\dim T_{\bar x}\fV_{si}^c\eeq
where $\bar x=\Pii(x)\in \cA^{int}_{si}/\cG_k\sub \cB$ and $\Pii:\cA\to \cB$ is the quotient map by the gauge group action.

Fix an integer $r\ge d_x$ and a hermitian metric $h$ on $E$. Let $\Xi_x$ be an $r$-dimensional CS framing at $x$ and 
let $V_x=V(x,h,\Xi_x)$ be the CS chart constructed in Theorem \ref{4.28} with respect to a sufficiently small $\varepsilon>0$. 
Let $f_x=CS|_{V_x}$ be the restriction of the Chern-Simons functional and let $U_x=\Crit(f_x)$. 
Then $\Pii(U_x)$ is an open neighborhood of $\bar x=\pi(x)$ in $\fV_{si}^c=\cA^{int}_{si}/\cG_k$.
Notice that we have added the subscript $x$ to emphasize the ``center'' of the CS chart.
We have to show the following:

\medskip

\noindent\emph{After shrinking $U_x$ if necessary, for any $y\in U_x$, there is a CS framing $\Xi_y$ at $y$ such that  
the CS chart $V_y=V(y,h,\Xi_y)$ together with $f_y=CS|_{V_y}$ is locally equivalent to $(V_x,f_x)$ at $y$.}

\bigskip

\subsection{A family CS framing}\label{S5.1.1} 
We first introduce the notion of a family CS framing.
\begin{defi}\label{5.0} A \emph{family CS framing} on $U_x$ extending $\Xi_x$ is a rank $r$ analytic subbundle 
\beq\label{5.1} \Xi\sub U_x\times \Omega^{0,1}(\End E)_k\eeq
of the trivial bundle $U_x\times \Omega^{0,1}(\End E)_k$ over $U$ such that
\begin{enumerate}
\item for each $y\in U_x$, the fiber $\Xi_y:=\Xi|_y$ over $y$ is a CS framing at $y$;
\item $\Xi|_x=\Xi_x$. 
\end{enumerate}
\end{defi}
Here an analytic vector bundle on $U_x=\Crit(f_x)\sub V_x$ means a vector bundle on $U_x$ that locally is the restriction of a holomorphic vector bundle on a complexification of $V_x$.

%\medskip

%We have the following existence result.

\begin{prop} \label{5.2}
After shrinking $U_x$ if necessary, we can find a family CS framing $\Xi$ on $U_x$ extending $\Xi_x$. 
\end{prop}

Our proof of Proposition \ref{5.2} uses the truncated eigenspaces.
For $q=1,2$ and $y=\dbar_y\in V_x\subset \cA_{si}$, we form
$$\Lap_y=\dbar_y\dbar_y\sta+\dbar_y\sta\dbar_y: \Omega^{0,q}(\End E)_{k-q+1}\lra \Omega^{0,q}(\End E)_{k-q-1}.
$$
where the adjoint is defined using the fixed hermitian metric $h$ (and a fixed K\"ahler metric on $Y$). We let
\beq\label{5.6}
\Theta^{0,q}_y(\eps)= \CC\text{-span}\,\{\fa\in \Omega^{0,q}(\End E)_{k-q+1}\mid 
\Lap_{y}\fa=\lambda \fa, \  \lambda< \eps\}, %\sub \ker(\dbar_y\sta)^{0,q}.
\eeq
and let
\beq\label{5.7}
\Theta^{0,q}_{V_x}(\eps)=\coprod_{y\in V_x}\{y\}\times \Theta^{0,q}_{y}(\eps)\sub V_x\times \Omega^{0,q}(\End E)_{k-q+1},
\eeq
considered as a subbundle of the trivial bundle $V_x\times \Omega^{0,q}(\End E)_{k-q+1}$ over $V_x$.
Since $\Lap_x$ has only nonnegative discrete eigenvalues, there is an $\eps>0$ such that $\Lap_x$ has no eigenvalue in the open interval $(0,2\eps)$. 

\begin{lemm}\label{5.4} 
With $\Theta^{0,1}_{V_x}(\eps)$ and $\eps>0$ as above,
after shrinking $V_x$ if necessary, 
$\Theta^{0,1}_{V_x}(\eps)$ is a rank $d_x$ analytic subbundle over $V_x$. % as in \eqref{5.7}
\end{lemm}
 
We first recall the notion of complexification. 
For $\delta>0$, we let 
$$D_\delta=\{(u_1,\cdots,u_n)\in \RR^n\mid |u_1|<\delta,\cdots,|u_n|<\delta\}$$  
to be the $\delta$-polydisk. We let $\CC^n$ have coordinate variables $w_1,\cdots, w_n$ and $w_j=u_j+iv_j$.
We let $\RR^n\sub \CC^n$ be the embedding via $(u_j)\mapsto (u_j)$.
The complex polydisk $D_\delta^\CC\sub\CC^n$ is called a complexification of $D_\delta$.

Since $V_x$ is a smooth complex manifold, after shrinking $x\in V_x$ if necessary, we can make it
analytically isomorphic to $D_\delta$ with $x\mapsto 0$, for $n=2\dim V_x$. We then call $D_\delta^\CC$ 
the induced complexification of $V_x$, denoted by $V_x^\CC$. 

In the following, we use $u$ to mean both a point in $V_x$ and its coordinates in $\RR^n$
via $V_x\cong D_\delta\sub \RR^n$.
We denote $\dbar_{u}=\dbar_0+\fa(u)$ the associated semiconnection of $u\in V_x$.
Then $\fa(u)\in \OadE$ and is analytic in $u$.
Thus by choosing $\delta$ small, we
can assume that $\fa(u)$ extends to a family of
forms $\fa(w)$, holomorphic in $w$; i.e. a holomorphic $\fa(\cdot): D_\delta^\CC\to \OadE$ extending the given $\fa(u)$.

We consider the formal adjoints $\dbar_{u}\sta=\dbar_0\sta+\fa(u)\udag$ of $\dbar_u$. 
Note that $\fa(u)\udag$ is also analytic in $u$. 
We then choose $\delta$ small enough so that
$\fa(u)\udag$ extends to a holomorphic $\fa(\cdot)\udag :D^\CC_\delta\to \Omega^{0,1}(\End E)_k$.

We let 
$\dbar_w=\dbar_0+\fa(w)$, $\dbar\sta_w=\dbar_0\sta+\fa(w)\udag$ and form
$$\Lap_w=\dbar_w\dbar_w\sta+\dbar_w\sta\dbar_w: \Omega^{0,1}(\End E)_k\lra \Omega^{0,1}(\End E)_{k-2},
\quad w\in V_x^\CC.
$$
It is a family of second order elliptic operators, holomorphic in $w\in V_x^\CC$, whose symbols are
identical to that of $\Lap_{0}$. %Then Lemma \ref{cont} follows from the following lemma.

%\begin{lemm} \label{Lhes}
%Let the notation be as stated. Suppose $\eps_0>0$ separates eigenvalues of $\dbar_0$. % and $\eps_0<\eps(z)$ for all $z\in O_0$.
%Then for sufficiently small $\delta$, % $D^\CC\supset D$ such that
%$$\Theta^{0,1}_{Z_\delta,\eta}(\eps_0)=\coprod_{x\in Z_\delta} x\times \Theta^{0,1}_{\eta(x)}(\eps_0)\sub 
%Z_\delta\times \Omega^{0,1}(ad E)_s
%$$
%extends to a holomorphic subbundle of $D_\delta^\CC\times \Omega^{0,1}(ad E)_s$. 
%\end{lemm}

\begin{proof}[Proof of Lemma \ref{5.4}]
%We extend $\Theta_D(\eps_0)$ to $D^\CC$ using the generalized eigenforms of $\Lap_w$.
In the following, we are free to shrink $x\in V_x$ (and $V_x^\CC$) whenever necessary.
Since $\Lap_w$ is holomorphic in $w\in D_\delta^\CC$, and since $1+\Lap_w$, $w\in V_x^\CC$, are invertible,
by \cite[page 365]{Kato} the family
$$(1+\Lap_w)\upmo: \Omega^{0,1}(ad E)_s\lra \Omega^{0,1}(ad E)_s, \quad w\in D_\delta^\CC,
$$ 
is a holomorphic family of bounded operators. 

We now extend $\Theta_{V_x}^{0,1}(\eps)$. First, note that
$\lambda$ is a spectrum of $\Lap_w$ if and only if $(1+\lambda)\upmo$ is a spectrum of $(1+\Lap_w)\upmo$,
and they have identical associated spaces of generalized eigenforms. 
Since $\Lap_0$ has discrete spectrum (eigenvalues) and $\eps$ is not an eigenvalue,
by the continuity of the spectrum, we can find
sufficiently small $\delta$ and $\delta'>0$
so that no eigenvalues $\lambda$ of $(1+\Lap_w)\upmo$ for $w\in V_x^\CC$ lie in 
$\big||\lambda|-(1+\eps_0)\upmo\big|<\delta'$.

% we can find an open $D_x^\CC\subset D^\CC$,
%$D_x^\CC\cap D=D_x$,
%such that no $(1+\Lap_w)\upmo$, $w\in D_x^\CC$, contains spectrum in the region $\big||\lambda|-(1+\eps_0)\upmo\big|<\delta/2$.
Applying \cite[Theorem VII-1.7]{Kato}, %for $w\in V_x^\CC$ 
we have decompositions
$$V_x^\CC\times \Omega^{0,1}(\End E)_k=E_1\oplus E_{2}
$$
into holomorphic subbundles of $V_x^\CC\times \Omega^{0,1}(\End E)_k$ 
such that $E_1|_w$ and $E_2|_w$ are 
invariant under $(1+\Lap_w)\upmo$, and $T_{i,w}\defeq (1+\Lap_w)\upmo|_{E_i|_w}: E_i|_w\to E_i|_w$ 
has spectrum in $|\lambda|<(1+\eps_0)\upmo$ for $i=1$, and in $|\lambda|>(1+\eps_0)\upmo$ for $i=2$.

For us, the key property is that $E_{1,u}=\Theta^{0,1}_{u}(\eps)$ when $u\in V_x$. 
%We define $\Theta_w(\eps_0)=E_{1,w}$, for $w\in D_\delta^\CC$.
Thus $\Theta_{V_x}^{1,0}(\eps)$ is the restriction of the holomoprhic bundle $E_1$.
%$$\Theta^{0,1}_{D_\delta^\CC,\eta}(\eps_0)\defeq \coprod_{w\in D_\delta^\CC} w\times 
%\Theta^{0,1}_{\eta^\CC(w)}(\eps_0)
%\sub D_\delta^\CC\times\Omega^{0,1}(adE)_s.
%$$
This proves the lemma.
\end{proof}

Let $y\in\cA^{int}_{si}$. 
Recall that $\cQ_y$ is the descent of $CS_{2}$ referenced at $y$ (cf. \eqref{4.14})
to $T_y\cA/\Lap_y\upmo(0)^{0,1}$ (cf. \eqref{4.23}).
From \eqref{4.37}, we also have a bilinear pairing 
$$\langle \fa,\fb\rangle=\frac{1}{8\pi^2}\int_Y \tr(\fa\wedge\fb)\wedge\Omega.$$

\begin{lemm}\label{5.5} 
Let $y\in \cA^{int}_{si}$. 
Let $W\sub \ker(\dbar_y\sta)^{0,1}_k$ be a subspace containing $\Lap_y\upmo(0)^{0,1}$, 
and let 
$$\hat{W}=\Lap_y^{-1}(0)^{0,2}\oplus \dbar_y(W).$$ 
Then $\cQ_{y}|_{W/\Lap_y\upmo(0)^{0,1}}$ is nondegenerate if and only if
the restricted pairing
$$\langle \cdot,\cdot\rangle: W\times \hat{W}\lra \CC$$
is nondegenerate.
\end{lemm}

\begin{proof}
The proof is identical to that of Lemma \ref{4.36} with $\Xi$ replaced by $W$. 
\end{proof}

\def\im{\mathrm{im} }
\begin{proof}[Proof of Proposition \ref{5.2}]
Let $l=r-d_x$ and pick $e_1,\cdots,e_l\in \Xi_x\subset \Omega^{0,1}(\End E)_k$ whose images in $\Xi_x/\Lap_x^{-1}(0)^{0,1}$ form a basis. 
By the decomposition $$\Omega^{0,1}(\End E)_k=\im (\dbar_y)^{0,0}_{k+1}\oplus  \ker(\dbar_y^*)^{0,1}_k$$ for $y\in U_x$, we have a surjective homomorphism 
$$\Omega^{0,1}(\End E)_k\lra \ker(\dbar_y^*)^{0,1}_k, \quad \forall y\in U_x.$$ 
Let $\ti e_1(y),\cdots,\ti e_l(y)\in \ker(\dbar_y^*)^{0,1}_k$ be the images of $e_1,\cdots,e_l$. 

Since $\Theta^{0,1}_{U_x}(\eps)$ is an analytic vector bundle of rank $d_x$ by Lemma \ref{5.4} and because $\ti e_1(x)=e_1,\cdots,\ti e_l(x)=e_l\in \ker(\dbar_x^*)^{0,1}_k$ are linearly independent,  after shrinking $U_x$ if necessary, 
$\ti e_1,\cdots,\ti e_l$ and $\Theta_{U_x}^{0,1}(\eps)$ span
a rank $r$ analytic subbundle 
$$\Xi\defeq \text{Span}\bigl\{\ti e_1,\cdots,\ti e_l, \Theta_{U_x}^{0,1}(\eps)\bigr\}
\sub U_x\times\OadE.
$$
such that $\Xi_y:=\Xi|_y\sub \ker(\dbar_y^*)^{0,1}_k$ for $y\in U_x$. 

By Lemma \ref{5.5}, we only need to check that the pairing $\langle \cdot,\cdot\rangle: \Xi_y\times\hat{\Xi}_y\to\CC$
is nondegenerate. Because this pairing is nondegenerate when $y=x$, by shrinking $x\in U_x$ if necessary, 
it is nondegenerate for
all $y\in U_x$. Therefore, $\Xi$ is a family CS framing over $U_x$.
\end{proof}

\begin{rema}\label{remark5.5}
Let $V_x^\CC$ be a complexification of $V_x$, constructed before the proof of Lemma \ref{5.4}.
Then we can make the $\Xi$ constructed in the proof of Proposition \ref{5.2}, to be the restriction of a holomorphic subbundle $\Xi^\CC\sub
V_x^\CC\times \OadE$  by the same proof. 
\end{rema}

\medskip

\subsection{Comparing CS charts}\label{S5.1.2}

Fixing the complexification $V_x^\CC$ of $V_x$ and the holomorphic
bundle $\Xi^\CC$ mentioned in Remark \ref{remark5.5}, for $w\in V_x^\CC$, we let $\Xi_w=\Xi^\CC|_w$, form $\ti\Xi_w$ as in \eqref{4.25}, and form 
\beq\label{5.8} V_w=V(w,h,\Xi_w),\quad w\in V_x^\CC,\eeq
using Theorem \ref{4.28}. Namely, $V_w$ are defined as in \eqref{4.27} by
$$V_w=\{\dbar\,|\, \bL_w(\dbar-\dbar_w)=0,\ |\!|\dbar-\dbar_w|\!|<\varepsilon\}$$
where $\bL_w(\fa)=P_w(\Lap_w\fa+\dbar_w^*\fa\wedge\fa)$ and 
$$P_w:\Omega^{0,1}(\End E)_{k-2}\lra \Omega^{0,1}(\End E)_{k-2}/\ti \Xi_w$$
is the projection.
Note that by shrinking $V_x^\CC$ if necessary and choosing $\varepsilon$ sufficiently small, all $V_w$ are smooth
complex submanifolds of $\cA$, and
%. Note that for $y\in V_x$, $y\in V_y$ if and only if $y\in U_x$.
\beq\label{5.9}
\cV=\coprod_{w\in V^\CC_x} \{w\}\times V_w\sub V_x^\CC\times\cA %\mapright{pr} V_x^\CC,\quad (w,z)\mapsto z ,
\eeq
is a complex submanifold. 
%considered as an open subspace in $V^\CC_x\times\cA$. 
%(Without lose of generality, we can assume $\cV$ is an open subspace in $V^\CC_x\times\cA$.) 
We let $\pr_1:\cV\to V_x^\CC$ and $\pr_2:\cV\to\cA$ be its first and the second projection.
For $w\in V_x^\CC$, we denote $f_w=f|_{\cV_w}$, and denote $U_w=\Crit(f_w)\sub V_w$.
%In the following, we will use $y$ (resp. $w$) to denote a general member in $V_x$ (resp. $V_x^\CC$).
We will use $(w,w')$ to denote the element in $\cV$ whose projection to $V_x^\CC$ and $\cA$ are $w$ and $w'$,
respectively.

\begin{lemm}\label{5.9j}
There is an open neighborhood $x\in\cU\sub V_x^\CC$ so that for any $w\in \cU$,
$\pi|_{U_w}: U_w\to \cB$ factors through $\fV^c_{si}$ and gives an open embedding 
$\pi|_{U_w}: U_w\to \fV^c_{si}$ whose image contains $x$.
\end{lemm}

\begin{proof}
Let $F^{0,2}: \cA\to \Omega^{0,2}(\End E)_{k-2}$ be the curvature operator. We let
$$\cX=(F^{0,2}\circ\pi_2=0)\sub \cV
$$
be defined by the vanishing of $F^{0,2}\circ\pi_2$. For $y\in V_x$, we know that $\pi$ restricted to
$\cX_y=\cX\cap \{y\}\times V_y$ is an open embedding into $\fV^c_{si}$. We claim that the same
is true if we replace $y$ by $w\in V_x^\CC$, sufficiently close to $x$, and replace $\cV$ by an open neighborhood
$\cV^\circ\sub\cV$ of $(x,x)\in\cV$.

To this end, we let $\hat\cX_x$ be the formal completion of $\cX_x$ at $x\in \cX_x$; i.e. the germ of $\cX_x$
at $x$. We let $\fm_0$ be the ideal sheaf of $x\in V_x^\CC$ and let $\iota_n: z_n\to V_x^\CC$ be the
$n$-th neighborhood of $x\in V_x^\CC$, defined by the ideal sheaf $\fm_0^{n+1}$. We claim that
we can extend the tautological $\hat \cX_x\to\cX_x$ to a (holomorphic) map $j_n$ and a commutative diagram
\beq\label{sq1}
\begin{CD}
\hat \cX_x\times z_n @>{j_n}>> \cV\\
@VV{\pr_2}V @VVV\\
z_n @>{\iota_n}>> V_x^\CC
\end{CD}
\eeq
We prove the claim. As before, we let $(u_1,\cdots, u_l)$ be (real) analytic coordinates of $V_x$ at $x$, and let
$(w_1,\cdots,w_l)$ be its complexification, which become complex coordinates of $V_x^\CC$ at $x$. For an
$n\ge 1$, we can extend $\hat \cX_x\to \cX_x$($=U_x$) to
$$%\beq\label{sq2}
\begin{CD}
\hat \cX_x\times \spec \RR[u_1,\cdots,u_l]/(u_1,\cdots,u_l)^{n+1} @>>> \coprod_{y\in V_x}\{y\}\times V_y @>{\sub}>> \cV\\
@VV{\pr_2}V @VVV @VVV\\
\spec \RR[u_1,\cdots,u_l]/(u_1,\cdots,u_l)^{n+1} @>>> V_x @>>>V_x^\CC
\end{CD}
$$%\eeq
Since the second square is a Cartesian square and is the real part of $\cV\to X_x^\CC$, by the property of
complexification, we conclude that the desired extension \eqref{sq1} exists.

Because the local ring of analytic functions on $V_x^\CC$ is noetherian, we conclude that there
is an open $(x,x)\in \cV^\circ\sub\cV$ so that for all $w\in V_x^\CC$, letting $\cV_w^\circ=\cV^\circ\cap\{w\}\times V_w$,
considered as a subspace in $V_w$, we have
\beq\label{emb}
\pi|_{\cV_w^\circ}: \cX_w\cap \cV_w^\circ\lra \fV^c_{si}
\eeq
is an open embedding. Thus we can find an open $x\in \cU\sub V_x^\CC$ so that for every $w\in\cU$,
$\pi|_{\cV_w^\circ}(\cX_w\cap \cV_w^\circ)$ is an open neighborhood of $\pi(x)\in \fV^c_{si}$.

Finally, we need to verify that we can choose $\cV^\circ\sub \cV$ so that for any
$w\in\cU$, $\cX_w\cap \cV^\circ=U_w\cap\cV^\circ$. (Recall $U_w=\Crit(f_w)\sub \cV_w$.)
We let $f=CS\circ \pr_2: \cV\to \CC$, a holomorphic function; we let $d_{w'}$ be the relative
differential, relative to $V_x^\CC$. We let
 $\Crit_{V_x^\CC}(f)\sub \cV$
be defined by the ideal sheaf $(d_{w'} f)$. Because $dCS=F^{0,2}$, we conclude that 
$\Crit_{V_x^\CC}(f)\sub\cX$. On the other hand, since 
$$\Crit_{V_x^\CC}(f)\cap V_x=\cX\cap V_x,
$$
we conclude that after shrinking $(x,x)\in\cV^\circ\sub\cV$, we have
$$\Crit_{V_x^\CC}(f)\cap \cV^\circ=\cX\cap\cV^\circ.$$ 
This proves the lemma. 
\end{proof}

After shrinking $(x,x)\in \cV$ and $x\in V_x$, we can assume $\cV^\circ=\cV$ and $\cU=V_x^\CC$
satisfy the conclusions of the previous lemma.
%For $w\in V_x^\CC$, we denote $U_w=\Crit(f_w)$.

\begin{prop}\label{5.10}
After shrinking $V_x$ if necessary, we can find an open $\cV^\circ\sub \cV$ containing 
$(x,x)\in \cV$, and a holomorphic $\zeta:\cV^\circ\to V_x$ such that
$\zeta(y,y)=y$ whenever $y\in U_x$, and 
%for any $y\in U_x$, there is a biholomorphic map $\psi:V_x^\circ\to V_y^\circ$ of open neighborhoods $V_x^\circ\sub V_x$ and $V_y^\circ\sub V_y$ of $y$ such that 
$$\Pii\circ\zeta|_{\cV^\circ\cap \{w\}\times U_w}=\Pii\circ\pr_2|_{\cV^\circ\cap \{w\}\times U_w}: 
\cV^\circ\cap \{w\}\times U_w\to \fV_{si}^c,\quad \forall\, w\in V_x^\CC.$$
%for all $y\in U_x$,
%where $\Pii:\cA_{si}\to \cB_{si}$ is the quotient map by $\cG_k$.
%{\red and $f_x-f_y\circ\psi\in\fm_y^3$.} % and $\fV_{si}^c=\cA^{int}_{si}/\cG_k$.
\end{prop}

\begin{proof}

Let $\bar\cG=\cG_k/\CC\cdot\mathrm{id}$.
For $g\in \bar\cG$, and $(w,w')\in\cV$, we let
\beq\label{5.12}\bR_{(w,w')}: \bar\cG\lra \Omega^{0,0}(\End E)_{k-1}/\CC\eeq
be defined by
$$\bR_{(w,w')}(g)=\dbar_x\sta\bl \dbar_{w'}\cdot g-\dbar_x\br.$$
Then $\bR_{(x,x')}(\id)=0$ for $x'\in V_x$ by the definition of $V_x$.

We calculate the linearization of the operator $\bR_{(w,w')}$ at $(x,x)$ and $g=\id$:
$$\delta \bR_{(x,x)}|_{g=\id}=-\dbar_x\sta\dbar_x: 
\Omega^{0,0}(\End E)_{k+1}/\CC\lra \Omega^{0,0}(\End E)_{k-1}/\CC.
$$
Because $\dbar_x$ is simple, $\delta \bR_{(x,x)}|_{g=\id}$ is an isomorphism.
By the implicit function theorem, there is an open $\cV^\circ\sub \cV$,
%\coprod_{y\in U_x}\{y\}\times V_y\sub U_x\times\cA$, 
and a smooth map $g: \cV^\circ\to \bar\cG$ that is the unique solution to
$\bR_{(w,w')}(g(w,w'))=0$, with $g(x,x)=1$.

%We let
%\beq\label{5.13} V_y^\circ=\cV\cap \{y\} \times V_y,\and g_{y}(\cdot)=g(y,\cdot): V_y^\circ\lra \bar\cG.
%\eeq
Because the equation $\bR_{(w,w')}(g)$ is holomorphic in $(w,w')$ and $g$, $g(w,w')$ is holomorphic
in $(w,w')\in\cV^\circ$, and is the solution to
\beq\label{5.14} \dbar_x^*(\dbar_{w'}\cdot g(w,w')-\dbar_x)=0,\quad  (w,w') \in\cV^\circ.
\eeq
Without loss of generality, we can assume that $U_x$ is connected. Then since for any $y\in U_x$, $y\in U_y$,
we conclude that $g(y,y)=1$ for $y\in U_x$.

For $w\in V_x^\CC$, we let $V_w^\circ=\cV^\circ\cap \{w\}\times V_w$, and let 
$f_w^\circ=CS|_{V^\circ_w}$. We define 
$U_w^\circ=\Crit(f_w^\circ)\sub V_x^\circ$.
We define a family
$$\eta_w: V_w^\circ\lra \cA,\quad \eta_w(w')=\dbar_{w'}\cdot g(w,w').
$$
By our construction, $\eta_w(V^\circ_w)$ is a 
smooth complex submanifold of $\cA_{si}$ and $\eta_w$ is a biholomorphism onto its image. Further, because CS is invariant
under gauge transformations, we have $f_w\circ\eta_w\upmo =CS|_{V^\circ_w}$. %We denote $f'_y=CS|_{V'_y}$.

We now show that for any $w\in V_x^\CC$, $\eta_w(U^\circ_w)\sub U_x$ and 
$\eta_w|_{U^\circ_w}: U^\circ_w\to U_x$ is an open embedding, as analytic subspaces. 
Indeed, consider the image family  $\eta_w(V_w^\circ)=\{\dbar_{z}\,|\,z\in \eta_w(V_w^\circ)\}$, which is a family of semiconnections whose restriction to $\eta_w(U^\circ_w)$ is a family of integrable semiconnections. Adding \eqref{5.14},
we conclude that the family  $\{\bL_x(\dbar_{z}-\dbar_x)\,|\, z\in \eta_w(V_w^\circ)\}$ 
restricted to $\eta_w(U_w^\circ)$ vanishes. This proves that $\eta_w(U_w^\circ)\sub U_x$. 
On the other hand, by Lemma \ref{5.9j}, both $\eta_w(U_w^\circ)$ and $U_x$ become open neighborhoods of $\fV^c_{si}$ via the projection $\cA_{si}^{int}\to \fV^c_{si}$. Therefore, we conclude that $\eta_w(U_w^\circ)\sub U_x$ is an open analytic subspace.
%Because $g(y,y)=1$, $y\in U_y'\cap U_x$.

%Our next step is to modify $\eta_1$ to accommodate the requirement $f_x-f_y\circ\psi\in \fm_y^3$.
%To this end, we let $l=\dim T_xV_x/T_xU_x$; we find $l$ holomorphic functions $h_1,\cdots,h_l$ on
%$V_y^\circ$, vanishing along $U^\circ_y=U_y\cap V_y^\circ$, such that
%$dh_1(x),\cdots,dh_l(x)$ span the conormal space $\ker\{T\dual_yV_y\to T_y\dual U_y\}$.
%We then find $v_1,\cdots,v_l\in \Omega^{0,1}(\End E)_{k}$ so that,
%possibly after shrinking $V_y^\circ$ if necessary, the map 
%$$\eta_2: V^\circ_y\lra \cA,\quad \eta_2(y')= \dbar_{y'}\cdot g_y(y')+h_v(y'),
%$$
%where $h_v(y')=h_1(y')v_1+\cdots+h_l(y')v_l$,
%is biholomorphic onto its image, and further if we let $V_y\dpri=\eta_2(V_y^\circ)$, then
%%be its image with $\eta_2: V_y'\to V_y\dpri$ the induced biholomorphism, then 
%\beq\label{bb}
%T_y V_y\dpri=T_y V_x,\and f_y-CS\circ \eta_2\in \fm_y^3.
%\eeq
%Obviously, once $h_i$ is chosen, the desired $v_1,\cdots,v_l$ do exist.
We now construct the desired map $\zeta$. We first
find holomophic $\vartheta: \cV^\circ\to \OadE$ satisfying the system 
\beq\label{5.15}
\bL_{x}\bl \dbar_{w'}\cdot g(w,w')-\dbar_x+\vartheta(w,w')\br=0,
\quad \langle\vartheta(w,w'),\fb\rangle=0, \ \forall\, \fb\in \hat{\Xi}_x,
\eeq
for $(w,w')\in \cV^\circ$, and such that $\vartheta(y,y)=0$ for all $y\in U_x$ when defined.

First, for $y\in U_x$ since $g(y,y)=\id$, $\vartheta(y,y)=0$ is a solution. We let
$$\hat{\Xi}_x^\perp=\{\fa\mid 
\langle \fa,\fb\rangle=0, \forall \fb\in \hat{\Xi}_x\}\sub \Omega^{0,1}(\End E)_k.$$
For $z\in\cA$, we define 
\beq\label{5.16}
\bM_{z}(\cdot)=\bL_{x}(\dbar_{z}-\dbar_x+\cdot): \hat{\Xi}_x^\perp\lra  \Omega^{0,1}(\End E)_{k-2}/\ti\Xi_x.
\eeq
By our construction, the linearization $\delta\bM_{x}$ at $0$ is
\beq\label{5.17}
\delta\bM_{x}|_{0}=P_{x}\circ \Lap_{x}: \hat{\Xi}_x^\perp\lra \Omega^{0,1}(\End E)_{k-2}/\ti \Xi_x,
\eeq
which by Lemma \ref{4.35} is an isomorphism. Applying the implicit function theorem, 
we conclude that there is an open $x\in\cU\sub \cA$ so that for any $z\in\cU$, $\delta\bM_{z}|_{0}$
is an isomorphism. We assume that $V_x^\CC\sub \cU$. Thus after shrinking $(x,x)\in\cV^\circ$ if necessary,
we can solve
\beq\label{varth}
\bL_x(\dbar_{w'}\cdot g(w,w')-\dbar_x+\vartheta(w,w'))=0,
\eeq
subject to $\vartheta(y,y)=0$ for all $y\in U_x$, uniquely and smoothly in $(w,w')$. Since the operator \eqref{varth} is holomorphic in $(w,w')$, $\vartheta(w,w')$ is holomorphic in $(w,w')$.

%We let $U_y^\circ=U_y\cap V_y^\circ$.
We verify that for any $w\in V_x^\CC$, $\vartheta(w,\cdot)|_{U_w^\circ}=0$. %(Note $\eta_2|_{U_y'}=\id_{U_y'}$.) 
This is obvious because the family $\{\dbar_{w'}\cdot g(w,w')\,|\,w'\in V_w^\circ\}$ restricted to $U_w^\circ$
is a subfamily in $U_x$, thus $\vartheta(w,\cdot)|_{U_w^\circ}=0$ solves
the system \eqref{5.15} restricted to $U_w^\circ$. By the uniqueness, we conclude that $\vartheta(w,\cdot)|_{U_w^\circ}=0$.

We define $\zeta: \cV^\circ\to \cA$ by
\beq\label{5.18}
\zeta(w,w')=\dbar_{w'}\cdot g(w,w')+\vartheta(w,w') .
\eeq
Since $\vartheta$ satisfies the system \eqref{varth}, 
%$$\bM_{(w,w')}(\vartheta(w,w'))=\bL_x(\dbar_{y'}\cdot g(w,w')-\dbar_x+\vartheta(w,w'))=0,$$
$\zeta(w,w')\in V_x$. %=\dbar_{z}\cdot g_y(z)+\vartheta(w,w')\in V_x,\quad \forall\, z\in V_y\dpri.
Adding that $\zeta(x,\cdot): \cV^\circ\cap \{x\}\times V_x\to V_x$ is the identity embedding,
after shrinking $(x,x)\in \cV^\circ$ if necessary, we conclude that for any $w\in V_x^\CC$,
$\zeta(w,\cdot)$ is an open embedding of $V_w^\circ$ %=\cV^\circ\cap\{w\}\times V_w$
into $V_x$, and $\zeta(w,\cdot)|_{U_w^\circ}$ is an open embedding of $U_w^\circ$ into $U_x$. We remark that
by the construction, $\zeta(w,w')$ is holomorphic in $(w,w')$.
This proves the proposition.
\end{proof}

%\vfil\break
%\footnote{
%. Let 
%\beq\label{5.15a}U_y=\Crit(f_y)\sub V_y^\circ,\quad \text{for } f_y=CS|_{V_y^\circ}.\eeq
%
%We claim that 
%\beq\label{5.11} \bL_{x}\bl \dbar_{y'}\cdot g_y(y')-\dbar_x\br=0,\quad \forall y'\in U_y.\eeq 
%Indeed, since $\{\dbar_{y'}\,|\,y'\in U_y\}$ is a family of integrable semiconnections, so is 
%$\{\dbar_{y'}\cdot g_y(y')\,|\,y'\in U_y\}$. Combining \eqref{5.14} with the integrability 
%$F^{0,2}_{\dbar_{y'}\cdot g_y(y')}=0$, we find that \eqref{5.11} holds. Thus $\vartheta_y(y')=0$
%is a solution to \eqref{5.15} along $U_y=\Crit(f_y)\sub V_y^\circ$.
%Also $\vartheta_x(x')=0$ is a solution to \eqref{5.15} along $V_x$ by the definition of $V_x$. 
%}
%

We let $\zeta_w=\zeta(w,\cdot)$. Let $V_{x,w}^\circ$ be the image of $V_w^\circ$ under $\zeta_w$.
After shrinking $(x,x)\in\cV^\circ$ if necessary, we know that 
$\zeta_w:V_w^\circ\to V_{x,w}^\circ$ is a biholomorphic map. 
We let
$$\psi_w=\zeta_w^{-1}:V_{x,w}^\circ\lra V_w^\circ.
$$ % is the desired biholomorphic map.
We continue to denote $f_w^\circ=CS|_{V_w^\circ}$, $f_x=CS|_{V_x}$ and $f^\circ_{x,w}=f_x|_{V_{x,w}^\circ}$. We introduce
$$f_{w,t}=(1-t) f^\circ_{x,w}+t f_w^\circ\circ \psi_w,\quad t\in [0,1].
$$
These are holomorphic functions on $V_{x,w}^\circ$.

\begin{lemm}\label{dd2} There are open neighborhoods $x\in\cU\sub V_x^\CC$ and $x\in V_x^\circ\sub V_x$
so that for any $w\in\cU$, $V^\circ_{x,w}\sub V_x^\circ$, and the ideal  $\bl df_{w,t}\br\sub\sO_{V_{x,w}^\circ}$ is independent of $t\in [0,1]$.
\end{lemm}

\begin{proof}
We let $F$ be the holomorphic function on
$$\cW=\CC\times (\coprod_{w\in V_x^\CC}V_{x,w}^\circ)\sub \CC\times V_x^\CC\times V_x,
$$
defined via $F(t,w,w')=f_{w,t}(w')$. We let $d_{w'}$ be the differential along $V_x$ (the last factor), and form the ideal sheaf 
$$(d_{w'} F)\sub\sO_{\cW}.
$$

We let $I_{U_x}\sub\sO_{V_x}$ be the ideal sheaf of $U_x\sub V_x$; we let $\pr_3:\cW\sub \CC\times V_x^\CC\times V_x
\to V_x$ be induced by the third projection, and form the ideal sheaf $\pr_3\sta I_{U_x}$. By our construction
of $f_w$ and $\psi_w$, we see that $(d_{w'}F)\sub \pr_3\sta I_{U_x}$. The lemma will follow if we can find an
open subset $\cW^\circ\sub \cW$, containing $[0,1]\times\{ (x,x)\}$, such that 
\beq\label{equal}
(d_{w'}F)\otimes_{\sO_\cW}\sO_{\cW^\circ} =\pr_3\sta I_{U_x}\otimes_{\sO_\cW}\sO_{\cW^\circ}.
\eeq

We consider the quotient sheaf $\pr_3\sta I_{U_x}/(d_{w'}F)$, which is a finitely generated $\sO_\cW$-module
and is zero when tensored with $\sO_{\CC\times\{(x,x)\}}$. 
Thus by Nakayama Lemma and using that $[0,1]$ is compact, we conclude that there is an open neighborhood $\cW^\circ$
of $[0,1]\times \{(x,x)\}$ in $\cW$ making \eqref{equal} hold.
This proves the lemma.
\end{proof}

\subsection{Local equivalence}\label{S5.1.3} 
The biholomorphic map $\psi_y:V_{x,y}^\circ\to V_y^\circ$ may not pull back $f_y^\circ=CS|_{V_y^\circ}$ to $f_x^\circ=CS|_{V_{x,y}^\circ}$. 
In this section, we use the vector field technique to complete our proof of Proposition \ref{4.42}.
To keep the notation simple, we denote $V_{x,y}^\circ$ by $V_x^\circ$ and let $\psi=\psi_y$. 

%We let the notation be as in the previous subsection: $V_x=V(x,h,\Xi_x)$ is the CS chart from Theorem \ref{4.28}; $f_x=CS|_{V_x}$ is the restriction of the Chern-Simons functional; $U_x=\Crit(f_x)\sub V_x$ so that $\Pii(U_x)$ is an open neighborhood of $\bar{x}=\Pii(x)\in \fV_{si}^c$ by the quotient map $\Pii:\cA_{si}\to \cB_{si}=\cA_{si}/\cG_k$. For $y\in U_x$, $V_y=V(y,h,\Xi_y)$ and $f_y=CS|_{V_y}$. We have a biholomorphic map $\psi:V_x^\circ\to V_y^\circ$ of open neighborhoods of $y$, from Proposition \ref{5.10}. For $y\in U_x$, $y'$ denotes a point in $V_y$ and let $f_y^\circ=f_y|_{V_y^\circ}$; $U_y^\circ=\Crit(f_y^\circ)\sub V_y^\circ$. As always, for psychological reasons, we write $\dbar_x$ for $x\in \cA$, i.e. $x=\dbar_x$. 
%\medskip

We first rephrase the construction of the map $\psi$. 
Since $\zeta(y')=\dbar_{y'}\cdot g_y(y')+\vartheta_y(y')$ by \eqref{5.18}, we have
\beq\label{5.20}
y'=\dbar_{y'}=(\zeta(y')-\vartheta_y(y'))\cdot g_y(y')^{-1}.\eeq
We let $$\tilde{g}_y(x')=g_y(\psi(x'))^{-1}\and \tilde{\vartheta}_y(x')=-\vartheta_y(\psi(x')).$$
Then we have
 $\tilde{\vartheta}_y|_{U_x\cap V_x^\circ}=0$ and 
\beq\label{5.19} \psi(x')=(x'+\tilde{\vartheta}_y(x'))\cdot \tilde{g}_y(x') \quad\text{ for all }x'\in V_x^\circ.\eeq

%We let $f_x^\circ=f_x|_{V_x^\circ}$ so that $U_x^\circ=U_x\cap V_x^\circ=\Crit(f_x^\circ)$. 
Let
$$\cI=(df_x^\circ)\sub\sO_{V_x^\circ}$$ 
be the ideal sheaf of $U_x^\circ\sub V_x^\circ$.

\begin{lemm}\label{5.21} 
Let $\psi$ be defined by \eqref{5.19}. Then ${f_y^\circ}\circ\psi-f_x^\circ\in \cI^2$.
\end{lemm}

\begin{proof}
Without loss of generality, we assume $U_x$ is connected so that the Chern-Simons functional $CS$ is independent of the reference point. Since $CS$ is invariant under infinitesimal gauge transformations, we have
$$f_y^\circ(\psi(x'))=CS(\psi(x'))=CS((x'+\tilde{\vartheta}_y(x'))\cdot \tilde{g}_y(x'))=CS(x'+\tilde{\vartheta}_y(x')).
$$
Since $f_x^\circ(x')=CS(x')$, it suffices to show that 
\beq\label{5.22}
CS(x')-CS(x'+\tilde{\vartheta}_y(x'))\in \cI^2.
\eeq

We use finite dimensional approximation to reduce this to a familiar problem in several complex variables. 
First, since $\tilde{\vartheta}_y$ takes values in $C^\infty$-forms in $\OadE$, we can find an integer $l$
so that $\tilde{\vartheta}_y$ factors through
$$\tilde{\vartheta}_y: V_x^\circ\lra \Omega^{0,1}(\End E)_{L^2_l}\sub\OadE.
$$
Since $\Omega^{0,1}(\End E)_{L^2_l}$ is a separable Hilbert space, we can approximate it by an increasing sequence of 
finite dimensional subspaces $R_m\sub \Omega^{0,1}(\End E)_{L^2_l}$. Let 
$$q_m: \Omega^{0,1}(\End E)_{L^2_l}\to R_m$$
be the orthogonal projections. Then we have 
a convergence of holomorphic functions 
$$\lim_{m\to\infty} CS(x'+ q_m\circ \tilde{\vartheta}_y(x'))=CS(x'+\tilde{\vartheta}_y(x')), 
$$
uniformly on every compact subset of $V_x^\circ$.

We claim that
\beq\label{5.23}
CS(x')-CS(x'+ q_m\circ \tilde{\vartheta}_y(x'))\in \cI^2.
\eeq
Note that the claim and the uniform convergence imply \eqref{5.22}.

We pick a basis $e_1,\cdots, e_{n}$ of $R_m$.
As $q_m\circ\tilde{\vartheta}_y: V_x^\circ\to R_m$ is holomorphic, we can find holomorphic functions
$w_i: V_x^\circ\to\CC$ so that  
$$q_m\circ\tilde{\vartheta}_y(x')=\sum_{i=1}^n w_i(x')e_i.
$$
Because $\tilde{\vartheta}_y|_{U_x^\circ}=0$, we have $w_i(x')\in \cI$ for all $i$.
Because the Chern-Simons functional $CS$ is a cubic polynomial in $\fa\in \OadE$, using the Taylor expansion, we have
$$CS(x'+ q_m\circ \tilde{\vartheta}_y(x'))-CS(x')\equiv \sum_{i=1}^n\frac{\partial CS}{\partial e_i}(x')\cdot w_i(x')\mod \cI^2.
$$
Because $U_x^\circ\sub V_x^\circ$ is defined by the vanishing of the partial derivatives of $CS$,
we have $\frac{\partial CS}{\partial e_i}(x')\in\cI$. Adding that $w_i(x')\in\cI$, we obtain \eqref{5.23}.
This proves the lemma.
\end{proof}

Now we can complete our proof of Proposition \ref{4.42}.

\begin{proof}[Proof of Proposition \ref{4.42}]
Let $f_0=f_x^\circ$ and $f_1={f_y^\circ}\circ\psi$. 
To prove Proposition \ref{4.42}, it suffices to find a
biholomorphic map $$\varphi:V_x^\circ\lra V_x^\circ$$ such that 
\beq\label{5.24}f_1\circ\varphi=f_0\eeq
because then
$\Psi=\psi\circ\varphi:V_x^\circ\to V_y^\circ$
satisfies $${f_y^\circ}\circ\Psi={f_y^\circ}\circ\psi\circ\varphi=f_1\circ\varphi=f_0=f_x^\circ$$
and so $\Psi$ is the desired local equivalence of the CS charts at $y$.

Let
$$f_t=(1-t)f_0+t f_1,\quad t\in [0,1].$$
By Lemma \ref{dd2}, we know that the ideal $(df_t)$ generated by the partial derivatives of $f_t$ is independent of $t\in [0,1]$. 

%Indeed, since $f_0-f_1\in \cI^2$ by Lemma \ref{5.21} and $\cI=(df_0)$, we find that $\cI\supset (df_1)$. By applying Lemma \ref{5.21} again after interchanging the roles of $x$ and $y$, we find that
%$$f_x^\circ \circ\psi^{-1}-f_y^\circ\in (df_y^\circ)^2$$
%and that $\cI=(df_0)\subset (df_1)$. Therefore we have
%$$\cI=(df_0)=(df_1)$$
%and hence $$(df_\gamma)\subset \cI, \quad \forall\gamma\in\CC.$$
%This means that the analytic space $$U_x^\circ\times \CC=\mathrm{zero}(\cI)\times\CC$$ defined by the ideal $\cI$ is contained in the analytic space
%$$\mathrm{zero}(df_\gamma)\sub V_x^\circ\times \CC$$
%defined by $(df_\gamma)$ for $\gamma\in \CC$. Since a flat deformation is a maximal deformation (cf. Lemma \ref{5.80}), we find that $(df_\gamma)=\cI$ holds for $\gamma\in \CC$
%which does not belong to a discrete subset of $\CC-\{0,1\}$. Hence we can find a path $\gamma$ from $0$ to $1$ so that we have the identity of ideals
%$$(df_{\gamma(t)})=\cI,\quad \forall t\in [0,1].$$

By Lemma \ref{5.21}, we can find a family of holomorphic vector field $\xi_t$ on $V_x^\circ$, smooth in $t\in[0,1]$, satisfying
\beq\label{5.25}
df_{t}(\xi_t)=f_0-f_1 %\cdot \frac{d\gamma}{dt}
\eeq
which vanishes along $U_x^\circ=\Crit(f_0)$.
If we let $\varphi_t(x')$ denote the integral curve of the vector field $\xi_t$ with the initial condition $\varphi_0(x')=x'$, we find that 
$$f_{t}(\varphi_t(x'))$$
is constant. Moreover $\varphi_t$ fixes points in $U_x^\circ$. Since the vector field $\xi_t$ is holomorphic, $\varphi_t$ is a biholomorphic map for all $t\in [0,1]$ in a neighborhood of $y$. 

Let $\varphi=\varphi_1$. Then 
$$f_1(\varphi(x'))=f_0(x'),\quad \forall x'$$
in a neighborhood of $y$ in $V_x^\circ$. This proves the proposition. 
\end{proof}
%
%{\blue
%\begin{lemm}\label{5.80}
%Let $V$ be a smooth complex manifold and $S$ an analytic space. Suppose $W_1\sub W_2$ are two
%closed analytic subspaces of $V\times S$ such that the induced projection $W_1\sub V\times S\to S$ is flat, and that  there is a closed 
%analytic subspace 
%$S_0\sub S$ such that $W_1\times_S S_0=W_2\times_S S_0$. Then there is an open subset $\cU_0\sub V\times S$ that contains
%$V\times S_0$ such that $W_1\cap \cU_0= W_2\cap \cU_0$ as analytic subspaces in $\cU_0$.
%\end{lemm}
%
%\begin{proof}
%Let $\cA$ be the kernel of the surjection $\cO_{W_2}\to \cO_{W_1}$. Since $\cO_{W_1}$ is flat over $\cO_S$, we have an exact sequence
%\[ 0\lra  \cA\otimes_{\cO_S}\cO_{S_0} \lra \cO_{W_2}\otimes_{\cO_S}\cO_{S_0} \mapright{\phi} \cO_{W_1}\otimes_{\cO_S}\cO_{S_0} \lra 0.
%\]
%Since $W_1\times_S S_0=W_2\times_S S_0$, $\phi$ is an isomorphism, thus $ \cA\otimes_{\cO_S}\cO_{S_0}=0$.
%As $\cA$ is a coherent sheaf of $\cO_{V\times S}$-modules, there is an open $\cU_0\sub V\times S$, containing $V\times S_0$,
%such that $\cA|_{\cU_0}=0$. This proves that $W_1\cap \cU_0=W_2\cap \cU_0$, as analytic subspaces of $\cU_0$.
%\end{proof}
%
%{\red [The Lemma can not be applied to the previous case. Though $\cU_0$ is open, it is merely analytic open,
%not Zariski open. This lemma can not be strengthened.]}
%}
%

\bigskip

\section{Proof of Proposition \ref{4.41}}\label{S5.2}

In this section, by applying parallel arguments to those in \S\ref{S5.1}, 
we will prove Proposition \ref{4.41} and complete our proof of Theorem \ref{4.45}.

We have to show that the CS chart $V(x,h,\Xi)$ is independent of the choice of $x$ in its orbit $x\cdot\cG_k$, the hermitian metric $h$ on $E$, and the CS framing $\Xi$ at $x$. 

We first observe that the choice of $x$ is irrelevant.
\begin{lemm}\label{5.26}
Let $x,x'\in \cA^{int}_{si}$ with $x'=x\cdot g$ for $g\in \cG_k$. Let $h'=h^g=g^{\dagger}hg$ be the hermitian metric obtained from $h$ by applying the action of $g$. Let $\Xi'=\Xi^g$ be the image of $\Xi$ by the tangent map $T_g:T_x\cA\to T_x\cA$ of $g$. Let $V=V(x,h,\Xi)$ and $V'=V(x',h',\Xi')$ be the CS charts defined by Theorem \ref{4.28}. Let $f=CS|_V$ and $f'=CS|_{V'}$. Then the LG pair $(V,f)$ is equivalent to $(V',f')$.  
\end{lemm}
See Definition \ref{1.4} (2) for the definition of an equivalence of LG pairs. 
\begin{proof}
By definition, $g$ sends the triple $(x,h,\Xi)$ to $(x',h',\Xi')$ and hence it sends $V$ to $V'$. Moreover the Chern-Simons functional $CS$ is gauge invariant so that $CS(\dbar\cdot g)=CS(\dbar)$. 
So we get the equivalence $g:V\to V'$. 
\end{proof}

By Lemma \ref{5.26}, we may assume $x=x'$ and prove that $V(x,h,\Xi)$ is independent of $h$ and $\Xi$ in a neighborhood of $x$. 

Let $h,h'$ be two hermitian metrics on $E$. Then
\beq\label{5.27} h_t=(1-t)h+th',\quad t\in [0,1]\eeq
is a real analytic family of hermitian metrics on $E$. Let $\Xi,\Xi'$ be two
CS framings at $x$. We can also find a (piecewise) real analytic family of CS framings at $x$ connecting $\Xi$ and $\Xi'$.
\begin{lemm}\label{5.28} 
There is a real analytic family $\Xi_t$, $t\in [1,2]$ of CS framings at $x$ such that $\Xi_1=\Xi$ and $\Xi_2=\Xi'$ with respect to the fixed hermitian metric $h_1=h'$. 
\end{lemm}
\begin{proof}
By Lemma \ref{5.5}, an $r$-dimensional subspace 
$$W\sub \ker (\dbar_x^*)^{0,1}_k$$ containing $\Lap_x^{-1}(0)^{0,1}$ is a CS framing if and only if 
the pairing \beq\label{5.29} W\times \hat{W}\lra \CC\eeq
defined by \eqref{4.22} is nondegenerate. Let $A$ be a finite dimensional subspace of 
$\ker (\dbar_x^*)^{0,1}_k$ containing $\Xi\cup \Xi'$. Let 
$$\bar{A}=A/ \Lap_x^{-1}(0)^{0,1}$$
and $l=r-\dim \Lap_x^{-1}(0)^{0,1}$. 
Then the Grassmannian $Gr(l,\bar A)$ parameterizes $r$-dimensional subspaces $W$ of $A$ containing 
$\Lap_x^{-1}(0)^{0,1}$. The locus $Gr(l,\bar A)^\circ$ of $r$-dimensional subspaces $W$ containing $\Lap_x^{-1}(0)^{0,1}$ for which \eqref{5.29} is nondegenerate is the complement of a divisor in $Gr(l,\bar A)$. In particular $Gr(l,\bar A)^\circ$ is a smooth connected variety and hence we can connect two points $\Xi$ and $\Xi'$ by a piecewise real analytic path $\Xi_t$, $t\in [1,2]$. 
\end{proof}

We define a family
$$(x,h_t,\Xi_t)=\left\{\begin{matrix} (x,h_t,\Xi) & t\in [0,1]\\
(x,h_1,\Xi_t) & t\in [1,2]\end{matrix}\right.
$$ of triples which is piecewise real analytic. 
Let $$V_t=V(x,h_t,\Xi_t),\quad t\in [0,2]$$
be the CS charts from Theorem \ref{4.28}. Let $f_t=CS|_{V_t}$. 
Since the interval $[0,2]$ is compact, Proposition \ref{4.41} follows from the following.
\begin{prop}\label{5.30}
For $t\in [0,2]$, there exists $\delta>0$ such that for any $s\in (t-\delta,t+\delta)$, the LG pair $(V_s,f_s)$ is locally equivalent to $(V_t,f_t)$ at $x$.
\end{prop}

The proof of Proposition \ref{5.30} parallels that of Proposition \ref{4.42} in \S\ref{S5.1}. 

\medskip

Firstly, the proof of Proposition \ref{5.10} gives us the following.
\begin{lemm}\label{5.31}
Let $U_t=\Crit(f_t)$. Then there exists $\delta>0$ such that for any $s\in (t-\delta,t+\delta)$, there exists a biholomorphic map $\psi:V_t^\circ\to V_s^\circ$ of open neighborhoods $V_t^\circ\sub V_t$ and $V_s^\circ\sub V_s$ of $x$ such that
$$\Pii|_{V_t^\circ\cap U_t}=\Pii\circ\psi|_{V_t^\circ\cap U_t}:V_t^\circ\cap U_t\lra \fV_{si}^c$$
where $\Pii:\cA_{si}\to \cB_{si}$ is the quotient map by the action of $\cG_k$. 
\end{lemm}

Secondly, the proof of Lemma \ref{5.21} gives us the following.
\begin{lemm}\label{5.32} Let $f_t^\circ=f_t|_{V_t^\circ}$ and $f_s^\circ=f_s|_{V_s^\circ}$. 
For $\delta>0$ in Lemma \ref{5.31} and for any $s\in (t-\delta,t+\delta)$, 
$$f_s^\circ\circ\psi-f_t^\circ\in \cI^2$$
where $\cI=(df_t^\circ)\sub \sO_{V_t^\circ}$ is the ideal generated by the partial derivatives of $f_t^\circ$. 
\end{lemm}

Finally, the proof of Proposition \ref{4.42} in \S\ref{S5.1.3} gives us a proof of Proposition \ref{5.30}.

\medskip

This completes our proof of Theorem \ref{4.45}.

\chapter{Moduli of sheaves and orientability}\label{ch6}

Recall that $\fV_{si}^c$ and $\SH^c$ denote the stack of simple vector bundles and simple sheaves respectively of total Chern class $c$, on the Calabi-Yau 3-fold $Y$.
In Chapters \ref{ch4} and \ref{ch5}, we proved that an (open) moduli space $X\sub \fV_{si}^c$ of simple vector bundles 
admits a critical virtual manifold structure if $d_x=\dim T_xX$ is bounded from above (cf. Theorem \ref{4.45}). 
In this chapter, we deal with two remaining issues.

\medskip

\noindent (1) We extend the existence of a critical virtual manifold structure to moduli spaces of simple sheaves on $Y$. For a bounded moduli space $X\sub \SH^c$ of simple sheaves, we use the Seidel-Thomas twists to construct an open embedding of $X$ into $\fV_{si}^{\ti c}$ for some total Chern class $\ti c$. We therefore obtain a critical virtual manifold structure on $X$. 

\medskip 

\noindent (2) We prove that the critical virtual manifold structure on an (open) moduli space $X\sub \SH^c$ of simple sheaves is orientable if the reduced space $X\ured$ admits a tautological family, i.e. a sheaf $\sE$ on $X\ured\times Y$, flat over $X\ured$, whose restriction to $\{x\}\times Y$ is the sheaf represented by $x$. 

\bigskip

\section{Critical virtual manifold structure on moduli of simple sheaves}\label{S6.1}

In this section, we use the Seidel-Thomas twists and show that an open analytic subspace $X\sub\SH^c$ admits a critical virtual manifold structure.

\medskip

We begin with an obvious corollary from Theorem \ref{4.28}.
\begin{coro}\label{6.1}
Let $X$ be an analytic space which has an open embedding into the moduli stack $\fV_{si}^{\ti c}$ of simple locally free sheaves on $Y$ for some total Chern class $\ti c$ such that $\dim T_x\fX$ is a bounded from above. Then $X$ has a critical virtual manifold structure.  
\end{coro}

We fix an ample line bundle $\sO_Y(1)$ on $Y$, and use $\sE(m)$ to
denote the twisting $\sE\otimes_{\sO_Y}\sO_Y(m)$. 
For a positive integer $m$,
let $\frak{ST}[m]\sub\SH^c$ be the open substack of
$[\sE]\in\SH^c$ such that \begin{enumerate}
\item $\sE(m)$ is spanned by global sections and 
\item $h^i(\sE(m))=0$ for $i>0$.
\end{enumerate}
We use the $m$-th Seidel-Thomas twist (\cite[\S2]{SeiTho})  defined by 
\beq\label{6.2}
\cT_m: \frak{ST}[m]\lra \SH^{c'},\quad \sE\mapsto \ker\{H^0(\sE(m))\otimes \sO_Y(-m)\to \sE\}=\sE'
\eeq
whose inverse functor is 
\beq\label{6.3}
\sE'\mapsto \coker\{\sE'\to \Hom(\sE',\sO_Y(-m))\otimes \sO_Y(-m)\}.
\eeq

\begin{prop}\label{6.4} \cite[\S2]{SeiTho}
The morphism $\cT_m$ is an open embedding of stacks. If $X\subset \frak{ST}[m]$ admits a tautological family, so does $\cT_m(X)$.
\end{prop}
\def\fZ{\mathfrak{Z}}
\def\by{\mathbf{y}}

\begin{defi}\label{6.5}  An open analytic space $X\sub \SH^c$ is called \emph{bounded} if 
there is an integer $m$ such that $X\sub \frak{ST}[m]$.
\end{defi}

\begin{prop}\label{6.7}
Let $\fX\sub \SH^c$ be a bounded open analytic space. Then there is a total Chern class $\ti c$ such that $X$ is isomorphic to an open analytic subspace $\tilde{X}$ of the moduli space $\fV_{si}^{\ti c}$ of simple vector bundles on $Y$ with total Chern class $\ti c$. Moreover, $X\ured$ admits a tautological family on $X\ured\times Y$ if and only if $\tilde{X}\ured$ does on $\tilde{X}\ured\times Y$.
\end{prop}

\begin{proof}
Because $X$ is bounded, we can find an $m$ so that $X\sub \frak{ST}[m]$.
Let $Quot_Y^{M, P(n)}$ be Grothendieck's quot scheme parameterizing quotient sheaves $\sO_Y^{\oplus M}\to \sF$
of Hilbert polynomial $P(n)$, where $M=h^0(\sE(m))$ and $P(n)=\chi(\sE(m+n))$ for an $\sE\in X$.
Since $X\sub \frak{ST}[m]$, the map 
$$X\lra Quot_Y^{M, P(n)}, \quad \sE\mapsto [H^0(\sE(m))\otimes\sO_Y\to \sE(m)]$$ 
is an embedding. As $Quot_Y^{M, P(n)}$ is projective, we can find 
integers $m_0$, $m_1$ and $m_2$ such that 
$$\fX_0\defeq\fX\sub \frak{ST}[m_0],\quad X_{i+1}:=\cT_{m_i}(X_i)\sub\frak{ST}[m_{i+1}],
\and X_3:=\cT_{m_2}(X_2).$$ 
By Proposition \ref{6.4}, we have open embeddings $$X_1\sub \SH^{c'},\quad X_2\sub \SH^{c''},\quad X_3\sub\SH^{c'''}$$ and isomorphisms $$X=X_0\cong X_1\cong X_2\cong X_3.$$
Since $\dim Y=3$, $X_3$ is an open analytic subspace of $\SH^{\ti c}$ for $\ti c:=c'''$ that consists entirely of locally free sheaves. Therefore $X\cong X_3\sub \fV^{\ti c}_{si}$ is an open analytic subspace of $\fV^{\ti c}_{si}$. 

The last statement is obvious from the definition of Seidel-Thomas twists \eqref{6.2} and \eqref{6.3}.
\end{proof}
Combining Corollary \ref{6.1} and Proposition \ref{6.7}, we obtain the following.
\begin{theo}\label{6.8}
Let $\fX\sub \SH^c$ be a bounded open analytic space. Then $\fX$ admits a critical virtual manifold structure.
\end{theo}

\begin{proof} 
We need to verify that $\dim Ext^1(\sE,\sE)$ is uniformly bounded for $\sE\in X$. But this 
follows from the embedding $X\to Quot_Y^{M, P(n)}$ constructed in the proof of Proposition \ref{6.7} and that $Quot_Y^{M, P(n)}$ is projective.
\end{proof}

In the algebraic case, the boundedness assumption is always satisfied if $X$ is of finite type. 
\begin{coro}\label{6.9}
An algebraic space $X\sub \SH^c$ of finite type admits a critical virtual manifold structure. 
\end{coro}

\bigskip

\def\tcE{\tilde{\cE}}
\def\cAi{{\cA^{int}_{si}}}
\def\tcL{{\tilde{\cL}}}
\def\tz{\tilde{z}}
\def\tg{\tilde{g}}
\def\tB{\tilde B}
For the orientation issue, we first give several equivalent constructions of the determinant line bundle on $X$. 
We let $\tilde{\cE}$ be the tautological family of holomorphic bundles over $\cA^{int}_{si}$.
% the space of simple integrable semiconnections. 
Namely, let $\pr_Y: Y\times\cA_{si}^{int}\to Y$ be the projection, then $\ti\cE$ is the complex vector bundle $\pr_Y\sta E$
coupled with the tautological family of semiconnections $\dbar_{\cA^{int}_{si}}$ characterized by
that for any $x\in \cA_{si}^{int}$, $\dbar_{\cA^{int}_{si}}|_{Y\times\{x\}}=\dbar_x$. 

Let $\bar\cG=\cG/\CC\sta$ be the
reduced gauge group so that $\bar\cG$ acts on $\cA_{si}^{int}$ freely, and this action canonically lifts
to an action on $\ti\cE$ making the pair $(\cA_{si}^{int}, \ti\cE)$ $\bar\cG$-equivariant.

Let $\pr_\cA: Y\times\cA_{si}^{int}\to\cA^{int}_{si}$ be the second projection. 
We define
\beq\label{6.91} \tcL=\det R\pr_{\cA\ast}R\cH om(\tcE,\tcE).
\eeq
% where $p:\cAi\times Y\to \cAi$ is the projection. 
The $\bar\cG$ action on $\ti\cE$ induces an action on $\ti\cL$. Since $\fV_{si}^c=\cA_{si}^{int}/\bar\cG$ is the
quotient by a free action, the $\bar \cG$-equivariant line bundle $\ti\cL$ descends to a line bundle $\cL$ on
$\fV^c_{si}$, such that $\ti\cL=\pr_\cA\sta\cL$, where the action on $\ti\cL$ is the obvious induced 
$\bar\cG$-action on $\pr_\cA\sta\cL$.
For a subspace $X\sub \cA_{si}^{int}$, we define $\cL_X=\cL|_X$.

%
%we consider the determinant line bundle
%\beq\label{6.11}\det Rp_*R\cH om(\cE,\cE)\eeq
%where $p:X\times Y\to X$ is the projection and $\cE$ is the tautological family of simple sheaves on the Calabi-Yau 3-fold $Y$ parameterized by 

The determinant line bundle is well-defined when $X\sub \SH^c$. 
Let $X=\cup_\alpha X\lalp$ be an open cover so that over each $X\lalp$ we have a tautological family 
$\cE_\alpha$ over $X\lalp$. We form the determinant line bundles 
%(abbreviated to $Det(\cE\lalp,\cE\lalp)$)
\beq\label{det-glue}
%\Det(\cE\lalp,\cE\lalp)\defeq 
\det R{p\lalp}_*R\cH om(\cE\lalp,\cE\lalp),
\eeq
where $p\lalp:X\lalp\times Y\to X\lalp$ is the projection. 
Using that $Y$ is a Calabi-Yau 3-fold, an easy argument shows that
% $\sum (-1)^i\Ext^i(\cE\lalp|_{x}, \cE\lalp|_{x})=0$, for any $x\in X\lalp$,
the collection \eqref{det-glue} glues to a line bundle $\cL_X$ on $X$,
which is independent of the choice of the tautological families $\cE\lalp$.

It is clear that in case $X\sub\fV^c_{si}$, the two constructions yield the identical line bundle.
We call $\cL_X$ the \emph{determinant line bundle} of $X$.

\bigskip

\section{Orientability for moduli of sheaves}\label{S6.2}

In this section, we deal with the orientation issue.  By Proposition \ref{6.7}, a bounded open $X\sub \SH^c$ is isomorphic to an open $\tilde{X}\sub \fV_{si}^{\tilde{c}}$. Moreover $X\ured$ admits a tautological family if and only if $\tilde{X}\ured$ does. We prove that the critical virtual manifold structure on $X$ induced from that on $\tilde{X}$ (cf. Theorem \ref{6.8}) is orientable if $X\ured$ admits a tautological family (cf. Theorem \ref{6.41}). By Theorems \ref{2.23} and \ref{2.32}, we then have a perverse sheaf $P$ and a mixed Hodge module $\cM$ on $X$ which are locally the perverse sheaf and mixed Hodge module of vanishing cycles.

In \S\ref{S6.2.1}, we prove that the critical virtual manifold structure on $\tilde{X}\sub \fV_{si}^{\tilde{c}}$ is orientable if the determinant line bundle $\cL_{\tilde{X}\ured}$ admits a square root (cf. Theorem \ref{6.20}). In \S\ref{S6.2.2}, we prove that $\cL_{\tilde{X}\ured}$ admits a square root if there is a tautological family on $\tilde{X}\ured$ (cf. Proposition \ref{6.40}), which is equivalent to saying that there is a tautological family on $X\ured$. The main result, Theorem \ref{6.41}, follows from these two results.

%We prove that the critical virtual manifold structure on an open $X\sub\SH^c$ constructed in Theorem \ref{6.8} is orientable if and only if the determinant line bundle $\cL_{X\ured}$ on the reduced space $X\ured$ admits a square root . Then we prove that $\cL_{X\ured}$ admits a square root if there is a tautological family $\cE$ on $X\ured\times Y$ for the reduced space $X\ured$ underlying $X$ . We therefore obtain an orientable critical virtual manifold structure on $X$ when there is a tautological family on $X\ured\times Y$. 
\bigskip

\subsection{A criterion for orientability}\label{S6.2.1}
The purpose of this subsection is to prove the following criterion for orientability. 

\begin{theo}\label{6.20}
Let $X\sub\fV^c_{si}$ be a bounded open analytic space. 
%Suppose that there is a tautological family $\cE$ of simple sheaves parameterized by $X$.
Then the critical virtual manifold structure on $X$ by Theorem \ref{4.45} 
is orientable if the determinant line bundle 
$\cL_{X\ured}=\cL_X|_{X\ured}$ on $X\ured$ (before \eqref{det-glue}) admits a square root. 
\end{theo}

The rest of this subsection is devoted to a proof of Theorem \ref{6.20}.

\medskip

\def\tU{\tilde{U}}
 
%By the Seidel-Thomas twists (Proposition \ref{6.7} and Lemma \ref{6.10}), it suffices to consider the case of vector bundles. 
We will use the notation in the proof of Proposition \ref{4.44} where 
$$
\{U_x=B(x,\varepsilon_x)\,|\,x\in X\}
$$
form an open cover of $X$ and 
\beq\label{a3}U_x=\Crit(f_x)\sub V_x\mapright{f_x}\CC.
\eeq
Moreover, we have an equivalence 
\beq\label{6.93} \varphi_{xx'}:V_{xx'}\lra V_{x'x}
\eeq
of open neighborhoods $V_{xx'}\sub V_x$ and $V_{x'x}\sub V_{x'}$ of $U_{xx'}=U_x\cap U_{x'}$. 

By Proposition \ref{1.12}, the transition functions 
\beq\label{6.92}
\xi_{xx'}=\det\, d\varphi_{xx'}:K^\vee_x|_{U_{xx'}\ured}\to K^\vee_{x'}|_{U_{xx'}\ured}
\eeq 
of the line bundles $K^\vee_x:=\det T_{V_x}|_{U_x\ured}$ satisfy the rigidity 
\beq\label{6.21} \xi_{xx'x''}=\xi_{x''x}\circ\xi_{x'x''}\circ\xi_{xx'} =\pm 1.\eeq 
The critical virtual manifold $(X=\cup U_x, V_x,f_x,\varphi_{xx'})$ is orientable if the cohomology class 
\beq\label{b1}  [\{\xi_{xx'x''}\}]\in H^2(X,\ZZ_2)
\eeq
defined by the 2-cocycle $\{\xi_{xx'x''}\in \ZZ_2=\{\pm 1\}\}$
is zero (cf. Definition \ref{1.17}). 

To prove the theorem, we will show that the class \eqref{b1} is the obstruction to the existence of a square root of
$\cL_{X\ured}$. Let $U=\Crit(f)\sub V\mapright{f}\CC$ be one of the charts in \eqref{a3}.
Following our construction, $U\sub V\sub \cA$ and $f=CS|_V$. We first construct a canonical isomorphism
\beq\label{b2}
(\det T_V|_{U})^{2}\cong \ti\cL^{-1}|_{U},
\eeq
where $\ti\cL$ is defined by \eqref{6.91}. Indeed, 
let $\fA^{0,i}(\End E)$ denote the sheaf of differentiable (i.e. $C^\infty$) $(0,i)$-forms on $Y$ with values in $\End E$.
% so that $\Gamma(Y,\fA^{0,i}(\End E))=\Omega^{0,i}(\End E)$. 
Let $\fA^{0,i}(\End E)_{U}$ denote the pullback of $\fA^{0,i}(\End E)$ by the projection $U\times Y\to Y$. 
Let 
$$\tilde A_U\ud=[\ \cdots \lra \fA^{0,i}(\End E)_U\mapright{\dbar_U} \fA^{0,i+1}(\End E)_U\lra \cdots\
%\mapright{\dbar} \fA^{0,2}(\End E)_U\mapright{\dbar} \fA^{0,3}(\End E)_U
]$$
denote the tautological complex where the differentials $\partial_U$ at $x\in U$ is $\dbar_x$. 
Then $\tilde A_U\ud$ is a flabby resolution of $\cH om(\tcE,\tcE)|_U$. Using 
$\Gamma(Y,\fA^{0,i}(\End E))=\Omega^{0,i}(\End E)$ and denoting
$\Omega^{0,i}(\End E)_U=U\times \Omega^{0,i}(\End E)$, we conclude that
$$A_U\ud:=p_{\ast}\tilde A_U\ud = [\ \cdots\lra \Omega^{0,i}(\End E)_U\mapright{\dbar_U} \Omega^{0,i+1}(\End E)_U\lra\cdots\,]
$$%\mapright{\dbar} \Omega^{0,2}(\End E)_U\mapright{\dbar} \Omega^{0,3}(\End E)_U]$$
is quasi-isomorphic to 
$Rp_*R\cH om(\tcE,\tcE)|_U$ where $p:X\times Y\to X$ is the projection. Hence $\tilde{\cL}|_U=\det Rp_*R\cH om(\tcE,\tcE)|_U$ is isomorphic to $\det A_U\ud$ and the isomorphism is compatible with the gauge action.

Because $V\subset \cA_{si}$, we have an induced
embedding $T_{V}|_{U}\subset \Omega^{0,1}(\End E)_U$. Following the proof of Theorem \ref{4.28}, 
we have constructed a subbundle $R|_U \sub \Omega^{0,2}(\End E)_U$ (cf. \eqref{4.30}) 
and an isomorphism of $R|_U \cong T^*_{V}|_{U}$ that fits into a commutative diagram %quasi-isomorphism
$$%\beq\label{6.25}
\xymatrix{
& [T_{V}|_{U}\ar[r]^{d(df)|_U}\ar[d] & T^*_{V}|_{U}]\ar[d]&\\
[\Omega^{0,0}(\End E)_{U}\ar[r]^\dbar & \Omega^{0,1}(\End E)_{U}\ar[r]^\dbar & \Omega^{0,2}(\End E)_{U}\ar[r]^\dbar &
\Omega^{0,3}(\End E)_{U}].
}$$%\eeq
(Note $f=CS|_V$.) This diagram induces a distinguished triangle 
\beq\label{6.25} [T_{V}|_{U}\to T\sta_{V}|_{U}]\lra A^\bullet_U\lra \sO_U\oplus \sO_U[-3]\mapright{+1},
\eeq
which induces the isomorphism \eqref{b2}, as desired.

We now consider two charts 
$$U_x=\Crit(f_x)\sub V_x\mapright{f_x}\CC \and U_{x'}=\Crit(f_{x'})\sub V_{x'}\mapright{f_{x'}}\CC,$$ as in \eqref{a3}, such that $\pi(U_x)\cap \pi(U_{x'})\ne \emptyset$. (Here we deviate from the convention
\eqref{a3} in that we view $U_x\sub V_x\sub\cA$, and use $\pi(U_x)$ to denote the open subspace in $X\sub \fV^c_{si}$,
where as usual $\pi: \cA_{si}^{int}\to\fV^c_{si}$ is the quotient by $\bar\cG$.)
We let $V_{xx'}\sub V_x$, $U_{xx'}=U_x\cap V_{xx'}$ and $\varphi_{xx'}$, etc., be as in \eqref{6.93}. 
We let $$\rho_{xx'}=\varphi_{xx'}|_{U_{xx'}\ured}: U_{xx'}\ured\to U_{x'x}\ured.$$ Then the $\xi_{xx'}$ in \eqref{6.92} induces an isomorphism
\beq\label{b4} (\xi_{xx'})^2=(\det d\varphi_{xx'}|_{U_{xx'}\ured})^2: (\det T_{V_x})^{\otimes 2}|_{U_{xx'}\ured}
\lra \rho_{xx'}\sta  (\det T_{V_{x'}})^{\otimes 2}|_{U_{x'x}\ured}.
\eeq
On the other hand, since $\rho_{xx'}$ commutes with the projections $\pi: U_{xx'}\ured\to \fV_{si}^c$ and $\pi:U_{x'x}\ured\to\fV_{si}^c$,
we can find a holomorphic $h_{xx'}: U_{xx'}\ured\to \bar \cG$ so that $\rho_{xx'}$ is induced by applying the gauge transformation
$h_{xx'}$. At the same time, the gauge transformations $h_{xx'}$ induces an isomorphism
$$[h_{xx'}]\lsta : \ti\cL^{-1}|_{U_{xx'}\ured}\lra \rho_{xx'}\sta\ti\cL^{-1}|_{U_{x'x}\ured}.
$$

\begin{lemm} The isomorphisms \eqref{b2} fit into the commutative square
\beq\label{b4}
\begin{CD}
(\det T_{V_x})^{\otimes 2}|_{U_{xx'}\ured}
@>{(\xi_{xx'})^2}>> \rho_{xx'}\sta  (\det T_{V_{x'}})^{\otimes 2}|_{U_{x'x}\ured}\\
@VV{\eqref{b2}}V @VV{\eqref{b2}}V\\
\ti\cL^{-1}|_{U_{xx'}\ured} @>{[h_{xx'}]\lsta}>> \rho_{xx'}\sta\ti\cL^{-1}|_{U_{x'x}\ured}.
\end{CD}
\eeq
\end{lemm}

\begin{proof}
By examining the construction of $\varphi_{xx'}$, we see that it is the composition of $\varphi'_{xx'}$ with
$\varphi_{xx'}\dpri$ shown below. We let $\tilde{h}_{xx'}: V_{xx'}\to \bar\cG$ be a holomorphic extension of
$h_{xx'} : U_{xx'}\ured\to \bar\cG$. We let 
$$V_{xx'}'=\{\dbar_y\cdot \tilde{h}_{xx'}(y)\mid y\in V_{xx'}\}, \and f_{xx'}'=CS|_{V_{xx'}'}.
$$
Then $\Crit(f'_{xx'})=U_{x'x}$. We let $\varphi_{xx'}': V_{xx'}\to V_{xx'}'$ be the induced map. The second 
map $\varphi_{xx'}\dpri:V_{xx'}'\to V_{x'x}$
is a biholomorphic map, fixing $\Crit(f'_{xx'})=U_{x'x}$, such that $f_{x'}\circ \varphi_{xx'}\dpri=f_{xx'}'$.
Let $\xi'_{xx'}=\det d\varphi'_{xx'}|_{U_{xx'}\ured}$ and  $\xi''_{xx'}=\det d\varphi''_{xx'}|_{U_{xx'}\ured}$

Because the two
horizontal arrows below are induced by gauge transformations $h_{xx'}$, by the proof of \eqref{b2}, we have the following
commutative square
$$\begin{CD}
(\det T_{V_x})^{\otimes 2}|_{U_{xx'}\ured}
@>{(\xi'_{xx'})^2}>> \rho_{xx'}\sta  (\det T_{V_{x'}})^{\otimes 2}|_{U_{x'x}\ured}\\
@VV{\eqref{b2}}V @VV{\eqref{b2}}V\\
\ti\cL^{-1}|_{U_{xx'}\ured} @>{[h_{xx'}]\lsta}>> \rho_{xx'}\sta\ti\cL^{-1}|_{U_{x'x}\ured}.
\end{CD}
$$
Therefore, to prove the lemma, %it remains to show that $U_{xx'}=U_{x'x}$ (as subspaces in $\cA^{int}_{si}$) in the first place. From now on, 
we may assume $U_{xx'}=U_{x'x}$, $h_{xx'}=\id$ and $\rho_{xx'}=\id$. 
%For notational simplicity, we check that

Let us check \eqref{b4} at any $y\in U_{xx'}$. 
%(Our application only need \eqref{b4} to hold on the reduced part of $U_{xx'}$.) 
Let $S_1\sub T_{V_{xx'}}|_y$ be a subspace complementary to $T_{U_{xx'}}|_y$.
Let $S_2=\dbar_y(S_1)\sub \Omega^{0,2}(\End E)$, which via the isomorphism $R|_U\cong T_V\sta|_U$ mentioned before
\eqref{6.25} is a subspace in $T\sta_{V_{xx'}}|_y$. Let $S_1'=d\varphi_{xx'}|_y(S_1)$ and $S_2'=d\varphi_{xx'}\upmo|_y(S_2)$. Because $f_{x'}\circ \varphi_{xx'}=f_x$, we conclude that the square
$$\begin{CD}
S_1 @>{\dbar_y}>> S_2\\
@VV{d\varphi_{xx'}|_y}V @VV{d\varphi_{xx'}\upmo|_y}V\\
S_1' @>{\dbar_y}>> S_2'
\end{CD}
$$
Adding that $d\varphi_{xx'}|_y|_{T_{U_{xx'}}|_y}=\id$, similarly to the proof of Lemma \ref{1.14}, we conclude that
\eqref{b4} holds true when restricted to $y$. This proves the lemma.
\end{proof}

\begin{proof}[Proof of Theorem \ref{6.20}] We see that $\det T_{V_x}|_{U_{x}\ured}$ are local square roots of $\cL_X^{-1}|_{\pi(U_x\ured)}$ and 
the isomorphisms $\xi_{xx'}$ are local transition functions whose squares are the transition functions of $\cL_X^{-1}|_{X\ured}$.
Thus the cycle \eqref{b1} is the obstruction to the existence of a square root of $\cL_{X\ured}=\cL_X|_{X\ured}$. This proves the theorem.
\end{proof}
\bigskip

\subsection{Orientability by a tautological family}\label{S6.2.2}

The orientability of the critical virtual manifold structure on the moduli $\fX$ of simple sheaves on the Calabi-Yau 3-fold $Y$ 
now follows from the following existence result of Okounkov.
\begin{prop}\label{6.40}
Let $X$ be a reduced complex analytic space and $\cE$ be a family of perfect complexes on $Y$ parameterized by $X$.  Let $p:X\times Y\to X$ be the projection. Then there exists a holomorphic line bundle on $X$ whose square is $\cL_X=\det Rp_*R\cH om(\cE,\cE)$.
\end{prop}
The proof is due to Okounkov. We thank him for informing us of this proof.

\begin{proof}
For $i,j\in \{1,2,3\}$, let $p_{ij}$ denote the projection from $X\times X\times Y$ to the product of the $i$th and $j$th factors. Let $$\cL=\det R{p_{12}}_*R\cH om(p_{13}^*\cE,p_{23}^*\cE).$$
By Serre duality, we have an isomorphism
$$\det R{p_{12}}_*R\cH om(p_{13}^*\cE,p_{23}^*\cE)\cong \det R{p_{12}}_*R\cH om(p_{23}^*\cE,p_{13}^*\cE).
$$
Thus for the isomorphism $\tau:X\times X\to X\times X$ defined by  $\tau(x,y)=(y,x)$, we have $\tau^*\cL\cong \cL$. 

Next, by the K\"unneth formula, we have the following exact sequence and commutative diagram
\[ \xymatrix{
0 \ar[r] & \bigoplus_{i=0}^2 H^i(X)\otimes H^{2-i}(X)\ar[r]\ar[dr]_{\cup} & H^2(X\times X)\ar[d]^{\Delta^*} \ar[r] & \text{Tor}%_1(H^1(X), H^2(X))^{\oplus 2}
=0 \\ %\ar[r] &0\\
&&H^2(X)
}\]
where the cohomology groups are with integer coefficients; $\Delta^*$ is the pullback by the diagonal embedding $\Delta:X\to X\times X$; $\cup$ is the cup product, and the vanishing $\text{Tor}=\text{Tor}_1(H^1(X), H^2(X))=0$ follows from 
that $H^1(X)$ is torsion-fee, by the universal coefficient theorem.

By the definition of $\cL$, we have $$\Delta^*c_1(\cL)\cong c_1(\det R{p}_*R\cH om(\cE,\cE))$$
where $p:X\times Y\to X$ is the projection. 
Since $\tau^*c_1(\cL)=c_1(\cL)$, $c_1(\cL)$ is the sum of classes of the form $a_0\otimes a_2+ a_2\otimes a_0$ and $a_1\otimes b_1+b_1\otimes a_1$ with $a_i,b_i\in H^i(X)$. Therefore $c_1(\det R{p}_*R\cH om(\cE,\cE))=\Delta^* c_1(\cL)$ is the sum of classes of the form $a_0\wedge a_2+a_2\wedge a_0=2 a_0\wedge a_2$ and $a_1\wedge b_1+b_1\wedge a_1=0$. Hence $c_1(\det R{p}_*R\cH om(\cE,\cE))\in H^2(X,\ZZ)$ is divisible by $2$. 
%which implies that there is a topologic line bundle whose tensor square is isomorphic to the underlying topological line bundle of $\det Ext\ud_{\pi}(\cE,\cE)$. 
This implies that there is a holomorphic line bundle on $X$ whose (tensor) square is isomorphic to $\det R{p}_*R\cH om(\cE,\cE)$.
\end{proof}

From Theorems \ref{4.45}, \ref{6.20} and Proposition \ref{6.40}, we obtain the following.
\begin{theo}\label{6.41}
Let $X$ be a bounded open analytic subspace of $\SH^c$ whose reduced space $X\ured$ is equipped with a tautological family. 
Then the critical virtual manifold structure from Theorem \ref{6.8} is orientable and thus $X$ admits a perverse sheaf $P$ and polarizable mixed Hodge module $\cM$ with $\mathrm{rat}(\cM)=P$, that are geometric gluings of local perverse sheaves and mixed Hodge modules of vanishing cycles.
\end{theo}
\begin{proof}
By the Seidel-Thomas twists, $X$ is isomorphic to an open $\tilde{X}\sub\fV_{si}^{\tilde{c}}$ and $\tilde{X}\ured$ admits a tautological family because $X\ured$ does. 
By Proposition \ref{6.40}, $\cL_{\tilde{X}\ured}$ admits a square root and hence the critical virtual manifold structure on $\tilde{X}$ is orientable by Theorem \ref{6.20}. Since the critical virtual manifold structure on $X$
was imported from that on $\tilde{X}$, $X$ is an orientable critical virtual manifold.
\end{proof}

As we mentioned before, if $X$ is (quasi-)compact or a quasi-projective scheme or a Noetherian scheme, then $X$ is bounded.

\medskip

For the Gopakumar-Vafa theory in \S\ref{ch7}, we will use the following.
\begin{theo}\label{6.55} Let $\cO_Y(1)$ be an ample line bundle on $Y$.
Let $\fX$ be the moduli space of stable sheaves $\cF$ on $Y$ whose Hilbert polynomial is $\chi(\cF\otimes \cO_Y(m))=dm+e$ for coprime integers $d,e$. Then $X$ is an orientable critical virtual manifold and hence there is a perverse sheaf $P$ (resp. polarizable mixed Hodge module $\cM$) on $X$ with $\mathrm{rat}(\cM)=P$, which is the gluing of local perverse sheaves (resp. mixed Hodge modules) of vanishing cycles.
\end{theo}
\begin{proof}
Since the coefficients $d$ and $e$ of the Hilbert polynomial are coprime, by the standard descent argument (\cite[Theorem 6.11]{Mar}), the universal family of Grothendieck's Quot scheme descends to 
$X$. By Theorem \ref{6.41}, we have the gluing of perverse sheaves and mixed Hodge modules.
\end{proof}

\chapter{Gopakumar-Vafa invariant via perverse sheaves}\label{ch7}
\def\fS{\mathfrak{S} }
\black

We proved that a bounded open moduli space $X\sub \SH^c$ of simple sheaves on a Calabi-Yau 3-fold
admits a critical virtual manifold structure if there is a tautological family. 
In this section, we provide a mathematical theory of the Gopakumar-Vafa invariant as an application.

\section{Intersection cohomology sheaf}\label{S7.1}
As before, let $Y$ be a smooth projective Calabi-Yau 3-fold over $\CC$. From the string theory (\cite{GoVa, Katz}), it is expected that \begin{enumerate}\item there are integers $n_h(\beta)$, called the Gopakumar-Vafa invariants (GV invariants, for short) which contain all the information about the Gromov-Witten invariants $N_g(\beta)$ of $Y$ in the sense that 
\beq\label{7.1}\sum_{g,\beta}N_g(\beta)q^\beta\lambda^{2g-2}=\sum_{k,h,\beta}n_h(\beta)\frac1{k} \left(2\sin (\frac{k\lambda}2)\right)^{2h-2}q^{k\beta} \eeq
where $\beta\in H_2(Y,\ZZ)$, $q^\beta=\exp(-2\pi\int_\beta c_1(\cO_Y(1)))$;
%\item the GV invariants are integer valued;
\item  $n_h(\beta)$ come from an $sl_2\times sl_2$ action on some cohomology theory of the moduli space $\fX$ of one dimensional stable sheaves on $Y$;
\item $n_0(\beta)$ should be the Donaldson-Thomas invariant of the moduli space $\fX$.
\end{enumerate}
By using the global perverse sheaf $P$ constructed above and the method of \cite{HST}, we can give a geometric theory for the GV invariants.

\medskip

We recall the following facts from \cite{Sai88, Sai90}. 
\begin{theo}\label{7.2} 
Let $f:X\to Y$ be a projective morphism with $\omega$ the first Chern class of a relative ample line bundle;
let $P$ be a perverse sheaf on $X$ underlying a polarizable Hodge module $\cM$, then
\begin{enumerate}
\item \emph{(Hard Lefschetz theorem)}  the cup product induces an isomorphism
\[ \omega^k:\, ^p\!\cH^{-k}Rf_*P\lra \,^p\!\cH^{k}Rf_*P ;\]

\item \emph{(Decomposition theorem)} we have the decomposition 
\[ Rf_*P\cong \oplus_k  \,^p\!\cH^{k}Rf_*P[-k] \]
and each summand $^p\cH^{k}Rf_*P[-k]$ is a perverse sheaf underlying a polarizable Hodge module.
\end{enumerate}
\end{theo}

We fix an algebraic curve class $\beta\in H_2(Y,\ZZ)$; let $\beta\dual\in H^4(Y,\ZZ)$ be its Poincare dual.
We form the moduli space $\fX=\cM_Y(-\beta\dual,\chi)$ of stable one-dimensional sheaves $F$ on $Y$ with $\chi(F)=\chi$ and $c_2(E)=-\beta\dual\in H^4(Y,\ZZ)$. (Consequently, $\rank F=c_1(F)=0$.) We assume 
$c_1(\sO_Y(1))\cdot \beta$ is coprime to $\chi$, so that
$\cM_Y(-\beta\dual,\chi)$ is projective.

%For instance, we can fix the Hilbert polynomial to be $P(m)=dm+\chi$ with $d$ and $\chi$ coprime. 

Let $\ti X$ be the semi-normalization of $X$. Applying \cite{Koll}, we have the Hilbert-Chow morphism
$\ti X\to Chow(Y)$; let $S$ be its image. By \cite{Koll}, the morphism $\ti X\to X$ is one-to-one and hence a homeomorphism because $\ti X$ is projective and $X$ is separated. The induced
morphism $f:\ti X\to S$ is projective and the intersection cohomology sheaf $IC\ud=IC_{\ti X}(\CC)\ud$ underlies a  simple polarizable mixed Hodge module. In \cite{HST}, Hosono-Saito-Takahashi showed that the hard Lefschetz theorem applied  to $f$ and $c:S\to \mathrm{pt}$ gives us an $sl_2\times sl_2$ action on the 
intersection cohomology $IH^*(\ti X)=\bbH^*(\ti X,IC\ud)$ as follows: The relative Lefschetz isomorphism $$\, ^p\!\cH^{-k}Rf_*(IC\ud)\lra \,^p\!\cH^{k}Rf_*(IC\ud)$$ for $f$ gives an action of $sl_2$, called the left action, via the isomorphisms
$$\bbH^*(\ti X, IC\ud)\cong \bbH^*(S,Rf_*IC\ud)\cong \oplus_k \bbH^*(S,  \,^p\!\cH^{k}Rf_*(IC\ud)[-k])$$
from the decomposition theorem.
On the other hand, since $\,^p\!\cH^{k}Rf_*(IC\ud)[-k]$ underlies a polarizable Hodge module, 
$$\bbH^*(S,  \,^p\!\cH^{k}Rf_*(IC\ud)[-k])$$ is equipped with another action of $sl_2$, called the right action, by hard Lefschetz again. Therefore we obtain an action of $sl_2\times sl_2$ on the intersection cohomology $IH^*(\ti X)=IH^*(X)$. 

If $C\in S$ is  a genus $h$ smooth curve in the class $\beta$, the fiber of $f$ over $C$ is the Jacobian of line bundles on $C$ whose cohomology is an $sl_2$-representation space $$\left((\frac12)\oplus 2(0)\right)^{\otimes h},
$$
where $(\frac12)$ denotes the 2-dimensional representation of $sl_2$ while $(0)$ is the trivial 1-dimensional representation. 
In \cite{HST}, the authors propose a theory of the Gopakumar-Vafa invariants  by using the $sl_2\times sl_2$ action on $IH^*(\ti X,\CC)$ as follows: By the Clebsch-Gordan rule, it is easy to see that one can uniquely write the $sl_2\times sl_2$-representation space $IH^*(\ti X,\CC)$ in the form 
\[ IH^*(\ti X,\CC)=\bigoplus_h \left((\frac12)_L\oplus 2(0)_L\right)^{\otimes h}\otimes R_h, \]
where $(\frac{k}2)_L$ denotes the $k+1$ dimensional irreducible representation of the left $sl_2$ action while $R_h$ is a representation space of the right $sl_2$ action. The authors of \cite{HST} define the GV invariant as the Euler number $Tr_{R_h}(-1)^{H_R}$ of $R_h$ where $H_R$ is the diagonal matrix in $sl_2$ with entries $1,-1$.

However it seems unlikely that the invariant $n_h(\beta)$ defined using the intersection cohomology as in \cite{HST} will relate to the GW invariants of a general Calabi-Yau 3-fold
$Y$ as proposed by Gopakumar-Vafa. %because intersection cohomology is unstable under deformation. 
We propose to use the perverse sheaf $P$ on $X$ constructed above instead of $IC\ud$.

\section{GV invariants from perverse sheaves} \label{S7.2}
We continue to denote by $\fX$ the moduli of stable sheaves on $Y$ of $\rank E=c_1(E)=0$, $c_2(E)=-\beta\dual$ and $\chi(E)=\chi$. By Theorem \ref{6.55},  when 
$\beta\cdot c_1(\cO_Y(1))$ and $\chi$ are coprime,  
there is a universal family $\cE$ on $\fX\times Y$ and thus we have a perverse sheaf $P$ that underlies a polarizable mixed Hodge module $\cM$.

%In this subsection, we assume that $\det \Ext_{\pi}\bul(\cE,\cE)$ admits a square root so that we have a perverse sheaf $P\ud$ and a MHM $M\ud$ which are locally the perverse sheaf and MHM of vanishing cycles for a restricted CS functional. 
%\begin{rema} In \cite{Hua}, it is proved that if the Calabi-Yau 3-fold $Y$ is simply connected and $H^*(Y,\ZZ)$ is torsion-free, then $\det \Ext_{\pi}\bul(\cE,\cE)$ admits a square root. For instance, when $Y$ is a quintic 3-fold, we have the desired perverse sheaf and MHM.
%\end{rema}
Since the semi-normalization $\gamma:\ti X\to X$ is a homeomorphism, the pullback $\ti P$ of $P$ is a perverse sheaf
such that $\gamma_*\ti P\cong P$, and 
%. By Theorem \ref{thmMHM}, $P\ud$ lifts to a MHM $M\ud$ and its 
the pullback $\ti \cM$ satisfies $\rat(\ti \cM)=\ti P$ since $\rat$ preserves Grothendieck's six functors (\cite{Sai90}). 
Let $\hat{\cM}=\mathrm{gr}^W \ti \cM$ be the graded object of $\ti \cM$ with respect to the weight filtration $W$. Then $\hat{\cM}$ is a direct sum of polarizable Hodge modules (\cite{SaiICM}). Let $\hat{P}=\rat(\hat{\cM})$; because $\rat$ is an exact functor (\cite{Sai90}), it is the gradation $\mathrm{gr}^W \ti P$ by the weight filtration of $\ti P$. 

%Since $Perv(\ti X)$ is both Artinian and Noetherian, there is a composition series
%$0=P_0\ud\subset P_1\ud\subset\cdots \subset P_n\ud=\ti P\ud$ and we can find a unique semisimple perverse sheaf $$\mathrm{gr} \ti P\ud =\bigoplus_i Q\ud_i,\quad Q\ud_i=P\ud_i/P\ud_{i-1}.$$
By \cite[\S5]{Sai88}, the hard Lefschetz theorem and the decomposition theorem hold for the polarizable Hodge module $\hat{\cM}$. Hence by applying the functor $\rat$, we obtain the hard Lefschetz theorem and the decomposition theorem for $\hat{P}$.
%For each simple perverse sheaf $Q\ud_i$
Therefore, we can apply the argument in \S\ref{S7.1} to obtain an action of $sl_2\times sl_2$  on the hypercohomology $\bbH^*(\ti X,\hat{P})$ to write $$\bbH^*(\ti X,\hat{P})
\cong \bigoplus_h \left((\frac12)_L\oplus 2(0)_L\right)^{\otimes h}\otimes R_{h}.$$ 

\begin{defi}\label{7.3} We define the Gopakumar-Vafa invariant as %the Euler number of $R_{h}$ 
\[ n_h(\beta):= Tr_{R_{h}}(-1)^{H_R}.\]
\end{defi}

The GV invariant $n_h(\beta)$ is integer valued and defined by an $sl_2\times sl_2$ representation space $\bbH^*(\ti X, \hat{P})$ as expected from \cite{GoVa}.
\begin{prop}\label{7.4} 
The number $n_0(\beta)$ is the Donaldson-Thomas invariant of $\fX$.
\end{prop}
\begin{proof} 
Recall that the DT invariant is the Euler number of $X$ weighted by the Behrend function $\nu_{\fX}$ on $X$,
and that $\nu_{\fX}(x)$ for $x\in X$ is the Euler number of the stalk cohomology of $P$ at $x$. Therefore the DT invariant of $\fX$ is the Euler number of $\bbH^*(X,P)$.

Since the semi-normalization $\gamma:\ti X\to X$ is a homeomorphism, $\gamma_*\ti P\cong P$ and $\bbH^*(X,P)\cong \bbH^*(X,\gamma_*\ti P)\cong \bbH^*(\ti X, \ti P)$. Therefore
$$ DT(\fX)=\sum_k (-1)^k\dim\bbH^k(X,P)=\sum_k (-1)^k\dim  \bbH^k(\ti X,\ti P).$$
Since $\ti P$ has a filtration $W$ with $\hat{P}=\mathrm{gr}^W \ti P$, we have the equality of alternating sums
$$\sum_k (-1)^k\dim  \bbH^k(\ti X,\ti P)=\sum_{k} (-1)^k\dim \bbH^k(\ti X,\hat{P}).$$
Since the Euler number of the torus part $\left((\frac12)_L\oplus 2(0)_L\right)^{\otimes h}$ is zero for $h\ne 0$, 
$$\sum_{k} (-1)^k\dim \bbH^k(\ti X,\hat{P})= Tr_{R_{0}}(-1)^{H_R}=n_0(\beta).$$
This proves the proposition. \end{proof}

We propose the following conjecture.
\begin{conj}\label{7.5} Let $n_h(\beta)$ be defined using the moduli of 1-dimensional stable sheaves with $c_2(E)=-\beta\dual$ and $\chi(E)=\chi$ such that $\beta\cdot c_1(\cO_Y(1))$ and $\chi$ are coprime. Then
\begin{enumerate}
\item $n_h(\beta)$ are invariant under deformation of the complex structure of $Y$.
\item The GV invariants $n_h(\beta)$ depend only on $\beta$ and are independent of the $\chi$.
\item $n_h(\beta)$ are independent of the choice of polarizations of $Y$.
\item The identity \eqref{7.1} holds.
\end{enumerate}
\end{conj}

Note that for $h=0$, (1) follows from Proposition \ref{7.4} and \cite{Tho}. Also by \cite{JoSo} and \cite{Tho}, (3) is known for $h=0$. Of course, (1)-(3) are consequences of (4). 
Furthermore, establishing the identity \eqref{7.1} will equate Definition \ref{7.3} with that introduced by Pandharipande-Thomas \cite{PT} for a large class of
Calabi-Yau 3-folds (cf. \cite{PP}).

\section{K3-fibered Calabi-Yau 3-folds}

In this last section, we show that Conjecture \ref{7.5} holds for a primitive fiber class of K3 fibered Calabi-Yau 3-folds.

Let $Y\to\PP^1$ be a K3 fibered smooth projective Calabi-Yau 3-fold, let $0\in \PP^1$ be a closed point
and let $\iota_0: Y_0\sub Y$ be the fiber over $0$, assumed to be smooth. Let
$\beta_0\in H_2(Y_0,\ZZ)$ be an algebraic curve class so that its Poincare dual $\beta_0\dual\in H^2(Y_0,\ZZ)$
ceases to be $(1,1)$ type in the first order deformation of $Y_0$ in the family $Y_c$, $c\in \PP^1$.\footnote{Namely, if we let $\Delta\sub \PP^1$ be a disk-like analytic neighborhood containing $0\in\PP^1$,
%Here we say $\beta$ cease to be $(1,1)$ in the first order deformation of 
%$Y_0$ in $Y$ if the following holds. L
let $\ti\beta\in \Gamma(\Delta, R^2p_* \ZZ_Y)$ be the continuous extension of $\beta_0$,
and let $\ti\omega\in\Gamma(\Delta, p_*\Omega_{Y/\Delta}^2)$ be a nowhere vanishing holomorphic section of 
relative $(2,0)$-form, then $p_*(\ti\omega\wedge\ti\beta)$ is a holomorphic function on $\Delta$ that has order one vanishing at $0$.}

We let $\beta=\iota_{0\ast}\beta_0\in H_2(Y,\ZZ)$. %, and let $\beta\dual\in H^4(Y,\ZZ)$ be its Poincare dual.
Fixing an ample $\cO_Y(1)$ on $Y$, we form the moduli $\cM_Y(-\beta\dual,\chi)$ of stable sheaves on $Y$ as before,
and form the moduli $\cM_{Y_0}(\beta_0\dual,\chi)$
of one dimensional $H|_{Y_0}$-stable sheaves $F$ of $\cO_{Y_0}$-modules with
$c_2(F)=-\beta\dual$ and $\chi(F)=\chi$. We assume $\beta\cdot c_1(\cO_Y(1))$ and $\chi$ are coprime. Then
both $\cM_Y(-\beta\dual,\chi)$ and $\cM_{Y_0}(\beta\dual_0,\chi)$ are projective schemes.
Because the polarization of $Y$ restricts to $H|_{Y_0}$, we have a natural closed embedding
\beq\label{7.6}
\cM_{Y_0}(\beta_0\dual,\chi)\mapright{\sub} \cM_Y(-\beta\dual,\chi).
\eeq

\begin{lemm}
Let the notation be as before, and suppose $\beta_0$ ceases to be $(1,1)$ in the first order deformation of 
$Y_0$ in $Y$. %and there are no $c\ne 0\in\Delta$ such that $\iota_c\sta \ti\beta\in H^{1,1}(Y_c,\RR)$. 
Then the embedding \eqref{7.6} is an open and closed embedding of schemes.
\end{lemm}

\begin{proof}
We pick a disk-like analytic neighborhood $0\in \Delta\sub \PP^1$ so that
for any $c\ne 0\in \Delta$, $Y_c=Y\times_{\PP^1}c$ is smooth and, letting $\iota_c: Y_c\to Y$ be the inclusion,
$\iota_c^!\beta\dual\in H^2(Y_c,\ZZ)$ is not a $(1,1)$-class. By our assumption that $\beta_0$ cease to be
$(1,1)$ in the first order, such $\Delta$ exists.

We claim that for any sheaf $[F]\in \cM_Y(-\beta\dual,\chi)$, $F=\iota_{c\ast}F'$ for a 
sheaf of $\cO_{Y_c}$-modules. Indeed, let $\text{supp}(F)$ be the scheme-theoretic support
of $F$. Since $F$ is stable, $\text{spt}(F)$ is connected and proper; because $\beta\in H_2(Y,\ZZ)$ is a fiber class,
the underlying set of $\text{spt}(F)$ is contained in 
a closed fiber $Y_c\sub Y$ for a $c\in\PP^1$. 

We let $u$ be a uniformizing parameter of $c\in \PP^1$. 
Since $F$ is coherent, there is a
positive integer $k$ so that $\text{spt}(F)\sub (u^k=0)$. 
%(Here we identify $t-c$ with its pullback in $\Gamma(\cO_Y)$.)
In particular, $F$ is annihilated
by $u^k$. Since $u|_{\text{spt}(F)}$ is nilpotent, multiplying by $u$ defines a sheaf homomorphism
$\cdot u: F\to F$, which has non-trivial kernel since $u^k$ annihilates $F$.
Since $F$ is stable, this is possible only if $F$ is annihilated by $u$. 
Therefore, letting $F'=F/u\cdot F$, which is a sheaf of $\cO_{Y_c}$-modules,
we have $F=\iota_{c\ast}F'$. Repeating the same argument for families in $\cM_Y(-\beta\dual,\chi)$,
we conclude that the assignment $F\mapsto c$ extends to a morphism $p: \cM_Y(-\beta\dual,\chi)\to \PP^1$.

For the $\Delta$ chosen, we claim that either $c\not\in\Delta$ or $c=0$. If not, then $c_1(F')=\iota_c^!\beta$ will be in $H^{1,1}(Y_c,\RR)$, and this is possible only when either $c=0$ or $c\not\in\Delta$. This prove the claim. Therefore, 
$p\upmo(\Delta-0)=\emptyset$. In particular, 
$$\hat{\cM}_{Y_0}(\beta_0\dual,\chi)\defeq p\upmo(\Delta)
$$
is an open and closed subscheme of $\cM_Y(-\beta\dual,\chi)$. Clearly, \eqref{7.6} induces a
closed embedding
\beq\label{7.7}
\cM_{Y_0}(\beta_0\dual,\chi)\mapright{\sub}\hat{\cM}_{Y_0}(\beta_0\dual,\chi).
\eeq
We prove that it is an isomorphism. Indeed, by the previous argument, we know that
\eqref{7.7} is a homeomorphism. To prove that it is an isomorphism, we need to show that for any
local Artin ring $A$ with quotient field $\CC$ and morphism $\varphi_A: \spec A\to \cM_Y(-\beta\dual, \chi)$
such that it sends the closed point in $\spec A$ into $M_{Y_0}(\beta_0\dual,\chi)$, then 
$\varphi_A$ factors through $\cM_{Y_0}(\beta_0\dual,\chi)$. 
By an induction on the length of $A$, we only need to consider the case where
there is an ideal $I\sub A$ such that $\dim_\CC I=1$ and the restriction $\varphi_{A/I}:\spec A/I\to \cM_Y(-\beta\dual,\chi)$ 
already factors through $\cM_{Y_0}(\beta_0\dual,\chi)$. 

Let $F$ be the sheaf of $A\times\cO_Y$-modules that is the pullback of the universal family of $\cM_Y(-\beta\dual,\chi)$
via $\varphi_A$. Let $t$ be a uniformizing parameter of $0\in \PP^1$. 
As $\varphi_{A/I}$ factors through $\cM_{Y_0}(\beta_0\dual,\chi)$, $t\cdot F\sub I\cdot F$. 
If $t\cdot F=0$, then $\varphi_A$ factors, and we are done. Suppose not. Then since $F_0=F\otimes_A\CC$ is
stable, there is a $c\in I$ so that $t\cdot F=c\cdot F\sub F$. Thus $(t-c)\cdot F=0$.
We now let $\psi: \spec A\to \Delta$ be the morphism defined by
$\psi\sta(t)=c$, and let $Y_{A}=Y\times_{\Delta, \psi} \spec A$.
Then $Y_A\to \spec A$ is a family of K3 surfaces with a tautological embedding $\iota_A: Y_A\to Y\times \spec A$.
Then $(t-c)\cdot F=0$ means that there is an $A$-flat family of sheaves $F'$ of $\cO_{Y_A}$-modules so
that $\iota_{A\ast}F'=F$.

Let $q$ be the projection and $\iota$ be the tautological morphism fitting into the Cartesian square
$$\begin{CD}
Y_A @>{\iota}>> Y_\Delta\defeq Y\times_{\PP^1}\Delta\\
@VV{q}V @VV{p}V\\
\spec A @>{\psi}>> \Delta
\end{CD}
$$
Let $\ti\beta\in \Gamma(\Delta, R^2p_* \ZZ_Y)$ be the continuous extension of $\beta_0\dual$ and
let $\ti\omega\in\Gamma(\Delta, p_*\Omega_{Y/\Delta}^2)$ be a nowhere vanishing holomorphic section of 
relative $(2,0)$-form. Notice that by our assumption, $p_*(\ti\omega\wedge\ti\beta)$ is not divisible by $t^2$.
On the other hand, by our construction, $c_1(F')=\iota\sta\ti\beta$. Thus for the $(2,0)$-form $\ti\omega$ mentioned, we have
$$0=q_\ast\bl c_1(F')\wedge \iota\sta\ti\omega\br=q_\ast\iota\sta(\ti\beta\wedge\ti\omega)
=\psi\sta p_\ast(\ti\beta\wedge\ti\omega).
$$
Since $p_\ast(\ti\beta\wedge\ti\omega)$ is not divisible by $t^2$, $\psi$ must factor through $0\in\Delta$. This proves that $c=0$ and $\varphi_A$ factors
through $\cM_{Y_0}(\beta_0\dual,\chi)$. This proves the proposition.
\end{proof}

\begin{prop}\label{7.10} Let the notation be as stated, and suppose $\beta_0\dual\in H^2(Y_0,\ZZ)$ 
ceases to be $(1,1)$ in the first order deformation of $Y_0$ in $Y\to \PP^1$.
Let $\chi\in\ZZ$ be coprime to $\beta\cdot c_1(\cO_Y(1))$ and let $\beta=\iota_{0\ast}\beta_0$. Then 
$$\fX_0\defeq \cM_{Y_0}(\beta_0\dual,\chi)\sub \fX\defeq \cM_Y(-\beta\dual, \chi)
$$ 
is a smooth open and closed (projective) subscheme. Suppose further that $Y_0$ is an elliptic K3 surface with
a section. Then \eqref{7.1} holds for the GV invariants of the perverse sheaf $P\bul$ on $\fX_0$ where 
$N_g(\beta_0)$ in \eqref{7.1} are the GW invariants contributed from the connected components of stable maps to $Y$ that factor through $Y_0\sub Y$.
\end{prop}

\begin{proof}
By the argument prior to the statement of the proposition, we know that $\fX_0\sub\fX$ is open and closed. 
Since $Y_0$ is a smooth K3 surface and $\chi$ is coprime to $\beta\cdot c_1(\cO_Y(1))$, 
$\cM_{Y_0}(\beta_0\dual,\chi)$ is a smooth
simply connected projective variety of dimension $\dim X=\beta_0^2+2$ (cf. \cite{Mukai}). 

By assumption, $\fX$ has a universal family. We let $P$ be a perverse sheaf on $X$ provided by Theorem \ref{6.55}. Since $\fX_0$ is smooth, connected and simply connected, and is open and closed in $\fX$, $X_0$,
and the restriction
$P|_{X_0}=\QQ_{X_0}$; thus $\bbH^*(X_0,P)=H^*(X_0,\QQ)$. 

Let 
$$Y_0=S\mapright{\pi} \PP^1
$$ 
be an elliptic K3 surface and $\beta_0=C+kF$ where $C$ is a section and $F$ is a fiber. We can calculate the GV invariants in this case by the same calculation as in \cite[Theorem 4.7]{HST}. 
Since the details are obvious modifications of those in \cite[\S4]{HST}, we briefly outline the calculation. Indeed by Fourier-Mukai transform, $X$ is isomorphic to the Hilbert scheme $S^{[k]}$ of points on $S$ and the Chow scheme in this case is a complete linear system $\PP^k$. The Hilbert-Chow morphism is 
\[ S^{[k]}\mapright{\mathfrak{hc}} S^{(k)} \mapright{\pi} (\PP^1)^{(k)}\cong \PP^k.\]
It is easy to see that the cohomology $H^*(S,\QQ)$ of $S$ as an $sl_2\times sl_2$ representation space by (relative) hard Lefschetz applied to $S\to \PP^1$ is 
\[ (\frac12)_L\otimes (\frac12)_R + 20\cdot (0)_L\otimes (0)_R.\]
If we denote by $t_L$ (resp. $t_R$) the weight of the action of the maximal torus for the left (resp. right) $sl_2$ action, we can write $H^*(S,\QQ)$ as $(t_L+t_L^{-1})(t_R+t_R^{-1})+20$. Hence $H^*(S^{(k)},\QQ)$ is the invariant part of 
\[ \left( (\frac12)_L\otimes (\frac12)_R + 20\cdot (0)_L\otimes (0)_R\right)^k
\]
by the symmetric group action. In terms of Poincar\'e series, we can write
\[ \sum_k P_{t_L,t_R}(S^{(k)})q^k=\frac{1}{(1-t_Lt_Rq)(1-t_L^{-1}t_Rq)(1-t_Lt_R^{-1}q)(1-t_L^{-1}t_R^{-1}q)(1-q)^{20}}. \]
Applying the decomposition theorem (\cite{BBD}) for the semismall map $S^{[k]}\to S^{(k)}$,  we find that $\sum_k P_{t_L,t_R}(S^{[k]})q^k$ is
$$\prod_{m\ge 1}\frac{1}{(1-t_Lt_Rq^m)(1-t_L^{-1}t_Rq^m)(1-t_Lt_R^{-1}q^m)(1-t_L^{-1}t_R^{-1}q^m)(1-q^m)^{20}}$$
which gives
\beq\label{7.8}\sum_k P_{t_L,t_R}(S^{[k]})|_{t_R=-1}q^k=\prod_{m\ge 1}\frac{1}{(1+t_Lq^m)^2(1+t_L^{-1}q^m)^2(1-q^m)^{20}}. \eeq
By definition, the GV invariants are defined by writing \eqref{7.8} as
\beq\label{7.9} \sum_{h,k} q^k(t_L+t_L^{-1}+2)^h\otimes R_h(S^{[k]})|_{t_R=-1} =\sum_{h,k}q^kn_h(k)(t_L+t_L^{-1}+2)^h \eeq
By equating \eqref{7.8} and \eqref{7.9} with $t_L=-y$, we obtain
$$%\beq\label{12173}
\sum_{h,k} (-1)^hn_h(k)(y^{\frac12}-y^{-\frac12})^{2h}q^{k-1}=\frac{1}{q\prod_{m\ge 1}(1-yq^m)^2(1-y^{-1}q^m)^2(1-q^m)^{20}}
$$
with $\beta_0^2=2k-2$. On the other hand, by \cite[Theorem 1]{MP}, we have
$$%\beq\label{12174}
\sum_{h,k} (-1)^hr_h(k)(y^{\frac12}-y^{-\frac12})^{2h}q^{k-1}=\frac{1}{q\prod_{m\ge 1}(1-yq^m)^2(1-y^{-1}q^m)^2(1-q^m)^{20}}
$$
where $r_h(k)$ are the BPS invariants from the Gromov-Witten theory for $Y_\Delta \to \Delta$ (of multiple class $k\beta_0$). 
Combining these two identities, we find that $$n_h(k)=r_h(k),
$$
which verifies Conjecture \ref{7.5} on invariants localized near $Y_0$ for the Calabi-Yau 3-fold $Y\to \PP^1$ . 
\end{proof}

%It will be interesting to extend this constructions to the setting of stable pairs. Then it may be possible to extend the theory of Gopakumar-Vafa invariants to the moduli scheme of stable pairs. Let $M$ be the moduli scheme of stable pairs $(F,s)$ with fixed topological type, where $F$ is a pure sheaf of one dimensional support and $s\in H^0(F)$ which is $\alpha$-stable for some $\alpha>0$. We consider the morphism $M\to  S$ to the Chow scheme possibly after taking the semi-normalization, sending $(F,s)$ to the support of $F$. The general fiber over a curve $C$ of genus $g$ is expected to be the symmetric product $C^{(d)}$ for a suitable $d$ defined by the topological type. Let $J_g=[\{(\frac12)_L+2g\,(0)_L\}^d]^{S_d}$ be the cohomology of $C^{(d)}$ as an $sl_2$ representation space. We can write the perverse hypercohomology of $M$ in the form $\bigoplus J_g\otimes R_g$ and define the GV invariant as the Euler number of $R_g$. 

\backmatter

%-----------------------------------------------------------------------------
% Beginning of biblio.tex
%-----------------------------------------------------------------------------

\bibliographystyle{amsalpha}

%-----------------------------------------------------------------------------
% End of biblio.tex
%-----------------------------------------------------------------------------

%\include{ind}
\end{document}